%% file: vi_ideal_mhd2d.tex
\documentclass[12pt,a4paper,reqno]{article}

\usepackage{amsmath}
\usepackage{amsfonts}
\usepackage{amssymb}
\usepackage{amsthm}
\usepackage{geometry}
\usepackage{graphicx}
\usepackage[unicode,bookmarks,breaklinks=true]{hyperref}
\usepackage[l2tabu, orthodox]{nag}
\usepackage{mathbbol}
\usepackage[numbers,square]{natbib}
\usepackage{pgf}
\usepackage{pgfplots}
\usepackage[subrefformat=parens,labelformat=parens]{subfig}
\usepackage{tikz}

\usetikzlibrary{arrows}
\usetikzlibrary{matrix}

\newcommand*{\doi}[1]{\href{http://dx.doi.org/\detokenize{#1}}{doi: \detokenize{#1}}}

\input{definitions_math}

\begin{document}

\title{Variational Integrators for\\Ideal Magnetohydrodynamics}

\author{
\large{Michael Kraus}\\
\small{(michael.kraus@ipp.mpg.de)}
\vspace{.5em}\\
\normalsize{Max-Planck-Institut f\"ur Plasmaphysik}\\
\normalsize{Boltzmannstra\ss{}e 2, 85748 Garching, Deutschland}%
\vspace{.5em}\\
\normalsize{Technische Universit\"at M\"unchen, Zentrum Mathematik}\\
\normalsize{Boltzmannstra\ss{}e 3, 85748 Garching, Deutschland}%
\vspace{1em}\\
\and
\large{Omar Maj}\\
\small{(omar.maj@ipp.mpg.de)}
\vspace{.5em}\\
\normalsize{Max-Planck-Institut f\"ur Plasmaphysik}\\
\normalsize{Boltzmannstra\ss{}e 2, 85748 Garching, Deutschland}%
\vspace{.5em}\\
\normalsize{Technische Universit\"at M\"unchen, Zentrum Mathematik}\\
\normalsize{Boltzmannstra\ss{}e 3, 85748 Garching, Deutschland}%
\vspace{1em}\\
}

\date{\today}

\maketitle

\begin{abstract}
A variational integrator for ideal magnetohydrodynamics is derived by applying a discrete action principle to a formal Lagrangian.
Discrete exterior calculus is used for the discretisation of the field variables in order to preserve their geometrical character.
The resulting numerical method is free of numerical resistivity, thus the magnetic field line topology is preserved and unphysical reconnection is absent.
In 2D numerical examples we find that important conservation laws like total energy, magnetic helicity and cross helicity are satisfied within machine accuracy.
\end{abstract}


\newpage

\tableofcontents

\newpage

\section{Introduction}

Magnetohydrodynamics (MHD) describes the dynamics of electrically conducting fluids like plasmas or liquid metals.
It is one of the most widely applied theories in laboratory as well as
astrophysical plasma physics
\cite{Schnack:2009,GoedbloedPoedts:2004,Biskamp:2003,Davidson:2001,Freidberg:1987},
used to describe, e.g., macroscopic phenomena like equilibrium states in tokamaks or stellarators, large scale turbulence in space plasmas, and dynamos that generate magnetic fields of stars and planets.
The structure of the equations is very similar to hydrodynamics, albeit in MHD the fluid equations are coupled with Faraday's law and Amp\`{e}re's law from electrodynamics, thereby allowing for an even richer variety of phenomena.

When numerically solving the MHD equations, it is important to preserve certain properties of the equations in order to obtain physically accurate and reliable solutions.
One such property is of topological nature, namely that in absence of electric resistivity the magnetic field line topology is preserved and magnetic field lines cannot open up and reconnect (frozen-in condition).
Another property is that the magnetic field is divergence-free, which implies the absence of magnetic monopoles.
Moreover, under suitable conditions the system satisfies several conservation laws, namely for energy, magnetic helicity and cross helicity.

The structure-preserving integration of ideal magnetohydrodynamics has attracted the interest of several researchers.
\citeauthor{LiuWang:2001} \cite{LiuWang:2001} approached the problem by coupling the MAC scheme \cite{Harlow:1965} for the Navier-Stokes equation with Yee's scheme \cite{Yee:1966} for the Maxwell equations.
Recently, more geometric motivated approaches were presented. \citeauthor{Gawlik:2011} \cite{Gawlik:2011} used a discrete Euler-Poincar\'{e} principle, which yields a similar scheme as that of \citeauthor{LiuWang:2001}, but with different time discretisation.
A variational integrator in Lagrangian variables, based on directly discretising Newcomb's Lagrangian \cite{Newcomb:1962}, has been derived by \citeauthor{Zhou:2014} \cite{Zhou:2014}.
Here, we propose a variational discretisation in Eulerian variables
\cite{MarsdenPatrick:1998}, based on a formal Lagrangian formulation
\cite{KrausMaj:2015} combined with ideas from discrete exterior
calculus~\cite{Robidoux:2011, Desbrun:2008, Hirani:2003}. While the
discretisation of the variational formulation leads to exact conservation of
energy, magnetic helicity and cross helicity, preserving the differential form
character of the physical variables ensures that the divergence of the magnetic
field is preserved and in combination with a staggered grid, prevents checker-boarding, a spurious phenomenon often observed with finite difference discretisations of incompressible fluid equations.

The outline of the paper is as follows.
In Section~\ref{sec:ideal_mhd}, the ideal MHD equations are reviewed. We manipulate the equations into a form suitable for our means and sketch the construction of formal Lagrangians as a starting point for deriving variational integrators.
In Section~\ref{sec:dec}, we provide the building blocks of discrete exterior calculus on two-dimensional cartesian meshes.
In Section~\ref{sec:vi}, we describe the variational discretisation on the staggered grid, which is motivated by the discrete differential forms of Section~\ref{sec:dec}.
In Section~\ref{sec:examples}, we provide several numerical examples which demonstrate the good conservation properties and long-time stability of the proposed scheme.

\section{Incompressible Magnetohydrodynamics}\label{sec:ideal_mhd}

Magnetohydrodynamics describes fluids that carry an electric current
but remain electrically neutral so that the fluid motion is coupled
to the magnetic field only. The equations of
incompressible magnetohydrodynamics (MHD) result from the combination of the
Navier--Stokes equation for an incompressible fluid with an appropriate form of
the induction equation for the magnetic field.
Specifically, the system of incompressible MHD equations is given by
\begin{subequations}\label{eq:mhd_eqs}
\begin{align}\label{eq:mhd_eqs_V}
\partial_t V + ( V \cdot \nabla) V  &= -\nabla p + 
( B \cdot \nabla ) B + \mu \, \Delta V, \\
\label{eq:mhd_eqs_B}
\partial_t B + ( V \cdot \nabla) B - ( B \cdot \nabla ) V  &= \eta \, \Delta B , \\
\label{eq:divV}
\nabla \cdot V &= 0 , \\
\label{eq:divB}
\nabla \cdot B &= 0 ,
\end{align}
\end{subequations}
for $t \in [0,T]$, with $T>0$ sufficiently small, and on a bounded spatial domain 
$\Omega \subset \mathbb{R}^d$. The three unknowns are the fluid velocity
$V: \Omega_T \to \mathbb{R}^d$, the Alfv\'{e}n velocity
$B: \Omega_T \to \mathbb{R}^d$, 
i.e., the magnetic field divided by $\sqrt{4\pi \rho}$ (in
c.g.s.~units) with a constant mass density $\rho>0$, and the effective pressure
$p : \Omega_T \to \mathbb{R}$, i.e., the sum of the fluid pressure per unit of
mass and the magnetic pressure per unit of mass.  
Here, $\Omega_T = [0,T] \times \Omega$ and $\Delta$ denotes the standard Laplace
operator on $\mathbb{R}^d$. 
With a common abuse of terminology, we shall refer to $B$ as the magnetic field,
thus implying the proportionality constant, and to $p$ as the pressure. The
parameters $\mu>0$ and $\eta>0$ determine the strength of viscosity and electric
resistivity of the fluid, respectively.  

Equation (\ref{eq:mhd_eqs_V}) is called the \emph{momentum equation}, while
equation (\ref{eq:mhd_eqs_B}) is the \emph{induction equation}.
Both $V$ and $B$ are divergence-free, $V$ as we consider an incompressible fluid, and $B$ as there are no magnetic monopoles.
But while $\nabla \cdot B = 0$ is implied by the induction equation, provided
that the initial magnetic field $B|_{t=0}$ is divergence-free, 
$\nabla \cdot V = 0$ is a dynamical constraint determining the pressure $p$. 

In this paper, we shall consider the two-dimensional case, $d=2$, where $\Omega$
is a rectangular domain with Cartesian coordinates $(x,y)$, and $V=(V^x, V^y)$, 
$B= (B^x, B^y)$. We shall make use of periodic boundary conditions so that
$\Omega \cong \mathbb{T}^2$, the flat two-dimensional torus. 

\subsection{Ideal Incompressible Magnetohydrodynamics}\label{subsec:ideal_mhd}

Ideal incompressible MHD equations are obtained by setting $\mu=0$ and $\eta=0$
in equations~(\ref{eq:mhd_eqs}), thus considering the case of an ideal fluid
with zero electrical resistivity.
Without resistivity, $\eta = 0$, equation (\ref{eq:mhd_eqs_B}) states that the magnetic field is advected with the fluid flow, which implies the conservation of the magnetic flux through a surface moving with the fluid \cite{ArnoldKhesin:1998}. In addition, the topology of magnetic field lines is conserved. They are not allowed to open up and reconnect, a property that is worthwhile to maintain on the discrete level.
In a resistive plasma, $\eta \Delta B$ describes diffusive effects, for which the magnetic field lines are not just dragged along with the field, but are free to change their topology.

In two spatial dimensions ($d=2$), the constraints (\ref{eq:divV}) and
(\ref{eq:divB}) can be used to recast the initial value problem for ideal
incompressible MHD in the form
\begin{subequations}\label{eq:ideal_mhd_equations}
\begin{align}
\partial_{t} V + \psi(V,V) &= \psi(B,B) - \nabla P , && \text{on } \Omega_T , \\
\partial_{t} B + \phi(V,B) &= 0 , && \text{on } \Omega_T , \\
\nabla \cdot V &= 0 , && \text{on } \Omega_T , \\
V = V_0, \quad  B &= B_0, && \text{on } \{t=0\} \times \Omega ,
\end{align}
\end{subequations}
where the initial data $(V_0,B_0)$ must be divergence-free, $P = p + \frac{1}{2}(|V|^2 - |B|^2)$,
and we have introduced two bi-linear operators, $\psi$ and $\phi$, defined
componentwise by 
\begin{align}
\psi^{x} (V,B) &\equiv V^{y} \, \big( \partial_{y} B^{x} - \partial_{x} B^{y} \big) , &
\psi^{y} (V,B) &\equiv V^{x} \, \big( \partial_{x} B^{y} - \partial_{y} B^{x} \big) , \\
\phi^{x} (V,B) &\equiv \partial_{y} \big( V^{y} B^{x} - V^{x} B^{y} \big) , &
\phi^{y} (V,B) &\equiv \partial_{x} \big( V^{x} B^{y} - V^{y} B^{x} \big).
\end{align}
Although this formulation might not appear natural at first, the operators
$\psi$ and $\phi$ have a geometrical meaning which will become clear in section~\ref{subsec:forms}.
The same bilinear operators were found by \citeauthor{Gawlik:2011}
\cite{Gawlik:2011} in their discretisation of the Euler-Poincar\'e variational
formulation of ideal MHD. 

A direct calculation shows that, for any triple of vector fields $B$, $V$, and
$W$, one has
\begin{subequations}\label{eq:properties_psi_phi}
\begin{align}
  \label{eq:property_psi1}
  W \cdot \psi(V, B) &= - V \cdot \psi(W, B) , \\  
  \label{eq:property_psi2}
  V \cdot \psi(V, B) &= 0, \\
  \label{eq:property_phi1}
  \phi(V, B) &= - \phi(B, V) , \\  
  \label{eq:property_psiphi}
  \int_\Omega B \cdot \phi(V, W) d\tau &= \int_\Omega V \cdot \psi(W, B)d\tau,
\end{align}
\end{subequations}
where $d\tau = dxdy$ denotes the volume element on $\Omega$ and periodic
boundary conditions have been used for (\ref{eq:property_psiphi}). By making use
of the foregoing identities, one readily obtains the three important conserved
quantities of ideal MHD in two dimensions \cite{ArnoldKhesin:1998}, namely, the
conservation of the total energy,  
\begin{align}\label{eq:ideal_mhd_energy}
E = \dfrac{1}{2} \int_\Omega \Big[ |V|^{2} + |B|^{2} \Big] d\tau,
\end{align}
of cross helicity
\begin{align}\label{eq:ideal_mhd_cross_helicity}
C_{\mrm{CH}} = \int_\Omega V \cdot B d\tau,
\end{align}
and of magnetic helicity
\begin{align}\label{eq:ideal_mhd_magnetic_helicity}
C_{\mrm{MH}} = \int_\Omega A_z \, d\tau,
\end{align}
where $A_z$ is the third component of the magnetic vector potential, 
$B = (-\partial_y A_z, \partial_x A_z)$; from the induction equation one has
that $\partial_t A_z + V \cdot \nabla A_z = 0$, hence $dC_{\mrm{MH}}/dt = 0$.
Conservation of all three quantities is desirable in numerical simulations in order to obtain qualitatively and quantitatively accurate results.
In the next step, we construct a formal Lagrangian for equations \eqref{eq:ideal_mhd_equations}.

\subsection{Formal Lagrangians}\label{subsec:formal_lagrangians}

A direct variational formulation of system~\eqref{eq:ideal_mhd_equations} does not exist since the functional derivative of the MHD equations is not symmetric~\cite{Vainberg:1964}.
Therefore, in order to derive variational integrators for these equations, we resort to a formal Lagrangian formulation~\cite{Ibragimov:2006, AthertonHomsy:1975}.
This approach to the development of numerical integration schemes is described in detail in reference~\cite{KrausMaj:2015}. Here we outline the procedure for the case under consideration without theoretical (geometrical and functional) details.
Essentially, we treat the ideal MHD system as part of a larger system, which admits a Lagrangian formulation.
Each equation of (\ref{eq:ideal_mhd_equations}) as well as the incompressibility constraint is multiplied with auxiliary variables, $\alpha$, $\beta$ and $\gamma$, respectively.
It is worth mentioning that $\nabla \cdot V = 0$ is treated as a dynamical equation determining the pressure.

The formal Lagrangian density is given as the sum of the resulting expressions,
\begin{align}\label{eq:ideal_mhd_formal_lagrangian}
\mcal{L} (\phy, \phy_{t}, \phy_{x}, \phy_{y})
\nonumber
&= \alpha \cdot \big[ \partial_{t} V + \psi(V,V) - \psi(B,B) + \nabla P \big] \\
&+ \beta  \cdot \big[ \partial_{t} B + \phi(V,B) \big]
 + \gamma \, \big[ \nabla \cdot V \big] ,
\end{align}
where $\phy$ denotes the tuple of all the variable of the extended system,
\begin{align}
\phy = (V, B, P, \alpha, \beta, \gamma) ,
\end{align}
and $\phy_{t}$, $\phy_{x}$, and $\phy_{y}$ denote the corresponding derivatives with respect to $t$, $x$, and $y$, respectively.
The ideal MHD equations (\ref{eq:ideal_mhd_equations}) are obtained from Hamilton's variational principle applied to the action
\begin{align}
\mcal{A} [\phy] = \int_{0}^{T} L (\phy (t), \phy' (t)) \, dt ,
\end{align}
where
\begin{align}
L (\phy (t), \phy' (t)) = \int \limits_{\Omega} \mcal{L} (\phy, \phy_{t}, \phy_{x}, \phy_{y}) \, d\tau ,
\end{align}
and $\phy (t) = \phy (t , \cdot)$, $\phy' (t) = \partial_{t} \phy (t , \cdot)$, with fixed initial and final points, i.e., $\phy(0) = \phy_0$, $\phy(T) = \phy_T$.
In addition we obtain a set of equations which determine the evolution of the auxiliary variables,
\begin{subequations}\label{eq:adjoint}
\begin{align}
\partial_{t} \alpha + \psi (\alpha,V) &= \psi(B, \beta)  + \phi (\alpha,V) - \nabla \gamma , \\
\partial_{t} \beta  + \phi (\alpha,B) &= \psi(\alpha, B) - \psi(V, \beta) , \\
\nabla \cdot \alpha &= 0 .
\end{align}
\end{subequations}
A remarkable property of the Lagrangian~(\ref{eq:ideal_mhd_formal_lagrangian})
is that, if $V$, $B$, and $P$ solve the physical
equations~(\ref{eq:ideal_mhd_equations}), then we can construct a solution of
the adjoint system~(\ref{eq:adjoint}) by setting $\alpha = V$, $\beta=B$, and
$\gamma=P$. 

\subsection{Differential Forms}\label{subsec:forms}

As the derivation of the proposed scheme will be based on discrete exterior calculus, we
need to reformulate the ideal incompressible MHD
equations~\eqref{eq:ideal_mhd_equations} and the formal
Lagrangian~\eqref{eq:ideal_mhd_formal_lagrangian} in terms of
exterior calculus and differential forms~\cite{AbrahamMarsdenRatiu:1988,
  Lee:2012, Lee:2009, Tu:2011}. We begin with a minimal simplified overview of the
necessary formalism.
For a more comprehensive review in the context of discretisation, the reader is referred to~\cite{Gerritsma:2014, Palha:2014}.

On a two-dimensional domain $\Omega$, there are three non-trivial spaces of
differential forms denoted by $\Lambda^p = \Lambda^p(\Omega)$ for $p \in \{0,1,2\}$. The space 
$\Lambda^0(\Omega)$ comprises scalar functions $f: \Omega \to \mathbb{R}$, the
zero-forms. The space of one-forms $\Lambda^1(\Omega)$ comprises line elements 
$\alpha = \alpha_x(x,y) \ext x + \alpha_y(x,y) \ext y$ where $\alpha_x$ and $\alpha_y$ are functions over $\Omega$.
Here, the differentials $\ext x$ and $\ext y$ are precisely defined as linear maps
acting on a vector $W=(W^x, W^y)$ tangent to $\Omega$ according to 
$\ext x(W) = W^x$ and $\ext y(W) = W^y$; hence a one-form corresponds to a map from
$\Omega$ to the dual space of the tangent vectors.  At
last, the space of two-forms $\Lambda^2(\Omega)$ comprises surface elements
$m(x,y) \ext x \wedge \ext y$, where we have defined the exterior product 
$\ext x \wedge \ext y = \ext x \otimes \ext y - \ext y \otimes \ext x$, 
with $\otimes$ being the tensor product of linear operators. This is a bi-linear
map over the tangent space,
$(W_1, W_2) \mapsto \ext x \wedge \ext y (W_1, W_2) = W_1^x W_2^y - W_1^y W_2^x$,
which corresponds to the Euclidean area of the quadrilateral spanned by the
vectors $W_1$ and $W_2$, so that $\ext x \wedge \ext y$ can be identified with
the surface element in $\Omega$ for the Euclidean metric.
The union of the spaces $\Lambda^p$ forms a graded algebra with respect to the
alternating product, which has been introduced above for $\ext x \wedge \ext y$ and can be defined
as a map $\wedge : \Lambda^p \times \Lambda^q \to \Lambda^{p+q}$ so that
$(\alpha, \beta) \mapsto \alpha \wedge \beta \in \Lambda^{p+q}$ for any pair of
forms $\alpha \in \Lambda^p$ and $\beta \in \Lambda^q$. This map can, for
simplicity, be characterised by its action on elementary objects, namely, $f \wedge \ext x = f \ext x$, $f \wedge \ext y = f \ext y$, 
for $f\in\Lambda^0$, $\ext x \wedge \ext x = \ext y \wedge \ext y = 0$, and associativity 
$(\alpha \wedge \beta) \wedge \gamma = \alpha \wedge (\beta \wedge \gamma)$.
The exterior product is not commutative and we have $\alpha \wedge \beta = (-1)^{pq} \beta \wedge \alpha$. 

We can introduce the exterior derivative $\ext : \Lambda^p \to \Lambda^{p+1}$ in
the same way: It acts on zero-forms as 
$\ext f(x,y) = \partial_x f(x,y) \ext x + \partial_yf(x,y) \ext y$ and on one-forms as
$\ext \alpha = (\partial_x \alpha_y - \partial_y \alpha_x) \ext x \wedge \ext y$.
Then, we have
$\ext^2 f = (\partial_x \partial_y f - \partial_y \partial_x f) \ext x \wedge \ext y = 0$
for any function $f \in \Lambda^0(\Omega)$ (assuming $C^2$ regularity), that is,  
$\text{Rng}(\ext: \Lambda^0 \to \Lambda^1) \subseteq \text{Ker}(\ext: \Lambda^1 \to \Lambda^2)$, 
where Rng and Ker denote the range and the null space (kernel)
of the operator, respectively. Constants are embedded in 
$\Lambda^0$ and $\mathbb{R} = \text{Ker}(\ext: \Lambda^0 \to \Lambda^1)$. All
spaces $\Lambda^p(\Omega)$ for $p > 2$ reduce to the trivial subspace $\{0\}$ of
the space of $p$-linear tensors, as the only alternating $p$-linear tensor with
$p>d$ is the trivial one. We therefore have $\ext : \Lambda^2(\Omega) \to \{0\}$.
The product rule for the exterior differentiation reads 
$\ext (\alpha \wedge \beta) = \ext \alpha \wedge \beta + (-1)^p \alpha \wedge \ext \beta$, 
for $\alpha \in \Lambda^p$ and $\beta \in \Lambda^q$, which is Leibnitz formula in exterior calculus. 

Together with the standard spaces of differential forms, we consider their
twisted counterparts $\widetilde{\Lambda}^p(\Omega)$, \cite{Burke:1985, Frankel:2011, Bossavit:1998, Kreeft:2011, Palha:2014}. In the simple domain under
consideration such spaces are just copies of $\Lambda^{p}(\Omega)$. The Hodge
operator $\hodge : \Lambda^p (\Omega) \to \widetilde{\Lambda}^{d-p}(\Omega)$  
acts on the various forms according to 
$\hodge f = f \ext x \wedge \ext y$ for $f \in \Lambda^0(\Omega)$,  
$\hodge \alpha = -\alpha_y \ext x + \alpha_x \ext y$ for 
$\alpha = (\alpha_x \ext x + \alpha_y \ext y) \in \Lambda^1(\Omega)$, and 
$\hodge \mu = m$ for $\mu = m \ext x \wedge \ext y \in \Lambda^2(\Omega)$.
Since in this case $\widetilde{\Lambda}^p = \Lambda^p$ the Hodge operator can be composed
with itself and one can check that $\hodge \hodge \alpha = (-1)^{p(d-p)} \alpha$
for any $\alpha \in \Lambda^p(\Omega)$. This structure is summarised in the
following diagram 
\begin{align}\label{eq:de_rham_complex_2d_forms}
\begin{matrix}
\rsp & \rightarrow & \Lambda^{0} (\Omega) & \xrightarrow{\ext} & \Lambda^{1} (\Omega) & \xrightarrow{\ext} & \Lambda^{2} (\Omega) & \rightarrow & 0 \\
&& \updownarrow \hodge && \updownarrow \hodge && \updownarrow \hodge && \\
0 & \leftarrow & \widetilde{\Lambda}^{2} (\Omega) & \xleftarrow{\ext} & \widetilde{\Lambda}^{1} (\Omega) & \xleftarrow{\ext} & \widetilde{\Lambda}^{0} (\Omega) & \leftarrow & \rsp .
\end{matrix}
\end{align}
Due to the fact that the range of $\ext$ defined on elements of $\Lambda^{p}$ is in the kernel of $\ext$ defined on elements of $\Lambda^{p+1}$, each row of the diagram forms a de Rham complex.

In addition, we can define the following operations. To a one-form 
$\alpha = \alpha_x \ext x + \alpha_y \ext y$ we can associate a vector field
$\alpha^\sharp = (\alpha_x, \alpha_y)$ and this map corresponds to the raising (or sharp)
operator $\sharp$ with the trivial metric. Conversely, to a vector field $W = (W^x, W^y)$ we can
associate a one form $W^\flat = W^x \ext x + W^y \ext y$ and this map
corresponds to the lowering (or flat) operator $\flat$ with the trivial
metric. For the contraction of a vector field and a one-form we write,
\begin{equation*}
  \iprod_W \alpha = \alpha(W) = \alpha_x W^x + \alpha_y W^y,
\end{equation*}
while the contraction with a two form $\mu = m \ext x \wedge \ext y$ we write
\begin{equation*}
  \iprod_W \mu = m (-W^y \ext x + W^x \ext y).
\end{equation*}

We can now write the MHD equations
in terms of such objects. In this framework, vector fields and one-forms are
isomorphic to each other, both representing two-component fields. We choose to
represent both $V$ and $B$ as one-forms, i.e., 
$V = V^x \ext x + V^y \ext y$ and $B = B^x \ext x + B^y \ext y$. 
The momentum balance suggests  $P \in \Lambda^0$ so that $\ext P$ represents the
pressure gradient as a one-form. At last we have the identities 
\begin{subequations}\label{eq:psi_phi_forms}
\begin{align}
\label{eq:psi_form}
\psi^x(V, B) \ext x + \psi^y(V, B) \ext y &= \iprod_{V^\sharp} \ext B, \\
\label{eq:phi_form}
\phi^x(V,B) \ext x + \phi^y(V,B) \ext y &= -\hodge \ext (\iprod_{V^\sharp} \hodge B).
\end{align}
\end{subequations}
The divergence-free constraint takes the form
\begin{equation*}
  \hodge \ext \hodge V = \partial_x V^x + \partial_y V^y = 0,
\end{equation*}
and analogously for $B$. All the foregoing identities can be checked by direct
computation.

The Cauchy problem \eqref{eq:ideal_mhd_equations} for ideal incompressible MHD
equations then writes
\begin{subequations}\label{eq:ideal_mhd_forms}
\begin{align}
\partial_{t} V + \iprod_{V^\sharp} \ext V &= \iprod_{B^\sharp} \ext B  - \ext P , && \text{on } \Omega_T , \\
\partial_{t} \hodge B + \ext (\iprod_{V^\sharp} \hodge B) &= 0 , && \text{on } \Omega_T , \\
\hodge \ext \hodge V &= 0 , && \text{on } \Omega_T , \\
V (0) = V_0 \in \Lambda^1, \quad  B (0) &= B_0 \in \Lambda^1, && \text{on } \Omega, 
\end{align}
\end{subequations}
for $V,B: [0,T] \to \Lambda^1(\Omega)$, $P : [0,T] \to \Lambda^0(\Omega)$, and with divergence-free initial data $(V_0, B_0)$. One should notice that we have applied the Hodge operator to
the induction equation, and used the fact that $\hodge$ commutes with
$\partial_t$ and $\hodge \hodge \alpha = -\alpha$ for any $\alpha \in \Lambda^1$. 

We remark that terms on the right hand side of (\ref{eq:psi_phi_forms}) are
related to the Lie derivative via Cartan's formula
\begin{equation}\label{eq:Cartan_magic}
\lie_{V^{\sharp}} \alpha = \ext ( \iprod_{V^{\sharp}} \alpha ) + \iprod_{V^{\sharp}} ( \ext \alpha ) .
\end{equation}
Particularly, $\psi(V,V)$ corresponds to the term in $\lie_{V^\sharp} V$, which
is not an exact differential (exact differentials can be reabsorbed by the
pressure gradient), while $\phi(V,B)$ is the full Lie derivative of
$\hodge B$, with $\hodge \ext \hodge B = 0$, and thus accounts for the
advection of $\hodge B$ along the flow of $V$.

Based on the formulation~(\ref{eq:ideal_mhd_forms}), we want to construct a formal
Lagrangian. To this end, we need to construct the $L^2$-scalar product for differential forms 
$\alpha, \beta \in \Lambda^p(\Omega)$, for all $p \in \{0,1,2\}$. This is achieved
on noticing that $\alpha \wedge \hodge \beta$ is always a surface element, and
thus it can be integrated over $\Omega$; moreover, it is symmetric, i.e., 
$\alpha \wedge \hodge \beta = \beta \wedge \hodge \alpha$. This allows us to
define 
\begin{equation}\label{eq:forms_inner_product}
  \pair{\alpha}{\beta} = \int \limits_\Omega \alpha \wedge \hodge \beta,
\end{equation}
which provides an inner product on $\Lambda^p(\Omega)$. With respect to this inner
product, the Hodge operator is either symmetric or anti-symmetric depending on
$p$, namely, for any two forms $\alpha \in \Lambda^p(\Omega)$ and $\beta \in \Lambda^{d-p}(\Omega)$ 
we have
\begin{equation*}
  \pair{\hodge \alpha}{\beta} = (-1)^{p(d-p)} \pair{\alpha}{\hodge \beta},
\end{equation*} 
and specializing this result to the case $\beta = \hodge \gamma$ with 
$\gamma \in \Lambda^p(\Omega)$ we obtain
\begin{equation*}
  \pair{\hodge \alpha}{\hodge \gamma} = \pair{\alpha}{\gamma},
\end{equation*}
that is, $\hodge$ is an $L^2$-isometry. We shall also use the fact that the
contraction operator $\iprod_{V^\sharp}$ with $V\in\Lambda^1(\Omega)$ is the formal
adjoint to the operator of exterior multiplication $V\wedge \cdot$, namely, for 
$\mu \in \Lambda^{p+1}(\Omega)$ and $\nu \in \Lambda^p(\Omega)$ we have
\begin{equation*}
  \pair{\iprod_{V^\sharp} \mu}{\nu} = \pair{\mu}{V \wedge \nu}.
\end{equation*}
Analogously we also have that, in two dimensions, $-\hodge \ext \hodge$ is the
formal adjoint of $\ext$, i.e., for $\mu \in \Lambda^p(\Omega)$ and $\nu \in \Lambda^{p+1}(\Omega)$
\begin{equation*}
  \pair{\ext \mu}{\nu} + \pair{\mu}{\hodge \ext \hodge \nu} = 0,
\end{equation*}
which follows from the product rule on noticing that
\begin{align*}
  \ext (\mu\wedge \hodge \nu) 
  &= \ext \mu \wedge \hodge \nu + (-1)^p \mu \wedge \ext \hodge \nu \\
  &= \ext \mu \wedge \hodge \nu + (-1)^{p(d-p+1)} \mu \wedge \hodge (\hodge \ext \hodge \nu) 
\end{align*}
and in two dimensions $(-1)^{p(d-p+1)} = 1$ while integrating on a periodic
domain annihilates the left-hand side.

With the inner product~\eqref{eq:forms_inner_product} we can construct the formal Lagrangian by pairing the
momentum equation with $\alpha \in \Lambda^1(\Omega)$, the induction equation with a
twisted one-form $\hodge \beta$, $\beta \in \Lambda^1(\Omega)$, and the divergence
constraint of the velocity with $\gamma \in \Lambda^0(\Omega)$,  
\begin{align}
L = \pair{V_t + \iprod_{V^{\sharp}} ( \ext V ) - \iprod_{B^{\sharp}} ( \ext B ) + \ext P}{\alpha}
  + \pair{\hodge B_t + \ext ( \iprod_{V^{\sharp}} \hodge B )}{\hodge \beta}
  + \pair{\hodge \ext \hodge V}{\gamma} , 
\end{align}
where $V_t = \partial_t V$ and $B_t = \partial_tB$. We want to recast the
Lagrangian in terms of forms only, thus removing the vector fields $V^\sharp$
and $B^\sharp$.
With that aim, we exploit the properties of the inner product to write
\begin{align*}
  L (\phy, \phy')
  &= \pair{V_t}{\alpha} + \pair{\ext V}{V\wedge\alpha}
  - \pair{\ext B}{B\wedge\alpha} + \pair{\ext P}{\alpha} \\
  &\qquad
  + \pair{\hodge B_t}{\hodge \beta} 
  - \pair{\iprod_{V^{\sharp}} \hodge B }{\hodge \ext \hodge \hodge \beta}
  + \pair{\hodge \ext \hodge V}{\gamma} \\
  &= \pair{V_t}{\alpha} + \pair{\ext V}{V\wedge\alpha}
  - \pair{\ext B}{B\wedge\alpha} + \pair{\ext P}{\alpha} \\
  &\qquad
  + \pair{\hodge B_t}{\hodge \beta} 
  + \pair{\iprod_{V^{\sharp}} \hodge B }{\hodge \ext \beta}
  + \pair{\hodge \ext \hodge V}{\gamma} \\
  &= \pair{V_t}{\alpha} + \pair{\ext V}{V\wedge\alpha}
  - \pair{\ext B}{B\wedge\alpha} + \pair{\ext P}{\alpha} \\
  &\qquad
  + \pair{B_t}{\beta} 
  + \pair{\hodge \iprod_{V^{\sharp}} \hodge B }{\ext \beta}
  + \pair{\hodge \ext \hodge V}{\gamma},
\end{align*}  
and at last we note that $\hodge \iprod_{V^{\sharp}} \hodge B = -V \wedge B$ with
the result that
\begin{align}\label{eq:ideal_mhd_formal_lagrangian_forms}
L (\phy, \phy')
\nonumber
&= \pair{ V_{t} }{ \alpha }
 + \pair{ \ext V }{ V \wedge \alpha } 
 - \pair{ \ext B }{ B \wedge \alpha } 
 + \pair{ \ext P }{ \alpha } \\
&\qquad + \pair{ B_{t} }{ \beta }
 - \pair{ V \wedge B }{ \ext \beta }
 + \pair{ \hodge \ext \hodge V }{ \gamma } .
\end{align}
We observe that all quantities are differential forms and no vector field is present anymore.
The Lagrangian~\eqref{eq:ideal_mhd_formal_lagrangian_forms} will be the basis for the discretisation outlined in the next section.

\section{Discrete Exterior Calculus on Cartesian Meshes}\label{sec:dec}

In this section we define basic elements of the discrete exterior calculus on two-dimensional Cartesian meshes.
With this aim we follow the work of~\citet{Robidoux:2011}, who developed a mimetic discretisation method for vector calculus on tensor-product meshes.
For discrete exterior calculus on simplicial meshes we refer to the work of~\citet{Hirani:2003} and~\citet{Desbrun:2008}.

In order to obtain a finite difference discretisation of spaces of differential forms, one can exploit the relationship between zero-forms and points, one-forms and line elements, two-forms and area elements, c.f. Section~\ref{subsec:forms}.
This allows us to identify discrete differential forms with their integral over geometric elements of the grid: vertices for zero-forms, edges for one-forms, and cells for two-forms.
Within this framework, it is convenient to have a linear space structure over the geometric elements of the grid, which leads to the concept of chains~\cite{Robidoux:2011}.

\begin{figure}[tb]
\centering
\subfloat[Primal Grid]{\label{fig:mhd_staggered_grid_momentum}
\includegraphics[width=.44\textwidth]{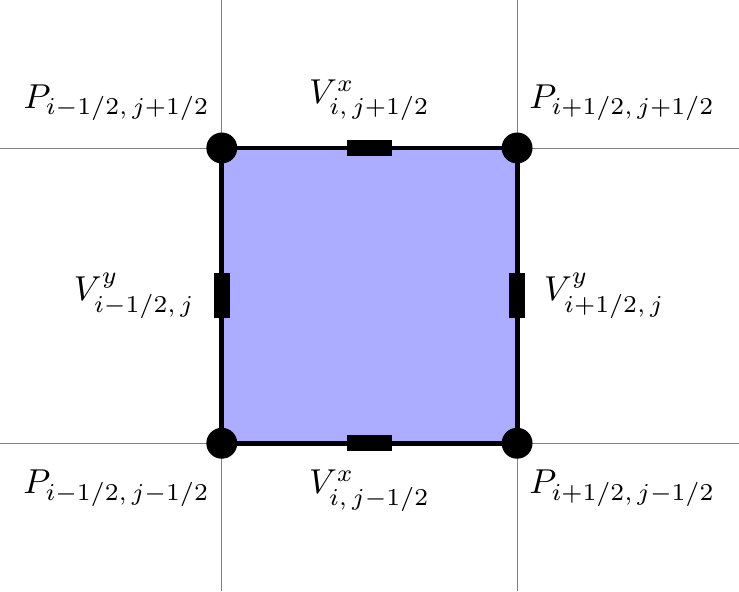}
}
\subfloat[Dual Grid]{\label{fig:mhd_staggered_grid_div}
\includegraphics[width=.44\textwidth]{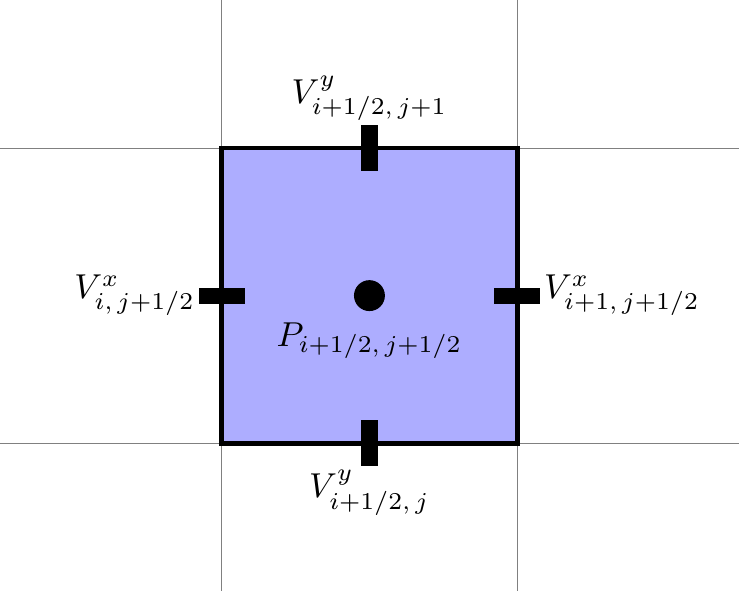}
}
\caption{Staggered grid in the $xy$-plane. Left: Primal grid with natural positions for the pressure and the velocity components for the computation of the advection operators. Right: Dual grid for the computation of the divergence constraint.}
\label{fig:mhd_staggered_grid}
\end{figure}

We introduce a staggered grid, where the pressure is collocated at the vertices of a grid cell and the velocity and magnetic field components at the edges, c.f. Figure~\subref*{fig:mhd_staggered_grid_momentum}.
The location of the physical quantities comes natural when viewed as differential forms. The pressure is a zero-form and is therefore collocated at the vertices of a cell of the primal grid. The velocity (and in two dimensions also the magnetic field) is a one-form and is therefore collocated at the edges of a cell, $x$-components on the horizontal edges and $y$-components on the vertical edges, c.f. Figure~\subref*{fig:mhd_staggered_grid_momentum}.
On the dual grid, the pressure becomes a two-form, collocated at the cell centre. The velocity and magnetic field are still one-forms, but twisted, so that $x$-components are collocated on the vertical edges and $y$-components on the horizontal edges, c.f. Figure~\subref*{fig:mhd_staggered_grid_div}.

\subsection{The Grid}

We consider a grid in the two-dimensional Euclidean space $\rsp^{2}$, given by equidistant points in each direction. On the primal grid, cell centers are labelled with integer indices, $(i,j)$, vertices are labelled with half-integer indices, $(i+1/2, \, j+1/2)$ and edges with mixed indices, $(i, \, j+1/2)$ and $(i+1/2, \, j)$. On the dual grid, the labelling is reversed.
Points in the grid correspond to $x$- and $y$-coordinates as follows,
\begin{align}
x_{i} &= i h_{x} , &
y_{j} &= j h_{y} , &
x_{i+1/2} &= (i + \tfrac{1}{2}) h_{x} , &
y_{j+1/2} &= (j + \tfrac{1}{2}) h_{y} ,
\end{align}
with $h_{x}, h_{y} \in \rsp^{+}$ and $i,j,k \in \zsp$.
Here, $h_{x}$ and $h_{y}$ denote the grid step size in $x$- and $y$-direction, respectively, which is assumed to be constant throughout the grid.

\subsection{Chains}

Let us start with the definition of cell chains.
Given $a_{i,j} \in \{ -1, 0, +1 \}$, the formal sum
\begin{align}
c_{h} = \sum \limits_{i,j} a_{i, \, j} \, c_{i, \, j} ,
\end{align}
where $c_{i,j}$ is a cell of the primal grid, is interpreted as the disjoint union of cells $c_{i,j}$.
Here, $a_{i,j} = \pm 1$ indicates the cell's orientation and $a_{i,j} = 0$ means that the cell $c_{i,j}$ is not present.
The set of such formal series can be extended to a linear space by taking the coefficients $a_{i,j}$ in $\mbb{R}$ and defining summation and multiplication by scalars according to
\begin{align}
\alpha c_{h} + \beta c_{h}' = \sum \limits_{i,j} \big( \alpha a_{i,j} + \beta a_{i,j}' \big) \, c_{i,j} .
\end{align}
Such formal sums are referred to as cell chains and their space is denoted by $\mcal{C}$.
They retain their geometrical meaning for $a_{i,j} \in \{ -1, 0, +1 \}$ and should be regarded as abstract algebraic objects for general coefficients.

Analogously, we define the space $\mcal{E}$ of edge chains
\begin{align}
e_{h} = \sum \limits_{i,j} \Big[ a^{x}_{i, \, j+1/2} \, e^{x}_{i, \, j+1/2} + a^{y}_{i+1/2, \, j} \, e^{y}_{i+1/2, \, j} \Big] ,
\end{align}
where $e^{x}_{i, \, j+1/2}$ and $e^{y}_{i+1/2, \, j}$ are the top and right edges of the cell $c_{i,j}$.
We always consider periodic boundary conditions, so that we have as many horizontal and vertical edges as cells.

At last, the space $\mcal{V}$ of vertex chains comprises formal sums of the form
\begin{align}
v_{h} = \sum \limits_{i,j} a_{i+1/2, \, j+1/2} \, v_{i+1/2, \, j+1/2} ,
\end{align}
where $v_{i+1/2, \, j+1/2}$ is the upper right corner of the cell $c_{i,j}$. 

One should notice that cells, edges and vertices of the grid act as basis of the corresponding spaces of cell, edge and vertex chains, respectively.
Among these spaces, we have boundary operators
\begin{align}
\begin{matrix}
0 & \xleftarrow{\partial} & \mcal{V} & \xleftarrow{\partial} & \mcal{E} & \xleftarrow{\partial} & \mcal{C} ,
\end{matrix}
\end{align}
which are the linear operators defined on the basis according to
\begin{subequations}\label{eq:dec_boundary}
\begin{align}
\partial v_{i+1/2, \, j+1/2} &= 0 , \\
\partial e^{x}_{i, \, j+1/2} &= v_{i+1/2, \, j+1/2} - v_{i-1/2, \, j+1/2} , \\
\partial e^{y}_{i+1/2, \, j} &= v_{i+1/2, \, j+1/2} - v_{i+1/2, \, j-1/2} , \\
\partial c_{i, \, j} &= e^{x}_{i, \, j+1/2} - e^{x}_{i, \, j-1/2} + e^{y}_{i+1/2, \, j} - e^{y}_{i-1/2, \, j} ,
\end{align}
\end{subequations}
and extended by linearity to the whole spaces.
For example, the boundary operator is applied to a cell chain $c_{h} \in \mcal{C}$ by applying $\partial$ to each term in the sum, giving an edge chain $\partial c_{h} \in \mcal{E}$ of the form
\begin{align}
\partial c_{h}
= \sum \limits_{i,j} a_{i, \, j} \, \partial c_{i, \, j} .
\end{align}
Note that $\partial \partial = 0$, which is verified in direct calculation of the boundaries of the expressions in \eqref{eq:dec_boundary}.

Analogously, chains on the dual grid are defined as
\begin{subequations}
\begin{align}
v_{h}^{\hodge} &= \sum \limits_{i,j} a_{i, \, j} \, v_{i, \, j}^{\hodge} , \\
e_{h}^{\hodge} &= \sum \limits_{i,j} \Big[ a^{x}_{i+1/2, \, j} \, e^{\hodge x}_{i+1/2, \, j} + a^{y}_{i, \, j+1/2} \, e^{\hodge y}_{i, \, j+1/2} \Big] , \\
c_{h}^{\hodge} &= \sum \limits_{i,j} a_{i+1/2, \, j+1/2} \, c_{i+1/2, \, j+1/2}^{\hodge} ,
\end{align}
\end{subequations}
where $v_{i, \, j}^{\hodge}$, $e^{\hodge x}_{i+1/2, \, j}$, $e^{\hodge y}_{i, \, j+1/2}$ and $c_{i+1/2, \, j+1/2}^{\hodge}$ denote vertices, edges, cells and thus basis elements of chains on the dual grid.
The corresponding spaces of vertex, edge and cell chains are denoted by $\mcal{V}^{\hodge}$, $\mcal{E}^{\hodge}$ and $\mcal{C}^{\hodge}$, respectively.
The action of the boundary operator on chains on the dual grid is defined in complete analogy to the boundary operator on the primal grid.

\subsection{Differential Forms}

The spaces of discrete differential forms can be defined as the algebraic duals of the spaces of chains (therefore discrete differential forms are also referred to as cochains).
Specifically, discrete zero-forms are linear operators from $\mcal{V}$ to $\mbb{R}$ and their space is denoted by $\Lambda_{h}^{0}$.
A basis for $\Lambda_{h}^{0}$ is then given by the linear operators $\obar{v}_{k+1/2, \, l+1/2}$, that act on the basis elements $v_{i+1/2, \, j+1/2}$ of $\mcal{V}$ as
\begin{subequations}\label{eq:dec_chain_basis_cochain_basis}
\begin{align}
\obar{v}_{k+1/2, \, l+1/2} (v_{i+1/2, \, j+1/2})
&= \int \limits_{v_{i+1/2, \, j+1/2}} \obar{v}_{k+1/2, \, l+1/2}
 = \delta_{ik} \delta_{jl} ,
\end{align}
where the formal integral denotes the action of a discrete form on a chain.
Analogously, the space $\Lambda_{h}^{1}$ of discrete one-forms is the dual of $\mcal{E}$. A basis for $\Lambda_{h}^{1}$ is given by the linear operators $\obar{e}^{x}_{k, \, l+1/2}$ and $\obar{e}^{y}_{k+1/2, \, l}$ acting on the basis elements $e^{x}_{i, \, j+1/2}$ and $e^{y}_{i+1/2, \, j}$ of $\mcal{E}$ by
\begin{align}
\obar{e}^{x}_{k, \, l+1/2} (e^{x}_{i, \, j+1/2})
&= \int \limits_{e^{x}_{i, \, j+1/2}} \obar{e}^{x}_{k, \, l+1/2}
 = \delta_{ik} \delta_{jl} h_{x} , \\
\obar{e}^{x}_{k, \, l+1/2} (e^{y}_{i+1/2, \, j})
&= \int \limits_{e^{y}_{i+1/2, \, j}} \obar{e}^{x}_{k, \, l+1/2}
 = 0 ,
\end{align}
and
\begin{align}
\obar{e}^{y}_{k+1/2, \, l} (e^{x}_{i, \, j+1/2})
&= \int \limits_{e^{x}_{i, \, j+1/2}} \obar{e}^{y}_{k+1/2, \, l}
 = 0 , \\
\obar{e}^{y}_{k+1/2, \, l} (e^{y}_{i+1/2, \, j})
&= \int \limits_{e^{y}_{i+1/2, \, j}} \obar{e}^{y}_{k+1/2, \, l}
 = \delta_{ik} \delta_{jl} h_{y} .
\end{align}
At last, the space $\Lambda_{h}^{2}$ of two-forms is the dual of the space $\mcal{V}$ of cell chains.
The basis for $\Lambda_{h}^{2}$ is given by the linear operators $\obar{c}_{k, \, l}$, defined by
\begin{align}
\obar{c}_{k, \, l} (c_{i, \, j})
= \int \limits_{c_{i, \, j}} \obar{c}_{k, \, l}
= \delta_{ik} \delta_{jl} h_{x} h_{y} .
\end{align}
\end{subequations}

Written in this basis, the discrete forms $\phi_{h} \in \Lambda_{h}^{0}$, $\alpha_{h} \in \Lambda_{h}^{1}$ and $\omega_{h} \in \Lambda_{h}^{2}$ take the following form,
\begin{subequations}\label{eq:discrete_diff_forms}
\begin{align}
\phi_{h}   &= \sum \limits_{i,j} \phi_{i+1/2, \, j+1/2} \, \obar{v}_{i+1/2, \, j+1/2} , \\
\alpha_{h} &= \sum \limits_{i,j} \Big[ \alpha^{x}_{i, \, j+1/2} \, \obar{e}^{x}_{i, \, j+1/2} + \alpha^{y}_{i+1/2, \, j} \, \obar{e}^{y}_{i+1/2, \, j} \Big] , \\
\omega_{h} &= \sum \limits_{i,j} \omega_{i, \, j} \, \obar{c}_{i, \, j} .
\end{align}
\end{subequations}

Spaces of twisted differential forms are defined as the algebraic duals of the chain spaces $\mcal{V}^{\hodge}$, $\mcal{E}^{\hodge}$ and $\mcal{C}^{\hodge}$ on the dual grid, and are denoted by $\widetilde{\Lambda}_{h}^{p}$ with $p = 0, 1, 2$ and with basis elements $\obar{v}^{\hodge}_{i, \, j}$, $\obar{e}^{\hodge x}_{i+1/2, \, j}$, $\obar{e}^{\hodge y}_{i, \, j+1/2}$ and $\obar{c}^{\hodge}_{i+1/2, \, j+1/2}$.

\subsection{Integrals}\label{sec:dec_integrals}

By linearity it follows from the basis representation~\eqref{eq:discrete_diff_forms}, that the zero-form $\phi_{h}$ applied to the vertex chain $v_{h}$ is
\begin{align}\label{eq:integral_one_form}
\phi_{h} (v_{h}) = \int \limits_{v_{h}} \phi_{h} = \sum \limits_{i,j} a_{i+1/2, \, j+1/2} \, \phi_{i+1/2, \, j+1/2} ,
\end{align}
the one-form $\alpha_{h}$ applied to the edge chain $e_{h}$ is
\begin{align}\label{eq:integral_two_form}
\alpha_{h} (e_{h}) = \int \limits_{e_{h}} \alpha_{h} = \sum \limits_{i,j} \Big[ h_{x} \, a^{x}_{i, \, j+1/2} \, \alpha^{x}_{i, \, j+1/2} + h_{y} \, a^{y}_{i+1/2, \, j} \, \alpha^{y}_{i+1/2, \, j} \Big] ,
\end{align}
and the two-form $\omega_{h}$ applied to the cell chain $c_{h}$ is
\begin{align}\label{eq:integral_three_form}
\omega_{h} (c_{h}) = \int \limits_{c_{h}} \omega_{h} = \sum \limits_{i,j} h_{x} h_{y} \, a_{i, \, j} \, \omega_{i, \, j} .
\end{align}
For the specific cases of discrete one- and two-forms, the formal integrals~\eqref{eq:integral_two_form} and~\eqref{eq:integral_three_form} can be identified with the integration~\cite{AbrahamMarsdenRatiu:1988,
  Lee:2012, Lee:2009, Tu:2011} of piecewise-defined one- and two-forms on the corresponding chains.
This justifies the use of the integral notation in the duality relation~\eqref{eq:dec_chain_basis_cochain_basis}.

\subsection{Hodge Operator}\label{sec:dec_hodge}

The discrete hodge operator $\star_{h}$ is a bijection which takes discrete forms on the primal grid to the dual grid as follows,
\begin{subequations}\label{eq:dec_hodge}
\begin{align}
\hodge_{h} \obar{v}_{i+1/2, \, j+1/2} &= \obar{c}^{\hodge}_{i+1/2, \, j+1/2} , &
\hodge_{h} \obar{c}_{i, \, j}         &= \obar{v}^{\hodge}_{i, \, j} , \\
\hodge_{h} \obar{e}^{x}_{i, \, j+1/2} &= \obar{e}^{\hodge y}_{i, \, j+1/2} , &
\hodge_{h} \obar{e}^{y}_{i+1/2, \, j} &= - \obar{e}^{\hodge x}_{i+1/2, \, j} .
\end{align}
\end{subequations}
The inverse of the hodge $\hodge_{h}$, taking discrete forms on the dual grid to the primal grid, is also denoted $\hodge_{h}$ and we have
\begin{subequations}
\begin{align}
\hodge_{h} \obar{c}^{\hodge}_{i+1/2, \, j+1/2} &= \obar{v}_{i+1/2, \, j+1/2} , &
\hodge_{h} \obar{v}^{\hodge}_{i, \, j}         &= \obar{c}_{i, \, j} , \\
\hodge_{h} \obar{e}^{\hodge y}_{i, \, j+1/2}   &= - \obar{e}^{x}_{i, \, j+1/2} , &
\hodge_{h} \obar{e}^{\hodge x}_{i+1/2, \, j}   &= \obar{e}^{y}_{i+1/2, \, j} .
\end{align}
\end{subequations}
Applying the discrete hodge twice, we find that $\hodge_{h} \hodge_{h} \phi_{h} = \phi_{h}$, $\hodge_{h} \hodge_{h} \alpha_{h} = - \alpha_{h}$ and $\hodge_{h} \hodge_{h} \omega_{h} = \omega_{h}$, for $\phi_{h} \in \Lambda_{h}^{0}$, $\alpha_{h} \in \Lambda_{h}^{1}$ and $\omega_{h} \in \Lambda_{h}^{2}$, respectively.
This shows that the identity $\hodge \hodge \alpha = (-1)^{p(d-p)} \alpha$ is preserved at the discrete level.

\subsection{Exterior Derivative}

The exterior derivative of a zero-form $\phi$ is given by
\begin{align}
\ext \phi = \partial_{x} \phi \, \ext x + \partial_{y} \phi \, \ext y ,
\end{align}
and discretised by
\begin{align}\label{eq:exterior_derivative_discrete_one_form}
\ext_{h} \phi_{h}
&= \sum \limits_{i,j} \Big[ 
   (\Delta_{x} \phi)_{i, \, j+1/2} \, \obar{e}^{x}_{i, \, j+1/2}
 + (\Delta_{y} \phi)_{i+1/2, \, j} \, \obar{e}^{y}_{i+1/2, \, j}
\Big] ,
\end{align}
with
\begin{subequations}
\begin{align}
(\Delta_{x} \phi)_{i, \, j+1/2} &= \dfrac{\phi_{i+1/2, \, j+1/2} - \phi_{i-1/2, \, j+1/2}}{h_{x}} , \\
(\Delta_{y} \phi)_{i+1/2, \, j} &= \dfrac{\phi_{i+1/2, \, j+1/2} - \phi_{i+1/2, \, j-1/2}}{h_{y}} .
\end{align}
\end{subequations}
As a consequence of this definition, we can state a discrete version of Stokes' theorem, namely,
\begin{align}\label{eq:discrete_stokes_zero_forms}
\int \limits_{e_{h}} \ext_{h} \phi_{h} = \int \limits_{\partial e_{h}} \phi_{h} .
\end{align}
In order to prove this relation, it suffices to consider the chain basis elements $e^{x}_{i, \, j+1/2}$ and $e^{y}_{i+1/2, \, j}$, for which we have
\begin{subequations}
\begin{align}
\partial e^{x}_{i, \, j+1/2} &= v_{i+1/2, \, j+1/2} - v_{i-1/2, \, j+1/2} , \\
\partial e^{y}_{i+1/2, \, j} &= v_{i+1/2, \, j+1/2} - v_{i+1/2, \, j-1/2} ,
\end{align}
\end{subequations}
and thus
\begin{subequations}
\begin{align}
\int \limits_{\partial e^{x}_{i, \, j+1/2}} \phi_{h} 
&= \phi_{i+1/2, \, j+1/2} - \phi_{i-1/2, \, j+1/2} , \\
\int \limits_{\partial e^{y}_{i+1/2, \, j}} \phi_{h} 
&= \phi_{i+1/2, \, j+1/2} - \phi_{i+1/2, \, j-1/2} .
\end{align}
\end{subequations}
At the same time, integrating \eqref{eq:exterior_derivative_discrete_one_form} over $e^{x}_{i, \, j+1/2}$ and $e^{y}_{i+1/2, \, j}$, respectively, we get
\begin{subequations}
\begin{align}
\int \limits_{e^{x}_{i, \, j+1/2}} \ext_{h} \phi_{h}
&= (\Delta_{x} \phi)_{i, \, j+1/2} \, \int \limits_{e^{x}_{i, \, j+1/2}} \obar{e}^{x}_{i, \, j+1/2}
 = \phi_{i+1/2, \, j+1/2} - \phi_{i-1/2, \, j+1/2} , \\
\int \limits_{e^{y}_{i+1/2, \, j}} \ext_{h} \phi_{h}
&= (\Delta_{y} \phi)_{i+1/2, \, j} \, \int \limits_{e^{y}_{i+1/2, \, j}} \obar{e}^{y}_{i+1/2, \, j}
 = \phi_{i+1/2, \, j+1/2} - \phi_{i+1/2, \, j-1/2} ,
\end{align}
\end{subequations}
which verifies~\eqref{eq:discrete_stokes_zero_forms}.

The exterior derivative of a one-form $\alpha$ is given by
\begin{align}
\ext \alpha = \big[ \partial_{x} \alpha^{y} - \partial_{y} \alpha^{x} \big] \ext x \wedge \ext y ,
\end{align}
and discretised by
\begin{align}
\ext_{h} \alpha_{h}
&= \sum \limits_{i,j} \Big[ (\Delta_{x} \alpha^{y})_{i, \, j} - (\Delta_{y} \alpha^{x})_{i, \, j} \Big] \, \obar{c}_{i, \, j} ,
\end{align}
with
\begin{subequations}
\begin{align}
(\Delta_{x} \alpha^{y})_{i, \, j} &= \dfrac{\alpha^{y}_{i+1/2, \, j} - \alpha^{y}_{i-1/2, \, j}}{h_{x}} , \\
(\Delta_{y} \alpha^{x})_{i, \, j} &= \dfrac{\alpha^{x}_{i, \, j+1/2} - \alpha^{x}_{i, \, j-1/2}}{h_{y}} ,
\end{align}
\end{subequations}
which also preserves Stokes' theorem.

\subsubsection*{Exactness}

Computing the discrete exterior derivative of $\ext_{h} \phi_{h}$ with periodic boundary conditions, we find
\begin{align}
\ext_{h} \ext_{h} \phi_{h}
\nonumber
&= \sum \limits_{i,j} \dfrac{1}{2 h_{y}} \Big[ (\Delta_{x} \phi)_{i, \, j+1/2} - (\Delta_{x} \phi)_{i, \, j-1/2} \Big] \, \obar{c}_{i, \, j} \\
&- \sum \limits_{i,j} \dfrac{1}{2 h_{x}} \Big[ (\Delta_{y} \phi)_{i+1/2, \, j} - (\Delta_{y} \phi)_{i-1/2, \, j} \Big] \, \obar{c}_{i, \, j} = 0 ,
\end{align}
which means that the result of the exterior derivative applied twice vanishes as it should.

\subsubsection*{Differential Operators}

Following the previous exposition, the divergence of a one-form $\alpha$, computed by
\begin{align}
\hodge \ext \hodge \alpha = \big[ \partial_{x} \alpha^{x} + \partial_{y} \alpha^{y} \big] ,
\end{align}
is discretised by
\begin{align}
\hodge_{h} \ext_{h} \hodge_{h} \alpha_{h}
\nonumber
&= \hodge_{h} \ext_{h} \hodge_{h} \sum \limits_{i,j} \Big[ \alpha^{x}_{i, \, j+1/2} \, \obar{e}^{x}_{i, \, j+1/2} + \alpha^{y}_{i+1/2, \, j} \, \obar{e}^{y}_{i+1/2, \, j} \Big] \\
\nonumber
&= \hodge_{h} \ext_{h} \sum \limits_{i,j} \Big[ \alpha^{x}_{i, \, j+1/2} \, \obar{e}^{\hodge y}_{i, \, j+1/2} + \alpha^{y}_{i+1/2, \, j} \, \obar{e}^{\hodge x}_{i+1/2, \, j} \Big] \\
\nonumber
&= \hodge_{h} \sum \limits_{i,j} \bigg[ \dfrac{\alpha^{x}_{i+1, \, j+1/2} - \alpha^{x}_{i, \, j+1/2}}{h_{x}} 
        - \dfrac{\alpha^{y}_{i+1/2, \, j+1} - \alpha^{y}_{i+1/2, \, j}}{h_{y}}
   \bigg] \, \obar{c}^{\hodge}_{i+1/2, \, j+1/2} \\
&= \sum \limits_{i,j} \bigg[ \dfrac{\alpha^{x}_{i+1, \, j+1/2} - \alpha^{x}_{i, \, j+1/2}}{h_{x}} 
        - \dfrac{\alpha^{y}_{i+1/2, \, j+1} - \alpha^{y}_{i+1/2, \, j}}{h_{y}}
   \bigg] \, \obar{v}_{i+1/2, \, j+1/2} .
\end{align}
The Laplacian of a zero-form, computed by
\begin{align}
\hodge \ext \hodge \ext \phi 
= \hodge \ext \big[ \partial_{x} \alpha^{x} \ext y - \partial_{y} \alpha^{y} \ext x \big]
= \big[ \partial_{x} \partial_{x} \alpha^{x} + \partial_{y} \partial_{y} \alpha^{y} \big] ,
\end{align}
is discretised by
\begin{align}
\hodge_{h} \ext_{h} \hodge_{h} \ext_{h} \phi_{h}
\nonumber
&= \hodge_{h} \ext_{h} \sum \limits_{i,j} \Big[ 
   (\Delta_{x} \phi)_{i, \, j+1/2} \, \obar{e}^{\hodge y}_{i, \, j+1/2}
 + (\Delta_{y} \phi)_{i+1/2, \, j} \, \obar{e}^{\hodge x}_{i+1/2, \, j}
\Big] \\
\nonumber
&= \sum \limits_{i,j} \Big[ 
   \dfrac{(\Delta_{x} \phi)_{i+1, \, j+1/2} - (\Delta_{x} \phi)_{i, \, j+1/2}}{h_{x}} \\
&\hspace{10em}
 + \dfrac{(\Delta_{y} \phi)_{i+1/2, \, j+1} - (\Delta_{y} \phi)_{i+1/2, \, j}}{h_{y}}
\Big] \, \obar{v}_{i+1/2, \, j+1/2} ,
\end{align}
which happens to be the usual second-order finite difference Laplacian.

\subsection{Exterior Product}\label{sec:exterior_product}

In order to discretise the exterior product of two forms, we have to introduce appropriate averaging as (except for the exterior product of zero-forms) the result will be a form of higher degree and therefore located on a different position on the grid. For each combination of forms, the particular averaging follows from the definition of the exterior product of the corresponding basis forms.

We need to define four different discrete exterior products: primal-primal, dual-dual, primal-dual and dual-primal. While the definition of the primal-primal and dual-dual exterior products are practically identical, the primal-dual and dual-primal exterior products require separate treatment.
In the first case, the result is a form on the primal grid, while the the second case, the result is a form on the dual grid.
This apparently arbitrary choice will become clear in Section~\ref{sec:dec_pairings} in relation to the definition of discrete pairings.

\subsubsection*{Primal-Primal and Dual-Dual Exterior Product}

The exterior product of two zero-forms $\phi$ and $\psi$ is a zero-form defined as
\begin{align}
\phi \wedge \psi
&= \phi \, \psi ,
\end{align}
and discretised point-wise as
\begin{align}\label{eq:dec_wedge_primal_one_primal_one}
\phi_{h} \wedge_{h} \psi_{h}
&= \sum \limits_{i,j} \phi_{i+1/2, \, j+1/2} \, \psi_{i+1/2, \, j+1/2} \, \obar{v}_{i+1/2, \, j+1/2} .
\end{align}
This can be expressed with respect to the basis forms as
\begin{align}
\obar{v}_{i+1/2, \, j+1/2} \wedge_{h} \obar{v}_{k+1/2, \, l+1/2} &= \delta_{ik} \delta_{jl} \, \obar{v}_{i+1/2, \, j+1/2} .
\end{align}
The exterior product of a zero-form $\phi$ and a one-form $\alpha$ is a one-form defined as
\begin{align}
\phi \wedge \alpha
&= \phi \, \alpha^{x} \ext x + \phi \, \alpha^{y} \ext y ,
\end{align}
and discretised as
\begin{align}
\phi_{h} \wedge_{h} \alpha_{h}
\nonumber
= \tfrac{1}{2} \sum \limits_{i,j} \Big[
  & \big( \phi_{i-1/2, \, j+1/2} + \phi_{i+1/2, \, j+1/2} \big) \, \alpha^{x}_{i, \, j+1/2} \, \obar{e}^{x}_{i, \, j+1/2} \\
 +& \big( \phi_{i+1/2, \, j-1/2} + \phi_{i+1/2, \, j+1/2} \big) \, \alpha^{y}_{i+1/2, \, j} \, \obar{e}^{y}_{i+1/2, \, j}
\Big] ,
\end{align}
which follows from the following relation of the basis forms,
\begin{subequations}
\begin{align}
\obar{v}_{i+1/2, \, j+1/2} \wedge_{h} \obar{e}^{x}_{k, \, l+1/2} &= \tfrac{1}{2} \delta_{ik} \delta_{jl} \, \big[ \obar{e}^{x}_{k, \, l+1/2} + \obar{e}^{x}_{k+1, \, l+1/2} \big] , \\
\obar{v}_{i+1/2, \, j+1/2} \wedge_{h} \obar{e}^{y}_{k+1/2, \, l} &= \tfrac{1}{2} \delta_{ik} \delta_{jl} \, \big[ \obar{e}^{y}_{k+1/2, \, l} + \obar{e}^{y}_{k+1/2, \, l+1} \big] .
\end{align}
\end{subequations}
As the one-form $\phi_{h}$ is defined on vertices while the one-form $\alpha_{h}$ and the resulting one-form are defined on edges, we need to average $\phi_{h}$ accordingly.
The exterior product of a zero-form $\phi$ and a two-form $\omega$ is a two-form defined as
\begin{align}
\phi \wedge \omega
&= \phi \, \omega \, \ext x \wedge \ext y ,
\end{align}
and discretised as
\begin{align}
\phi_{h} \wedge_{h} \omega_{h}
&= \tfrac{1}{4} \sum \limits_{i,j} \Big[
   \phi_{i-1/2, \, j-1/2}
 + \phi_{i+1/2, \, j-1/2}
 + \phi_{i-1/2, \, j+1/2}
 + \phi_{i+1/2, \, j+1/2}
\Big] \, \omega_{i, \, j} \, \obar{c}_{i, \, j} ,
\end{align}
where we computed cell averages of the zero-forms, which are defined on the vertices.
This relation follows from
\begin{align}
\obar{v}_{i+1/2, \, j+1/2} \wedge_{h} \obar{c}_{k, \, l} = \tfrac{1}{4} \delta_{ik} \delta_{jl} \, \big[ \obar{c}_{i,j} + \obar{c}_{i,j+1} + \obar{c}_{i+1,j} + \obar{c}_{i+1,j+1} \big] .
\end{align}
The exterior product of two one-forms $\alpha$ and $\beta$ is a two-form defined as
\begin{align}
\alpha \wedge \beta
&= \big[ \alpha^{x} \ext x + \alpha^{y} \ext y \big] \wedge \big[ \beta^{x} \ext x + \beta^{y} \ext y \big]
 = \big[ \alpha^{x} \beta^{y} - \alpha^{y} \beta^{x} \big] \ext x \wedge \ext y .
\end{align}
While the one-forms are defined on edges, the resulting two-form is defined on cells and therefore discretised as
\begin{align}
\alpha_{h} \wedge_{h} \beta_{h}
\nonumber
 = \tfrac{1}{4} \sum \limits_{i,j} \Big[
 & \big( \alpha^{x}_{i, \, j-1/2} + \alpha^{x}_{i, \, j+1/2} \big) \, 
   \big( \beta^{y}_{i-1/2, \, j} + \beta^{y}_{i+1/2, \, j} \big) \\
-& \big( \alpha^{y}_{i-1/2, \, j} + \alpha^{y}_{i+1/2, \, j} \big) 
   \big( \beta^{x}_{i, \, j-1/2} + \beta^{x}_{i, \, j+1/2} \big)
\Big] \, \obar{c}_{i, \, j} ,
\end{align}
where the $x$-components are averaged in $y$-direction and the $y$-components are averaged in $x$-direction, as follows from
\begin{align}
  \obar{e}^{x}_{i, \, j+1/2} \wedge_{h} \obar{e}^{y}_{k+1/2, \, l} &= \tfrac{1}{4} \delta_{ik} \delta_{jl} \, \big[ \obar{c}_{i,j} + \obar{c}_{i,j+1} + \obar{c}_{i+1,j} + \obar{c}_{i+1,j+1} \big] .
\end{align}
The dual-dual exterior products follow in complete analogy to the primal-primal definitions. The primal-dual and dual-primal exterior products, however, need special treatment.

\subsubsection*{Primal-Dual and Dual-Primal Exterior Products}

With two separate meshes for primal and dual forms, we have different exterior products between primal and dual forms and dual and primal forms, respectively.
For the primal-dual and dual-primal exterior products, we only consider those combinations which lead to a volume form.
As before, the exterior products of forms are defined through the exterior products of the basis forms.
For zero- and two-forms, we have
\begin{align}
\obar{v}_{i+1/2, \, j+1/2} \wedge_{h} \obar{c}^{\hodge}_{k+1/2, \, l+1/2} &= \tfrac{1}{4} \delta_{ik} \delta_{jl} \, \big[ \obar{c}_{i, \, j} + \obar{c}_{i, \, j+1} + \obar{c}_{i+1, \, j} + \obar{c}_{i+1, \, j+1} \big] ,
\end{align}
and
\begin{align}
\obar{v}^{\hodge}_{i, \, j} \wedge_{h} \obar{c}_{k, \, l} &= \tfrac{1}{4} \delta_{ik} \delta_{jl} \, \big[ \obar{c}^{\hodge}_{i-1/2, \, j-1/2} + \obar{c}^{\hodge}_{i-1/2, \, j+1/2} + \obar{c}^{\hodge}_{i+1/2, \, j-1/2} + \obar{c}^{\hodge}_{i+1/2, \, j+1/2} \big] ,
\end{align}
so that the exterior product of a primal zero-form $\psi_{h}$ and a dual two-form $\omega_{h}^{\hodge}$ becomes
\begin{multline}\label{eq:dec_wedge_primal_one_form_dual_two_form}
\psi_{h} \wedge_{h} \omega_{h}^{\hodge} = \tfrac{1}{4} \sum \limits_{i,j} \Big[ 
  \psi_{i-1/2, \, j-1/2} \, \omega^{\hodge}_{i-1/2, \, j-1/2}
+ \psi_{i-1/2, \, j+1/2} \, \omega^{\hodge}_{i-1/2, \, j+1/2} \\
+ \psi_{i+1/2, \, j-1/2} \, \omega^{\hodge}_{i+1/2, \, j-1/2}
+ \psi_{i+1/2, \, j+1/2} \, \omega^{\hodge}_{i+1/2, \, j+1/2} \Big] \, \obar{c}_{i, \, j} ,
\end{multline}
and the exterior product of a dual zero-form $\psi_{h}^{\hodge}$ and a primal two-form $\omega_{h}$ becomes
\begin{multline}
\psi_{h}^{\hodge} \wedge_{h} \omega_{h} = \tfrac{1}{4} \sum \limits_{i,j} \Big[
  \psi^{\hodge}_{i,   \, j  } \, \omega_{i,   \, j  } 
+ \psi^{\hodge}_{i,   \, j+1} \, \omega_{i,   \, j+1} \\
+ \psi^{\hodge}_{i+1, \, j  } \, \omega_{i+1, \, j  } 
+ \psi^{\hodge}_{i+1, \, j+1} \, \omega_{i+1, \, j+1} 
\Big] \, \obar{c}^{\hodge}_{i+1/2, \, j+1/2} .
\end{multline}
For exterior products of basis one-forms, we have
\begin{subequations}
\begin{align}
  \obar{e}^{x}_{i, \, j+1/2} \wedge_{h} \obar{e}^{\hodge y}_{k, \, l+1/2} &= + \tfrac{1}{2} \delta_{ik} \delta_{jl} \, \big[ \obar{c}_{i,j} + \obar{c}_{i,j+1} \big] , \\
  \obar{e}^{y}_{i+1/2, \, j} \wedge_{h} \obar{e}^{\hodge x}_{k+1/2, \, l} &= - \tfrac{1}{2} \delta_{ik} \delta_{jl} \, \big[ \obar{c}_{i,j} + \obar{c}_{i+1,j} \big] ,
\end{align}
\end{subequations}
and
\begin{subequations}
\begin{align}
  \obar{e}^{\hodge x}_{i+1/2, \, j} \wedge_{h} \obar{e}^{y}_{k+1/2, \, l} &= + \tfrac{1}{2} \delta_{ik} \delta_{jl} \, \big[ \obar{c}^{\hodge}_{i+1/2,j-1/2} + \obar{c}^{\hodge}_{i+1/2,j+1/2} \big] , \\
  \obar{e}^{\hodge y}_{i, \, j+1/2} \wedge_{h} \obar{e}^{x}_{k, \, l+1/2} &= - \tfrac{1}{2} \delta_{ik} \delta_{jl} \, \big[ \obar{c}^{\hodge}_{i-1/2,j+1/2} + \obar{c}^{\hodge}_{i+1/2,j+1/2} \big] ,
\end{align}
\end{subequations}
so that the exterior products between primal and dual one-forms, $\alpha_{h}$ and $\beta_{h}$ as well as $\alpha_{h}^{\hodge}$ and $\beta_{h}^{\hodge}$, become
\begin{align}
\label{eq:wedge_primal_dual_one_forms}
\nonumber
\alpha_{h} \wedge_{h} \beta_{h}^{\hodge}
 = \tfrac{1}{2} \sum \limits_{i,j} \big[
    & \alpha^{x}_{i, \, j-1/2} \, \beta^{\hodge y}_{i, \, j-1/2}
    + \alpha^{x}_{i, \, j+1/2} \, \beta^{\hodge y}_{i, \, j+1/2} \\
    & \hfill
    - \alpha^{y}_{i-1/2, \, j} \, \beta^{\hodge x}_{i-1/2, \, j}
    - \alpha^{y}_{i+1/2, \, j} \, \beta^{\hodge x}_{i+1/2, \, j}
\big] \, \obar{c}_{i, \, j} , \\
\label{eq:wedge_dual_primal_one_forms}
\nonumber
\alpha_{h}^{\hodge} \wedge_{h} \beta_{h}
 = \tfrac{1}{2} \sum \limits_{i,j} \big[
    & \alpha^{\hodge x}_{i-1/2, \, j} \, \beta^{y}_{i-1/2, \, j}
    + \alpha^{\hodge x}_{i+1/2, \, j} \, \beta^{y}_{i+1/2, \, j} \\
    & \hfill
    - \alpha^{\hodge y}_{i, \, j-1/2} \, \beta^{x}_{i, \, j-1/2}
    - \alpha^{\hodge y}_{i, \, j+1/2} \, \beta^{x}_{i, \, j+1/2}
\big] \, \obar{c}^{\hodge}_{i+1/2,j+1/2} .
\end{align}
For exterior products of basis two-forms and basis zero-forms, we have
\begin{align}
\obar{c}_{i, \, j} \wedge_{h} \obar{v}^{\hodge}_{k, \, l} &= \delta_{ik} \delta_{jl} \, \obar{c}_{i, \, j} ,
\end{align}
and
\begin{align}
\obar{c}^{\hodge}_{i+1/2, \, j+1/2} \wedge_{h} \obar{v}_{k+1/2, \, l+1/2}  &= \delta_{ik} \delta_{jl} \, \obar{c}^{\hodge}_{i+1/2, \, j+1/2} ,
\end{align}
so that
\begin{align}
\omega_{h} \wedge_{h} \phi_{h}^{\hodge}
&= \sum \limits_{i,j} \omega_{i, \, j} \, \phi^{\hodge}_{i, \, j} \, \obar{c}_{i, \, j} , \\
\label{eq:dec_wedge_dual_two_primal_zero}
\omega_{h}^{\hodge} \wedge_{h} \phi_{h}
&= \sum \limits_{i,j} \omega^{\hodge}_{i+1/2, \, j+1/2} \, \phi_{i+1/2, \, j+1/2} \, \obar{c}^{\hodge}_{i+1/2, \, j+1/2} .
\end{align}
With these definitions, we can define the pairings necessary to discretise the Lagrangian.

\subsection{Pairings}\label{sec:dec_pairings}

The discrete version of the pairing~\eqref{eq:forms_inner_product} follows on noting that $\alpha_{h} \wedge_{h} \hodge_{h} \beta_{h}$ is a discrete two-form for arbitrary discrete forms $\alpha_{h}, \beta_{h} \in \Lambda_{h}^{p}$ with $p \in \{ 0, 1, 2 \}$ and can be integrated over cell chains.
Analogously, for any $\alpha_{h}^{\hodge}, \beta_{h}^{\hodge} \in \widetilde{\Lambda}_{h}^{p}$ the expression $\alpha_{h}^{\hodge} \wedge_{h} \hodge_{h} \beta_{h}^{\hodge}$ is a discrete two-form on the dual grid.
This is a consequence of the specific choices in the definition of the discrete exterior product in Section~\ref{sec:exterior_product}.
Upon identifying the computational domain $\Omega$ with primary and dual chains with unit coefficients, 
\begin{align}
\Omega_{h} &= \sum \limits_{i,j} c_{i, \, j} &
& \text{and} &
\Omega_{h}^{\hodge} &= \sum \limits_{i,j} c_{i+1/2, \, j+1/2}^{\hodge} ,
\end{align}
respectively, it is natural to define the discrete pairing by
\begin{align}\label{eq:dec_wedge_primal_zero_primal_zero}
\bracket{ \alpha_{h} , \beta_{h} }_{\Omega_{h}}
&= \int \limits_{\Omega_{h}} \alpha_{h} \wedge_{h} \hodge_{h} \beta_{h} &
& \text{and} &
\bracket{ \alpha_{h}^{\hodge} , \beta_{h}^{\hodge} }_{\Omega_{h}}
&= \int \limits_{\Omega_{h}} \alpha_{h}^{\hodge} \wedge_{h} \hodge_{h} \beta_{h}^{\hodge} .
\end{align}
From Equations~\eqref{eq:dec_chain_basis_cochain_basis}, \eqref{eq:dec_hodge} and~\eqref{eq:dec_wedge_primal_one_form_dual_two_form} we obtain explicit expressions for the pairing of discrete zero-forms $\gamma_{h}$ and $\phi_{h}$,
\begingroup
\allowdisplaybreaks
\begin{align}\label{eq:dec_pairing_primal_zero_forms}
\bracket{ \gamma_{h} , \phi_{h} }_{\Omega_{h}}
\nonumber
&= \int \limits_{\Omega_{h}} \sum \limits_{i,j} \dfrac{1}{4} \Big[
  \gamma_{i-1/2, \, j-1/2} \, \phi_{i-1/2, \, j-1/2} + \gamma_{i-1/2, \, j+1/2} \, \phi_{i-1/2, \, j+1/2} \\
\nonumber
&\hspace{6em}
 + \gamma_{i+1/2, \, j-1/2} \, \phi_{i+1/2, \, j-1/2} + \gamma_{i+1/2, \, j+1/2} \, \phi_{i+1/2, \, j+1/2} 
\Big] \,  \obar{c}_{i, \, j} \\
\nonumber
&= \dfrac{h_x h_y}{4} \sum \limits_{i,j} \Big[
  \gamma_{i-1/2, \, j-1/2} \, \phi_{i-1/2, \, j-1/2} + \gamma_{i-1/2, \, j+1/2} \, \phi_{i-1/2, \, j+1/2} \\
&\hspace{6em}
 + \gamma_{i+1/2, \, j-1/2} \, \phi_{i+1/2, \, j-1/2} + \gamma_{i+1/2, \, j+1/2} \, \phi_{i+1/2, \, j+1/2} 
\Big] ,
\end{align}
\endgroup
the pairing of discrete one-forms $\alpha_{h}$ and $\beta_{h}$,
\begin{align}\label{eq:dec_pairing_primal_one_forms}
\bracket{ \alpha_{h} , \beta_{h} }_{\Omega_{h}}
\nonumber
&= \int \limits_{\Omega_{h}} \sum \limits_{i,j} \dfrac{1}{2} \Big[ 
   \alpha^{x}_{ i, \, j-1/2 } \, \beta^{x}_{i, \, j-1/2 }
 + \alpha^{x}_{ i, \, j+1/2 } \, \beta^{x}_{i, \, j+1/2 } \\
\nonumber
&\hspace{6em}
 + \alpha^{y}_{ i-1/2, \, j } \, \beta^{y}_{i-1/2, \, j }
 + \alpha^{y}_{ i+1/2, \, j } \, \beta^{y}_{i+1/2, \, j } 
\Big] \, \obar{c}_{i, \, j} \\
\nonumber
&= \dfrac{h_x h_y}{2} \sum \limits_{i,j} \Big[ 
   \alpha^{x}_{ i, \, j-1/2 } \, \beta^{x}_{i, \, j-1/2 }
 + \alpha^{x}_{ i, \, j+1/2 } \, \beta^{x}_{i, \, j+1/2 } \\
&\hspace{6em}
 + \alpha^{y}_{ i-1/2, \, j } \, \beta^{y}_{i-1/2, \, j }
 + \alpha^{y}_{ i+1/2, \, j } \, \beta^{y}_{i+1/2, \, j } 
\Big] ,
\end{align}
and the pairing of discrete two-forms $\sigma_{h}$ and $\omega_{h}$,
\begin{align}\label{eq:dec_pairing_primal_two_forms}
\bracket{ \sigma_{h} , \omega_{h} }_{\Omega_{h}}
&= h_x h_y \sum \limits_{i,j} \sigma_{i, \, j} \, \omega_{i, \, j} .
\end{align}
For the convenience of the reader, we provide explicit expressions of the pairings on the dual grid, namely for discrete zero-forms $\gamma_{h}$ and $\phi_h$, 
\begin{align}\label{eq:dec_pairing_dual_zero_forms}
\bracket{ \hodge_{h} \gamma_{h} , \hodge_{h} \phi_{h} }_{\Omega_{h}^{\hodge}}
&= h_x h_y \sum \limits_{i=1}^{n_{x}} \sum \limits_{j=1}^{n_{y}} \gamma_{i+1/2, \, j+1/2} \, \phi_{i+1/2, \, j+1/2} ,
\end{align}
for discrete one-forms $\alpha_{h}$ and $\beta_{h}$,
\begin{multline}\label{eq:dec_pairing_dual_one_forms}
\bracket{ \hodge_{h} \alpha_{h} , \hodge_{h} \beta_{h} }_{\Omega_{h}^{\hodge}}
= \dfrac{h_x h_y}{2} \sum \limits_{i,j} \Big[ 
  \alpha^{x}_{ i,   \, j+1/2 } \, \beta^{x}_{i,   \, j+1/2 }
+ \alpha^{x}_{ i+1, \, j+1/2 } \, \beta^{x}_{i+1, \, j+1/2 } \\
+ \alpha^{y}_{ i+1/2, \, j   } \, \beta^{y}_{i+1/2, \, j   } 
+ \alpha^{y}_{ i+1/2, \, j+1 } \, \beta^{y}_{i+1/2, \, j+1 } 
\Big] ,
\end{multline}
and for discrete two-forms $\sigma_{h}$ and $\omega_{h}$
\begin{multline}\label{eq:dec_pairing_dual_two_forms}
\bracket{ \hodge_{h} \sigma_{h} , \hodge_{h} \omega_{h} }_{\Omega_{h}^{\hodge}}
= \dfrac{h_x h_y}{4} \sum \limits_{i,j} \Big[
  \sigma_{i, \, j} \, \omega_{i, \, j} + \sigma_{i+1, \, j} \, \omega_{i+1, \, j} \\
+ \sigma_{i, \, j+1} \, \omega_{i, \, j+1} + \sigma_{i+1, \, j+1} \, \omega_{i+1, \, j+1} 
\Big] .
\end{multline}

At the continuous level, the Hodge operator $\hodge$ is an isometry with respect to the pairing, namely for $\alpha, \beta \in \Lambda^{p} (\Omega)$ we have that
\begin{align}
\bracket{\hodge \alpha , \hodge \beta}_{\Omega}
= \bracket{\alpha , \beta}_{\Omega} .
\end{align}
As a consequence of the definitions of discrete exterior products it can be seen that this property is preserved at the discrete level, in particular $\hodge_{h}$ is an isometry between the primal and dual grid,
\begin{align}\label{eq:dec_hodge_isometry}
\bracket{ \hodge_{h} \alpha_{h} , \hodge_{h} \beta_{h} }_{\Omega_{h}^{\hodge}}
= \bracket{ \alpha_{h} , \beta_{h} }_{\Omega_{h}} .
\end{align}
This follows from
\begin{align}\label{eq:dec_pairing_isometry}
\int \limits_{\Omega_{h}^{\hodge}} \hodge_{h} \alpha_{h} \wedge_{h} \beta_{h} = (-1)^{p (d-p)} \int \limits_{\Omega_{h}} \alpha_{h} \wedge_{h} \hodge_{h} \beta_{h} ,
\end{align}
which can be proved directly by inspection of the expressions for the discrete exterior products in Section~\ref{sec:exterior_product}, together with $\hodge_{h} \hodge_{h} \alpha = (-1)^{p (d-p)} \alpha$, cf. Section~\ref{sec:dec_hodge}.
The identity~\eqref{eq:dec_hodge_isometry} can also be checked directly by comparing Equations~\eqref{eq:dec_pairing_dual_zero_forms}-\eqref{eq:dec_pairing_dual_two_forms} to Equations~\eqref{eq:dec_pairing_primal_zero_forms}-\eqref{eq:dec_pairing_primal_two_forms} and accounting for periodic boundary conditions.

\section{Variational Discretisation}\label{sec:vi}

In order to obtain a numerical method for the ideal MHD equations (\ref{eq:ideal_mhd_equations}), we discretise the action functional $\mcal{A}$ and apply a discrete version of Hamilton's principle of stationary action \cite{MarsdenPatrick:1998, KrausMaj:2015}.
Standard variational discretisations on cartesian meshes usually lead to centred finite difference schemes, which are problematic for the momentum equation.
Such schemes, for which the components of the velocity vector and the pressure are located at the same grid points, are known to be prone to instabilities (see e.g.~\citet{Langtangen:2002} or~\citet{McDonough:2007}).
The pressure often becomes highly oscillatory as a symmetric difference operator, e.g., with stencil $[ -1 \; \hphantom{-}0 \; +1 ]$, annihilates pressures which oscillate between~$+1$ and~$-1$ between neighbouring grid points. This is often referred to as \emph{checker-boarding}.

An effective remedy for this problem within the finite difference method is the introduction of a staggered grid. The notion of discrete exterior calculus as introduced in the previous section provides a guiding principle on how to locate the various quantities on the grid so that their geometric character and desirable properties like the identities from vector calculus are retained on the discrete level.
Using these notions, the discrete divergence-free constraint of the velocity field, $\hodge_{h} \ext_{h} \hodge_{h} V_{h} = 0$, becomes
\begin{align}\label{eq:mhd_staggered_grid_divergence}
  \dfrac{V^{x}_{i, \, j+1/2, \, n} - V^{x}_{i-1, \, j+1/2, \, n}}{h_{x}}
+ \dfrac{V^{y}_{i+1/2, \, j, \, n} - V^{y}_{i+1/2, \, j-1, \, n}}{h_{y}} = 0 ,
\end{align}
which is defined in such a way that the natural location of the divergence coincides with the location of the pressure.
The function of the pressure in incompressible fluid dynamics can be described as a Lagrange multiplier enforcing the divergence-free constraint of the velocity field.
By the above discretisation, only the pressure at a single grid point enforces the divergence of the velocity of the surrounding grid points to vanish. As the divergence is computed by one-sided finite differences (and not e.g. by a combination of forward and backward differences), checker-boarding will not be an issue.

\subsection{Temporal Discretisation}\label{sec:temporal_discretisation}

The discrete exterior calculus defined in the previous section fully determines the spatial discretisation. In addition we need to define discrete time derivatives.
The velocity and magnetic field are collocated at full time steps $t_{n} = t_{0} + n h_{t}$, with $h_{t}$ the time step size, and their time derivatives are defined point-wise on the primal grid as
\begin{subequations}
\begin{alignat}{5}
& V^x_t = \partial_{t} V^{x} &
& \quad \rightarrow \quad &
& (\Delta_{t} V^{x})_{i, \, j+1/2, \, n+1/2} &
& \equiv \; &
& \dfrac{V^{x}_{i, \, j+1/2, \, n+1} - V^{x}_{i, \, j+1/2, \, n}}{h_{t}} ,
\\
& V^y_t = \partial_{t} V^{y} &
& \quad \rightarrow \quad &
& (\Delta_{t} V^{y})_{i+1/2, \, j, \, n+1/2} &
& \equiv \; &
& \dfrac{V^{y}_{i+1/2, \, j, \, n+1} - V^{y}_{i+1/2, \, j, \, n}}{h_{t}} ,
\end{alignat}
\end{subequations}
so that the time derivative does not change the differential form character of the variables.
For the pressure $P$, the staggering approach is applied also with respect to time, i.e., the pressure nodes are $(i+1/2, \, j+1/2, \, n+1/2)$. 
The adjoint variables will be collocated at the same spatial grid positions as the corresponding physical variables but staggered with respect to time. That is, $\alpha$ and $\beta$ share their spatial grid positions with  $V$ and $B$, but their temporal position is at half time steps $t_{n+1/2} = t_{0} + (n+1/2) h_{t}$, just as the time derivatives of $V$ and $B$. As $\gamma$ is a scalar field, it is collocated at the same spatial position as the pressure, but its temporal position is at full time steps $t_{n}$.

For concise notation, we introduce averages of the vector fields $V$ and $B$ on the primal grid, where averaging is applied with respect to both space and time,
\begin{subequations}
\begin{align}
\bar{V}^{x}_{i, \, j, \, n+1/2}
&\equiv \dfrac{1}{4} \Big[
  V^{x}_{ i, \, j-1/2 , \, n   } + V^{x}_{ i, \, j+1/2, \, n   }
+ V^{x}_{ i, \, j-1/2 , \, n+1 } + V^{x}_{ i, \, j+1/2, \, n+1 } \Big] ,
\\
\bar{V}^{y}_{i, \, j, \, n+1/2}
&\equiv \dfrac{1}{4} \Big[
  V^{y}_{ i-1/2 , \, j, \, n   } + V^{y}_{ i+1/2, \, j, \, n   }
+ V^{y}_{ i-1/2 , \, j, \, n+1 } + V^{y}_{ i+1/2, \, j, \, n+1 } \Big] .
\end{align}
\end{subequations}
The averages of the adjoint variables $\alpha$ and $\beta$ in~\eqref{eq:ideal_mhd_formal_lagrangian_forms} do not involve time, therefore we have
\begin{subequations}
\begin{align}
\bar{\alpha}^{x}_{i, \, j, \, n+1/2}
&\equiv \dfrac{1}{2} \Big[
 \alpha^{x}_{ i, \, j-1/2 , \, n+1/2 } + \alpha^{x}_{ i, \, j+1/2, \, n+1/2 } \Big] ,
\\
\bar{\alpha}^{y}_{i, \, j, \, n+1/2}
&\equiv \dfrac{1}{2} \Big[
 \alpha^{y}_{ i-1/2 , \, j, \, n+1/2 } + \alpha^{y}_{ i+1/2, \, j, \, n+1/2 } \Big] .
\end{align}
\end{subequations}
With these definitions we will now construct the discrete Lagrangians.

\subsection{Momentum Equation}

We start the derivation of the variational integrator by considering the momentum equation, which corresponds to the incompressible Euler equation with the Lorentz force due to the magnetic field.
The terms of the formal Lagrangian~\eqref{eq:ideal_mhd_formal_lagrangian_forms}, corresponding to the momentum equation, are
\begin{align}\label{eq:mhd_navier_stokes_lagrangian}
L^{\text{M}}
 = \pair{ V_{t} }{ \alpha }
 + \pair{ \ext V }{ V \wedge \alpha } 
 - \pair{ \ext B }{ B \wedge \alpha } 
 + \pair{ \ext P }{ \alpha } .
\end{align}
This expression is discretised on the primal grid as depicted in Figure~\subref*{fig:mhd_staggered_grid_momentum}, yielding the discrete Lagrangian for the momentum equation,
\begin{align}
L^{M}_{n+1/2} = \sum \limits_{i,j} \mcal{L}^{\text{M}}_{i,j,n+1/2} ,
\end{align}
with
\begin{align}
\mcal{L}^{\text{M}}_{i,j,n+1/2}
\nonumber
= h_{t} h_{x} h_{y} \bigg\{
 &   \dfrac{1}{2} \Big[
     \alpha^{x}_{ i, \, j-1/2, \, n+1/2 } \, (\Delta_{t} V^{x})_{i, \, j-1/2, \, n+1/2}
   + \alpha^{x}_{ i, \, j+1/2, \, n+1/2 } \, (\Delta_{t} V^{x})_{i, \, j+1/2, \, n+1/2} 
   \Big] \\
\nonumber
+&{} \dfrac{1}{2} \Big[
      \alpha^{y}_{ i-1/2, \, j, \, n+1/2 } \, (\Delta_{t} V^{y})_{i-1/2, \, j, \, n+1/2}
    + \alpha^{y}_{ i+1/2, \, j, \, n+1/2 } \, (\Delta_{t} V^{y})_{i+1/2, \, j, \, n+1/2} 
   \Big] \\
\nonumber
+&{} \dfrac{1}{2} \Big[
      \alpha^{x}_{ i, \, j-1/2, \, n+1/2 } \, (\Delta_{x} P)_{i, \, j-1/2, \, n+1/2}
    + \alpha^{x}_{ i, \, j+1/2, \, n+1/2 } \, (\Delta_{x} P)_{i, \, j+1/2, \, n+1/2} 
   \Big] \\
\nonumber
+&{} \dfrac{1}{2} \Big[
      \alpha^{y}_{ i-1/2, \, j, \, n+1/2 } \, (\Delta_{x} P)_{i-1/2, \, j, \, n+1/2}
    + \alpha^{y}_{ i+1/2, \, j, \, n+1/2 } \, (\Delta_{y} P)_{i+1/2, \, j, \, n+1/2} 
   \Big] \\
\nonumber
+&{} \bar{\alpha}^{x}_{i, \, j, \, n+1/2} \,
     \bar{V}^{y}_{i, \, j, \, n+1/2} \, \Big[
         (\Delta_{y} V^{x})_{i, \, j, \, n+1/2}
       - (\Delta_{x} V^{y})_{i, \, j, \, n+1/2}
   \Big] \\
\nonumber
+&{} \bar{\alpha}^{y}_{i, \, j, \, n+1/2} \,
     \bar{V}^{x}_{i, \, j, \, n+1/2} \, \Big[
         (\Delta_{x} V^{y})_{i, \, j, \, n+1/2}
       - (\Delta_{y} V^{x})_{i, \, j, \, n+1/2}
   \Big] \\
\nonumber
-&{} \bar{\alpha}^{x}_{i, \, j, \, n+1/2} \, 
     \bar{B}^{y}_{i, \, j, \, n+1/2} \, \Big[
         (\Delta_{y} B^{x})_{i, \, j, \, n+1/2}
       - (\Delta_{x} B^{y})_{i, \, j, \, n+1/2}
   \Big] \\
-&{} \bar{\alpha}^{y}_{i, \, j, \, n+1/2} \, 
     \bar{B}^{x}_{i, \, j, \, n+1/2} \, \Big[
         (\Delta_{x} B^{y})_{i, \, j, \, n+1/2}
       - (\Delta_{y} B^{x})_{i, \, j, \, n+1/2}
   \Big]
\bigg\} .
\end{align}
The part of the formal Lagrangian~\eqref{eq:ideal_mhd_formal_lagrangian_forms}, corresponding to the divergence term,
\begin{align}\label{eq:mhd_divergence_lagrangian}
L^{\text{D}} = \pair{ \hodge \ext \hodge V }{ \gamma } ,
\end{align}
is discretised on the dual grid in Figure~\subref*{fig:mhd_staggered_grid_div} as
\begin{align}
L^{D}_{n+1/2} = \sum \limits_{i,j} \dfrac{1}{2} \Big[ \mcal{L}^{\text{D}}_{i,j,n} + \mcal{L}^{\text{D}}_{i,j,n+1} \Big] ,
\end{align}
with
\begin{align}
\mcal{L}^{\text{D}}_{i,j,n}
&= h_{t} h_{x} h_{y} \,
   \gamma_{ i+1/2, \, j+1/2, \, n } \, \Big[
       (\Delta_{x} V^{x})_{i+1/2, \, j+1/2, \, n}
     + (\Delta_{y} V^{y})_{i+1/2, \, j+1/2, \, n}
\Big] ,
\end{align}
which will give (\ref{eq:mhd_staggered_grid_divergence}) as desired.

\subsection{Induction Equation}

Now we consider those terms of the formal Lagrangian~\eqref{eq:ideal_mhd_formal_lagrangian_forms} that will yield the induction equation,
\begin{align}
L^{\text{I}}
 = \pair{ B_{t} }{ \beta }
 - \pair{ V \wedge B }{ \ext \beta } .
\end{align}
These expressions are again discretised on the primal grid, leading to the discrete Lagrangian for the induction equation,
\begin{align}
L^{I}_{n+1/2} = \sum \limits_{i,j} \mcal{L}^{\text{I}}_{i,j,n+1/2} ,
\end{align}
with
\begin{align}
\mcal{L}^{\text{I}}_{i,j,n+1/2}
\nonumber
&= h_{t} h_{x} h_{y} \bigg\{ \\
\nonumber
&\hspace{1em}\hphantom{+}\;\,
  \dfrac{1}{2} \Big[
  \beta^{x}_{ i, \, j-1/2, \, n+1/2 } \, (\Delta_{t} B^{x})_{i, \, j-1/2, \, n+1/2}
+ \beta^{x}_{ i, \, j+1/2, \, n+1/2 } \, (\Delta_{t} B^{x})_{i, \, j+1/2, \, n+1/2} 
  \Big] \\
\nonumber
&\hspace{1em}
+ \dfrac{1}{2} \Big[
  \beta^{y}_{ i-1/2, \, j, \, n+1/2 } \, (\Delta_{t} B^{y})_{i-1/2, \, j, \, n+1/2}
+ \beta^{y}_{ i+1/2, \, j, \, n+1/2 } \, (\Delta_{t} B^{y})_{i+1/2, \, j, \, n+1/2} 
  \Big] \\
\nonumber
&\hspace{1em}
- (\Delta_{y} \beta^{x})_{i, \, j, \, n+1/2} \, \bigg[
     \bar{V}^{y}_{i, \, j, \, n+1/2} \, \bar{B}^{x}_{i, \, j, \, n+1/2}
   - \bar{V}^{x}_{i, \, j, \, n+1/2} \, \bar{B}^{y}_{i, \, j, \, n+1/2}
\bigg] \\
&\hspace{1em}
- (\Delta_{x} \beta^{y})_{i, \, j, \, n+1/2} \, \bigg[
     \bar{V}^{x}_{i, \, j, \, n+1/2} \, \bar{B}^{y}_{i, \, j, \, n+1/2}
   - \bar{V}^{y}_{i, \, j, \, n+1/2} \, \bar{B}^{x}_{i, \, j, \, n+1/2}
\bigg] \bigg\} .
\end{align}
Now we have all the ingredients for a complete discretisation of the action integral corresponding to~(\ref{eq:ideal_mhd_equations}).

\subsection{Variational Integrator}

\begin{figure}[bt]
\centering
\includegraphics[width=.44\textwidth]{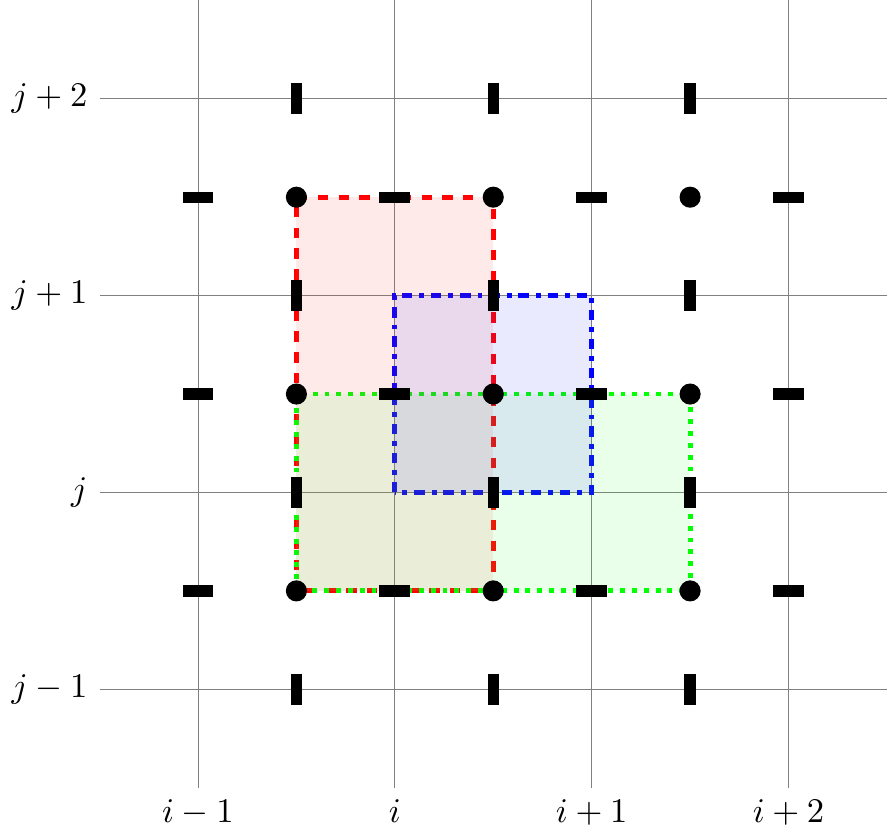}
\caption{Stencils on the staggered grid in the $xy$-plane: $x$-component (red) and $y$-component (green) of momentum and induction equation, divergence constraint (blue).}
\label{fig:mhd_staggered_grid_scheme}
\end{figure}

The discrete action amounts to 
\begin{align}
\mcal{A}_{d} [\phy_{d}] = \delta \sum \limits_{n=0}^{n_{t}-1} \Big[ L^{\text{M}}_{n+1/2} + L^{\text{D}}_{n+1/2} + L^{\text{I}}_{n+1/2} \Big] , 
\end{align}
where $\phy_{d}$ denotes the discrete solution, that is
\begin{align}
\phy_{d} = \Big\{
\nonumber
& V^{x}_{i,j+1/2,n}, \, V^{y}_{i+1/2,j,n}, \, B^{x}_{i,j+1/2,n}, \, B^{y}_{i+1/2,j,n}, \, P_{i+1/2,j+1/2,m+1/2}, \\
\nonumber
& \alpha^{x}_{i,j+1/2,m+1/2}, \, \alpha^{y}_{i+1/2,j,m+1/2}, \, \beta^{x}_{i,j+1/2,m+1/2}, \, \beta^{y}_{i+1/2,j,m+1/2}, \, \gamma_{i+1/2,j+1/2,n} \\
& \Big\vert \; 1 \leq i \leq n_{x} , \, 1 \leq j \leq n_{y} , \, 0 \leq m < n_{t} - 1 , \, 0 \leq n \leq n_{t} \Big\} .
\end{align}
A direct calculation of the variations gives the discrete ideal MHD equations,
\begingroup
\allowdisplaybreaks
\begin{subequations}\label{eq:variational_integrator_ideal_mhd}
\begin{align}
\nonumber
0 &= (\Delta_{t} V^{x})_{i, \, j+1/2, \, n+1/2}
+ \psi^{x}_{i, \, j+1/2, \, n+1/2} (V,V) \\
& \hspace{10em}
- \psi^{x}_{i, \, j+1/2, \, n+1/2} (B,B)
+ (\Delta_{x} P)_{i, \, j+1/2, \, n+1/2} ,
\\
\nonumber
0 &= (\Delta_{t} V^{y})_{i+1/2, \, j, \, n+1/2}
+ \psi^{y}_{i+1/2, \, j, \, n+1/2} (V,V) \\
& \hspace{10em}
- \psi^{y}_{i+1/2, \, j, \, n+1/2} (B,B)
+ (\Delta_{y} P)_{i+1/2, \, j, \, n+1/2} ,
\\
0 &= (\Delta_{t} B^{x})_{i, \, j+1/2, \, n+1/2}
+ \phi^{x}_{i, \, j+1/2, \, n+1/2} (V,B) ,
\\
0 &= (\Delta_{t} B^{y})_{i+1/2, \, j, \, n+1/2}
+ \phi^{y}_{i+1/2, \, j, \, n+1/2} (V,B) ,
\\
0 &= (\Delta_{x} V^{x})_{i+1/2, \, j+1/2, \, n+1/2}
   + (\Delta_{y} V^{y})_{i+1/2, \, j+1/2, \, n+1/2} ,
\end{align}
\end{subequations}
\endgroup
with the discrete operators defined by
\begin{subequations}
\begin{align}
\psi^{x}_{i, \, j+1/2, \, n+1/2} (V,V)
\nonumber
&= \dfrac{1}{2} \bar{V}^{y}_{i, \, j,   \, n+1/2} \Big[
     (\Delta_{y} V^{x})_{i, \, j,   \, n+1/2}
   - (\Delta_{x} V^{y})_{i, \, j,   \, n+1/2}
\Big] \\
&+ \dfrac{1}{2} \bar{V}^{y}_{i, \, j+1, \, n+1/2} \Big[
     (\Delta_{y} V^{x})_{i, \, j+1, \, n+1/2}
   - (\Delta_{x} V^{y})_{i, \, j+1, \, n+1/2}
\Big] , \\
\psi^{y}_{i+1/2, \, j, \, n+1/2} (V,V)
\nonumber
&= \dfrac{1}{2} \bar{V}^{x}_{i,   \, j, \, n+1/2} \Big[
     (\Delta_{x} V^{y})_{i,   \, j, \, n+1/2}
   - (\Delta_{y} V^{x})_{i,   \, j, \, n+1/2}
\Big] \\
&+ \dfrac{1}{2} \bar{V}^{x}_{i+1, \, j, \, n+1/2} \Big[
     (\Delta_{x} V^{y})_{i+1, \, j, \, n+1/2}
   - (\Delta_{y} V^{x})_{i+1, \, j, \, n+1/2}
\Big] ,
\end{align}
\end{subequations}
and
\begin{subequations}
\begin{align}
\phi^{x}_{i, \, j+1/2, \, n+1/2} (V,B)
\nonumber
&= \dfrac{1}{2} \Big[
       \bar{V}^{x}_{i, \, j+1, \, n+1/2} \, \bar{B}^{y}_{i, \, j+1, \, n+1/2}
     - \bar{V}^{y}_{i, \, j+1, \, n+1/2} \, \bar{B}^{x}_{i, \, j+1, \, n+1/2}
\Big] \\
&- \dfrac{1}{2} \Big[
       \bar{V}^{x}_{i, \, j, \, n+1/2} \, \bar{B}^{y}_{i, \, j, \, n+1/2}
     - \bar{V}^{y}_{i, \, j, \, n+1/2} \, \bar{B}^{x}_{i, \, j, \, n+1/2}
\Big] ,
\\
\phi^{y}_{i+1/2, \, j, \, n+1/2} (V,B)
\nonumber
&= \dfrac{1}{2} \Big[
       \bar{V}^{y}_{i+1, \, j, \, n+1/2} \, \bar{B}^{x}_{i+1, \, j, \, n+1/2}
     - \bar{V}^{x}_{i+1, \, j, \, n+1/2} \, \bar{B}^{y}_{i+1, \, j, \, n+1/2}
\Big] \\
&- \dfrac{1}{2} \Big[
       \bar{V}^{y}_{i,   j, \, n+1/2} \, \bar{B}^{x}_{i,   j, \, n+1/2}
     - \bar{V}^{x}_{i,   j, \, n+1/2} \, \bar{B}^{y}_{i,   j, \, n+1/2} 
\Big] .
\end{align}
\end{subequations}
Figure~\ref{fig:mhd_staggered_grid_scheme} shows the cells covered by the stencils of each equation.
The discretisation of the operators $\psi$ and $\phi$ is the very same as the one obtained by \citet{Gawlik:2011} and \citet{LiuWang:2001}.
\citeauthor{Gawlik:2011} follow a different but related path in their derivation, based on discrete Euler-Poincar\'{e} reduction. In this approach, only the velocity field is treated variationally while the magnetic field is a quantity passively advected with the velocity field. Instead in our method the velocity field and the magnetic field are treated on equal footing and fully variationally, leading to a different temporal discretisation. With the Euler-Poincar\'{e} integrator, cross helicity is preserved exactly but energy only approximately (cf. Figure~6.1 in Reference~\cite{Gawlik:2011}). We will see that in the proposed scheme, energy, magnetic helicity and cross helicity are preserved exactly (up to machine accuracy).
The scheme of \citeauthor{LiuWang:2001} uses an explicit Runge-Kutta method for time integration, so that conservation laws are broken, possibly posing problems in long time simulations. From the discrete action principle, we obtain the implicit midpoint method for time integration, which is a symmetric integrator that exhibits favourable long-time stability~\cite{HairerLubichWanner:2006}.

Note the absence of any spatial averaging of the time derivatives and the pressure gradient. This is on purpose, as we wanted to prevent the emergence of grid-scale oscillations in the fields by introducing the staggered grid. Spatial averages of the time derivatives or the pressure gradient might result in symmetric finite difference operators in the discrete equations of motion, which support spurious velocity or pressure oscillations between neighbouring grid points.

\section{Numerical Examples}\label{sec:examples}

In this section, we consider four standard test cases of ideal magnetohydrodynamics taken from the previous literature~\cite{CordobaMarliani:2000, GardinerStone:2005, Gawlik:2011}:
Alfv\'{e}n waves (Section~\ref{sec:ex_alfven_waves}), the development of current sheets in an Orszag-Tang vortex (Section~\ref{sec:ex_orszag_tang_vortex}), the passive advection of a magnetic loop (Section~\ref{sec:ex_loop_advection}), and the perturbation of a current sheet (Section~\ref{sec:ex_current_sheet}).

The variational integrator~(\ref{eq:variational_integrator_ideal_mhd}) is implemented using Python~\cite{ScopatzHuff:2015, Langtangen:2014}, Cython~\cite{Behnel:2010}, PETSc~\cite{petsc-web-page, petsc-user-ref} and petsc4py~\cite{Dalcin:2011}. Visualisation was done using NumPy~\cite{vanDerWalt:2011}, SciPy~\cite{scipy-web-page} and matplotlib~\cite{Hunter:2007}.
The nonlinear system is solved with Newton's method, where in each iteration the linear system is solved either via LU decomposition with MUMPS~\cite{Amestoy:2000, Amestoy:2001} or via GMRES with ASM preconditioning~\cite{Gander:2008ua}. The tolerance of the nonlinear solver is set to $10^{-10}$ or smaller, which is usually reached after $3-5$ iterations.

\subsection{Diagnostics}

In the following we give discrete expressions of the conserved quantities, energy~(\ref{eq:ideal_mhd_energy}), cross helicity~(\ref{eq:ideal_mhd_cross_helicity}) and magnetic helicity~(\ref{eq:ideal_mhd_magnetic_helicity}), which are monitored in the simulations, as well as the discrete equations for the reconstruction of the vector potential and the current density.
We do not derive these quantities from a discrete Noether theorem as described in~\cite{KrausMaj:2015}, as the corresponding generating vector fields have horizontal components, which are not supported by the current state of the discrete theory (see~\cite{Kraus:2016:Evolution} for an extension).

\subsubsection*{Energy}

The total energy of the system is the sum of kinetic energy and magnetic energy, which are computed by
\begingroup
\allowdisplaybreaks
\begin{align}
E_{\text{kin}}^{n}
\nonumber
&= \dfrac{1}{2} \bracket{ V_{h} , V_{h} }_{\Omega_{h}}
 = \dfrac{h_{x} h_{y}}{4} \sum \limits_{i, \, j} \Big[
       \big( V^{x}_{i, \, j-1/2, \, n} \big)^{2}
     + \big( V^{x}_{i, \, j+1/2, \, n} \big)^{2} \\
& \hspace{18em}
     + \big( V^{y}_{i-1/2, \, j, \, n} \big)^{2}
     + \big( V^{y}_{i+1/2, \, j, \, n} \big)^{2}
\Big] , \\
E_{\text{mag}}^{n}
\nonumber
&= \dfrac{1}{2} \bracket{ B_{h} , B_{h} }_{\Omega_{h}}
 = \dfrac{h_{x} h_{y}}{4} \sum \limits_{i, \, j} \Big[
       \big( B^{x}_{i, \, j-1/2, \, n} \big)^{2}
     + \big( B^{x}_{i, \, j+1/2, \, n} \big)^{2} \\
& \hspace{18em}
     + \big( B^{y}_{i-1/2, \, j, \, n} \big)^{2}
     + \big( B^{y}_{i+1/2, \, j, \, n} \big)^{2}
\Big] .
\end{align}
\endgroup
As there is no dissipation term in the ideal MHD equations, the total energy should always be preserved.
As usual for incompressible flows, internal energy due to pressure is not accounted for~\cite{MarsdenChorin:1993}.

\subsubsection*{Cross Helicity}

The cross helicity is the $L^{2}$-product of the velocity and the magnetic field,
\begin{align}
C_{\mrm{CH}}^{n}
\nonumber
&= \bracket{ V_{h} , B_{h} }_{\Omega_{h}}
 = \dfrac{h_{x} h_{y}}{2} \sum \limits_{i, \, j} \Big[
       V^{x}_{i, \, j-1/2, \, n} \, B^{x}_{i, \, j-1/2, \, n}
     + V^{x}_{i, \, j+1/2, \, n} \, B^{x}_{i, \, j+1/2, \, n} \\
& \hspace{14em}
     + V^{y}_{i-1/2, \, j, \, n} \, B^{y}_{i-1/2, \, j, \, n}
     + V^{y}_{i+1/2, \, j, \, n} \, B^{y}_{i+1/2, \, j, \, n}
\Big] .
\end{align}
In ideal MHD, the parallel components of the velocity and magnetic fields do not interact, so that the integral of their product over the spatial domain stays constant.

\subsubsection*{Magnetic Helicity}

In two dimensions, magnetic helicity reduces to the integral of the magnetic potential,
\begin{align}
C_{\mrm{MH}}^{n}
&= \int \limits_{\Omega_h} \hodge_{h} A_h
 = h_{x} h_{y} \sum \limits_{i, \, j} A_{i, \, j, \, n} ,
\end{align}
where $A_{i,j,n}$ is reconstructed as described below.

\subsubsection*{Magnetic Potential}

In two dimensions, the magnetic field is given in terms of the magnetic potential $A \in \widetilde{\Lambda}^{0} (\Omega)$ by $B = \hodge \ext A$ or in components by
\begin{align}
B^{x} = \partial_{y} A
\hspace{3em}
\text{and}
\hspace{3em}
B^{y} = - \partial_{x} A ,
\end{align}
where $A$ is the $z$-component of the magnetic vector potential, here treated as a twisted zero-form.
The magnetic potential is collocated at the vertices of Figure~\subref*{fig:mhd_staggered_grid_div}.
Therefore these equations are discretised as $B_{h} = \hodge_{h} \ext_{h} A_{h}$, namely
\begin{align}\label{eq:diagnostics_B}
B^{x}_{i, \, j+1/2, \, n} &=   \dfrac{A_{i, \, j+1, \, n} - A_{i, \, j, \, n}}{h_{y}}   &
& \text{and} &
B^{y}_{i+1/2, \, j, \, n} &= - \dfrac{A_{i+1, \, j, \, n} - A_{i, \, j, \, n}}{h_{x}} .
\end{align}
Equations \eqref{eq:diagnostics_B} can be rewritten as recurrence relations for $A_{i,j}$, namely
\begin{align}
A_{i, \, j+1, \, n} &= A_{i, \, j, \, n} + h_{y} \, B^{x}_{i, \, j+1/2, \, n} &
& \text{and} &
A_{i+1, \, j, \, n} &= A_{i, \, j, \, n} - h_{x} \, B^{y}_{i+1/2, \, j, \, n} .
\end{align}
The vector potential can be obtained by fixing the value of $A$ in the point $(i,j)=(1,1)$ and looping over the whole grid, using the first equation to compute columns and the second to jump between rows, or the other way around.
To which value $A_{1,1}$ is fixed is not important as $A$ is determined only up to a constant.
In a two-dimensional domain, the contour lines of the magnetic potential $A$ correspond to field lines of the magnetic field $B$. Hence, $A$ is an important diagnostic.

\subsubsection*{Current Density}

The current density $J \in \Lambda^{2} (\Omega)$ is given by the exterior derivative of the magnetic field, $J = \ext B$. The discrete version of that is $J_{h} = \ext_{h} B_{h} \in \Lambda^{2}_{h}$, or explicitly,
\begin{align}\label{eq:mhd_diagnostics_current_density}
J_{i, \, j, \, n}
= \dfrac{B^{y}_{i+1/2, \, j, \, n} - B^{y}_{i-1/2, \, j, \, n}}{h_{x}}
- \dfrac{B^{x}_{i, \, j+1/2, \, n} - B^{x}_{i, \, j-1/2, \, n}}{h_{y}} .
\end{align}
Like the vector potential, the current is collocated at cell centres.

\subsection{Alfv\'{e}n Waves}\label{sec:ex_alfven_waves}

In the first example, we consider an Alfv\'{e}n wave traveling along $x$, with initial conditions
\begin{align*}
V^{x} &= 0 , &
V^{y} &= V_{0} \, \sin (\pi x) , &
B^{x} &= B_{0} , &
B^{y} &= B_{0} \, \sin (\pi x) , &
P &= 0.1 ,
\end{align*}
with $V_{0} = 1$ and $B_{0} = 1$.
The simulation domain is $\Omega = [0,2] \times [0,2]$ with periodic boundaries and a resolution of $n_{x} \times n_{y} = 32 \times 32$. The time step is $h_{t} = 0.1$ in units of the Alfv\'{e}n time (i.e., the Alfv\'{e}n velocity is one).

\begin{figure}[tb]
\centering
\subfloat[Alfv\'{e}n Wave]{
\label{fig:mhd_alfven_wave_travelling}
\begin{minipage}{.48\textwidth}
\includegraphics[width=\textwidth]{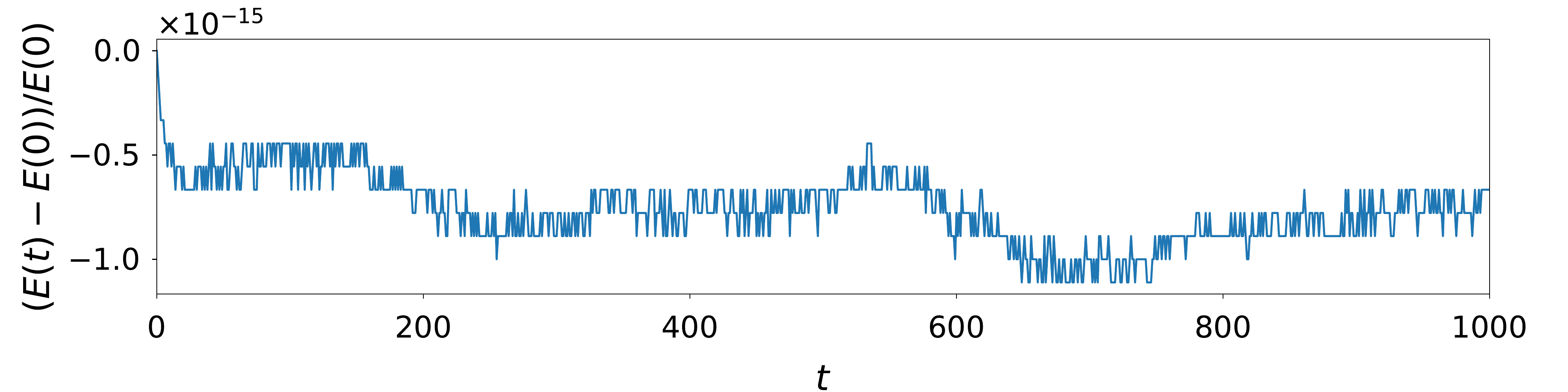}

\includegraphics[width=\textwidth]{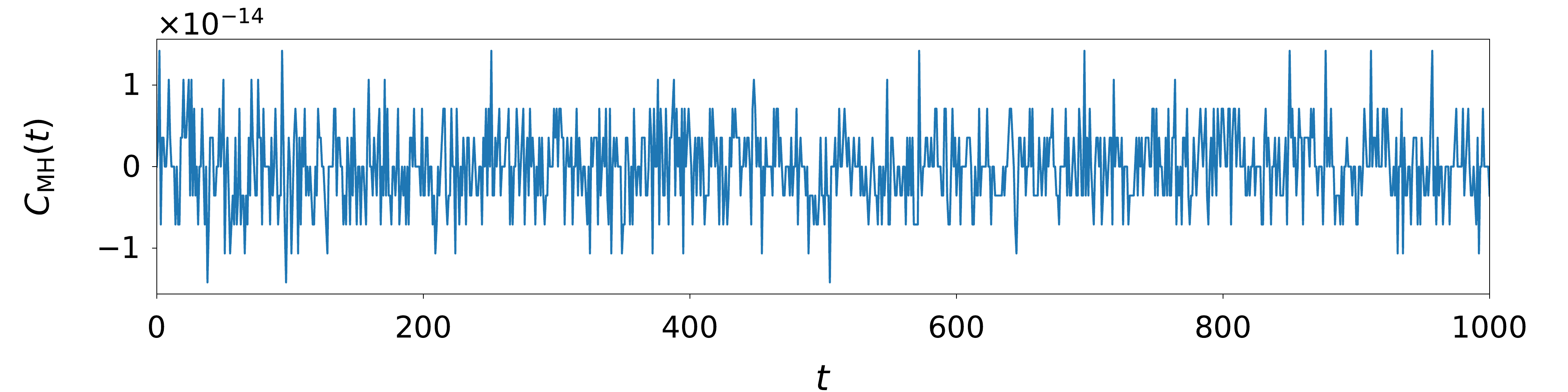}

\includegraphics[width=\textwidth]{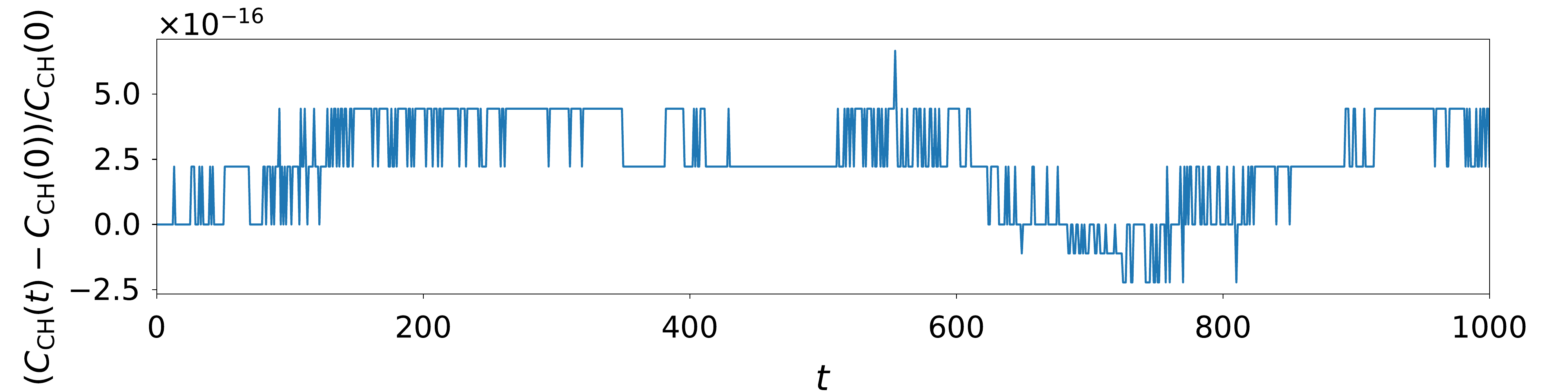}
\end{minipage}
}
\subfloat[Orszag Tang Vortex]{
\label{fig:orszag_tang_64x64_errors}
\begin{minipage}{.48\textwidth}
\includegraphics[width=\textwidth]{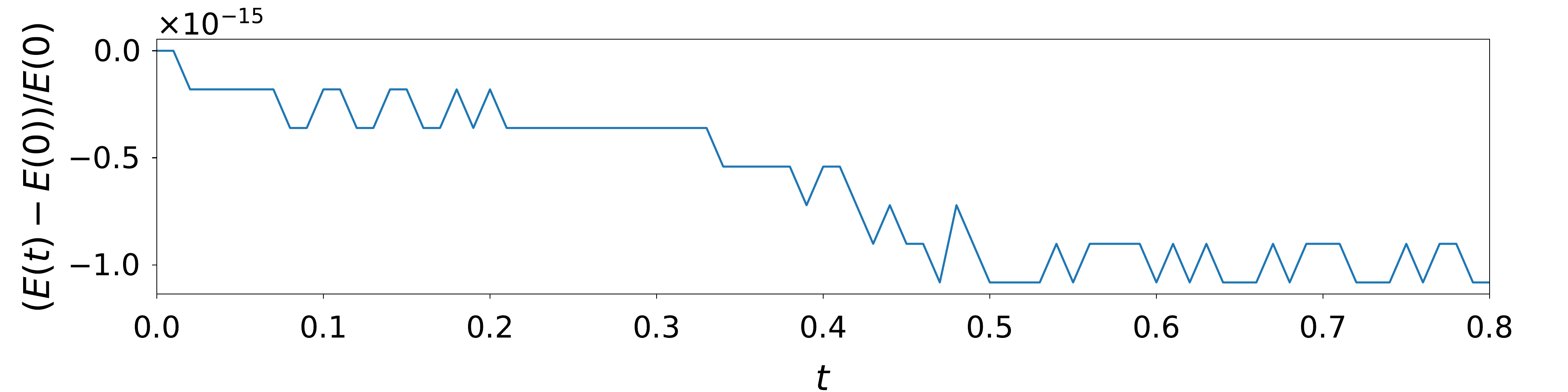}

\includegraphics[width=\textwidth]{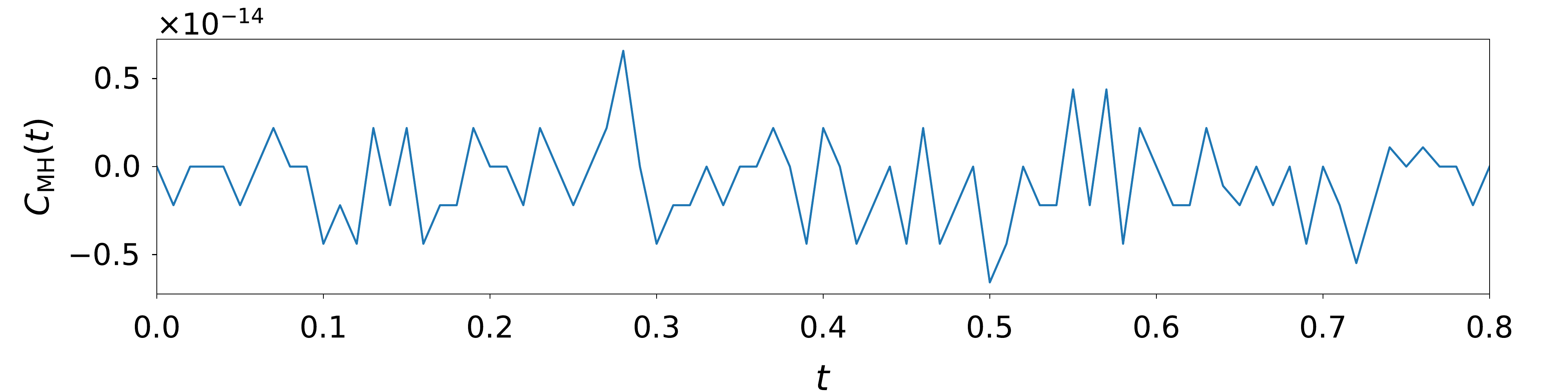}

\includegraphics[width=\textwidth]{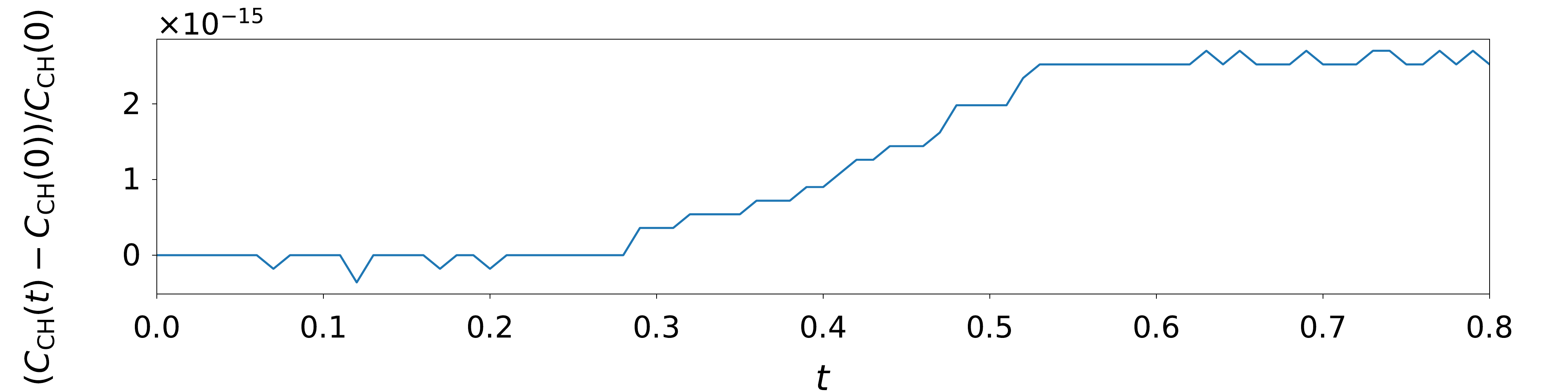}
\end{minipage}
}

\caption{Conservation of energy, magnetic helicity and cross helicity for a travelling Alfv\'{e}n wave (left) and Orszag Tang vortex (right).}
\end{figure}

\begin{figure}[tb]
	\centering
	\subfloat{
		\includegraphics[width=.35\textwidth]{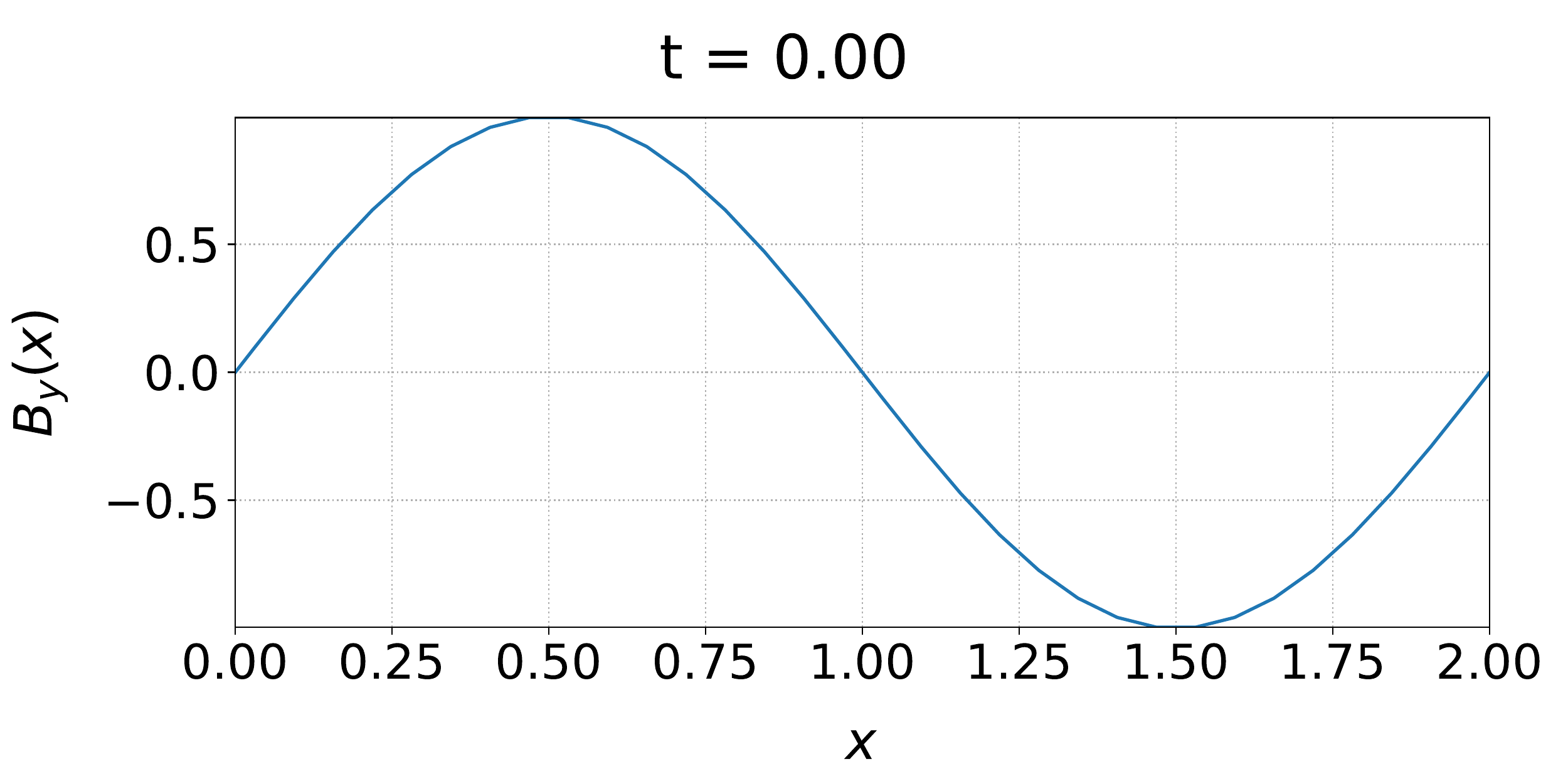}
	}
	\subfloat{
		\includegraphics[width=.35\textwidth]{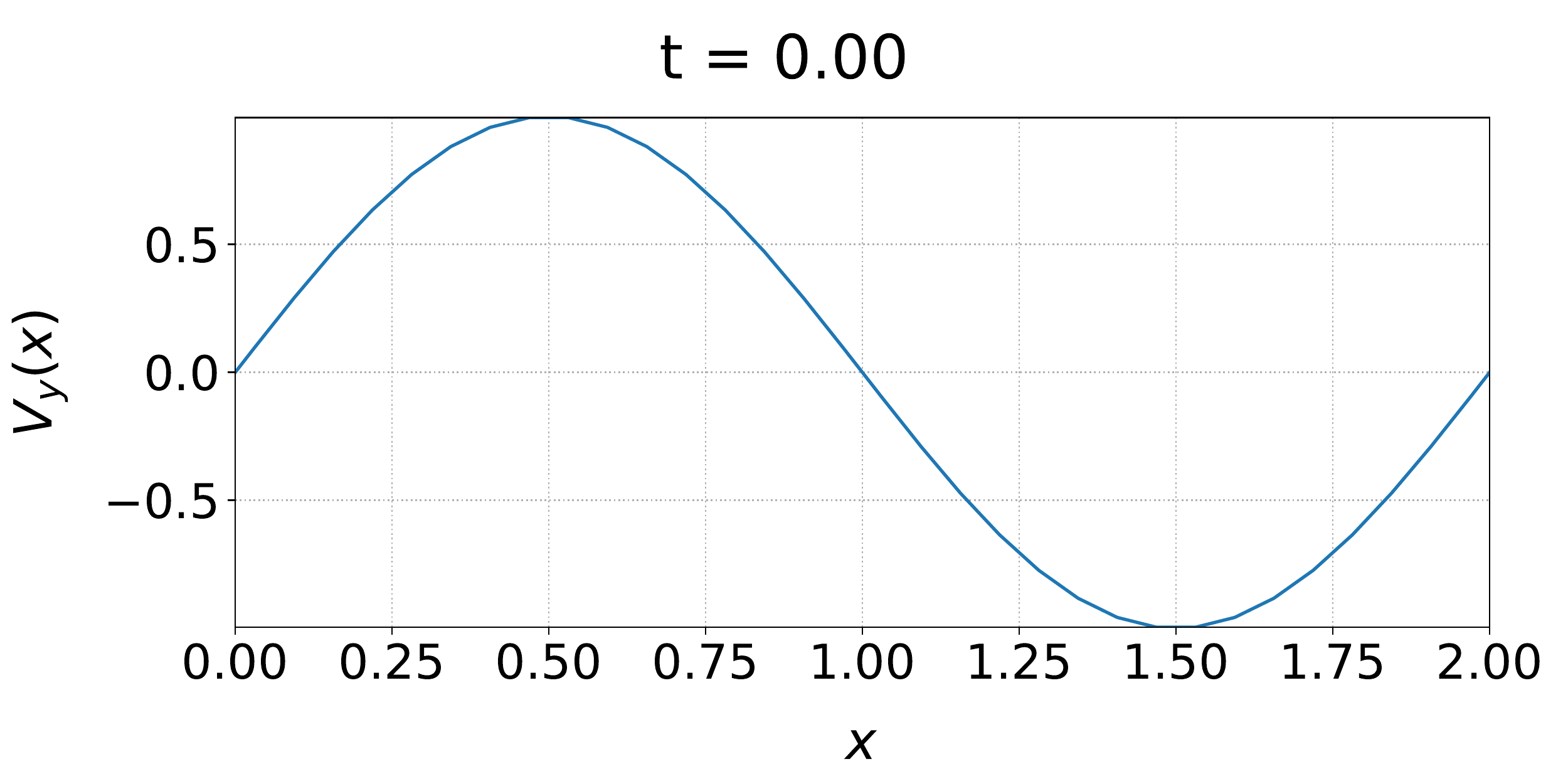}
	}

	\subfloat{
		\includegraphics[width=.35\textwidth]{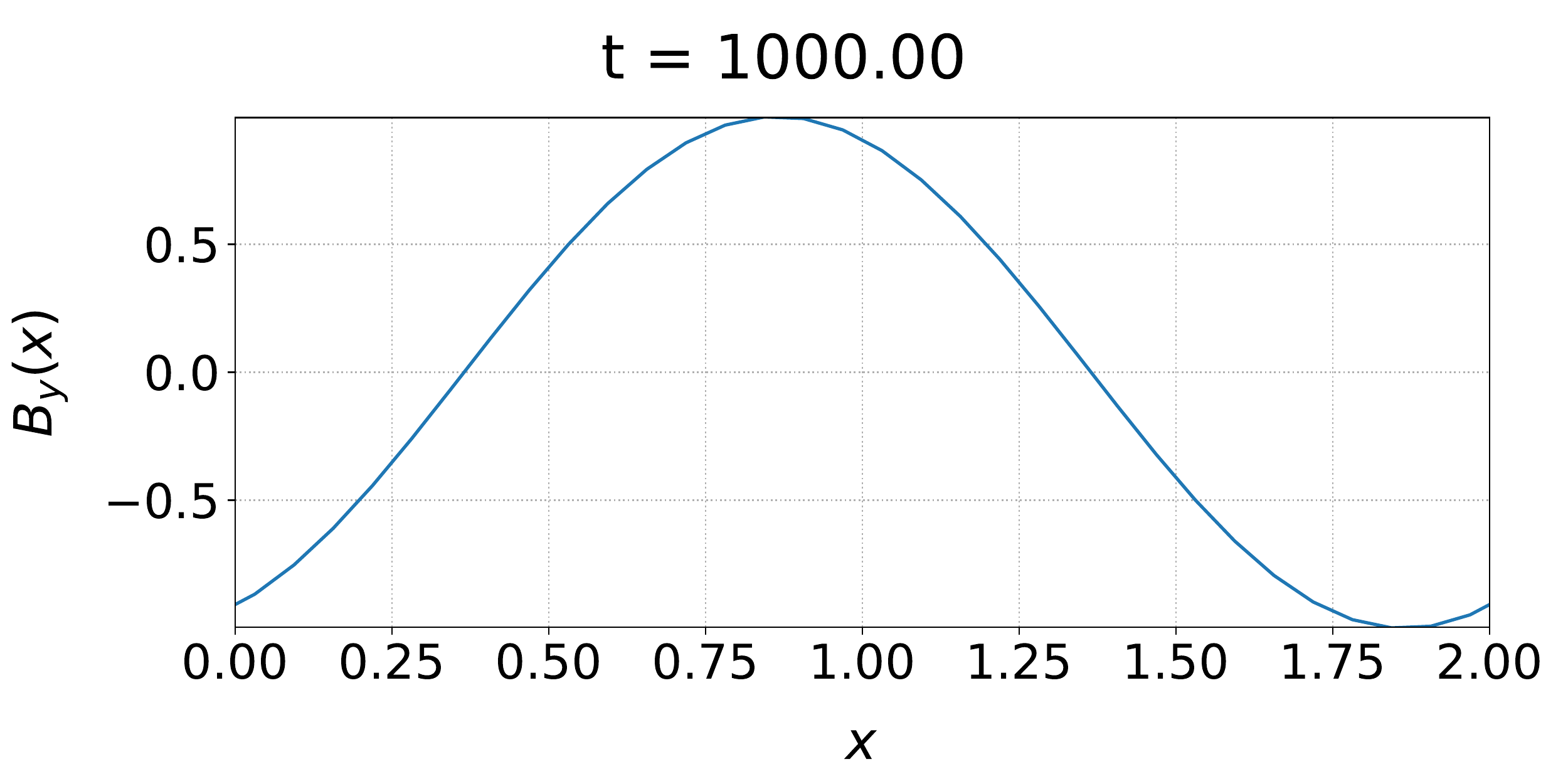}
	}
	\subfloat{
		\includegraphics[width=.35\textwidth]{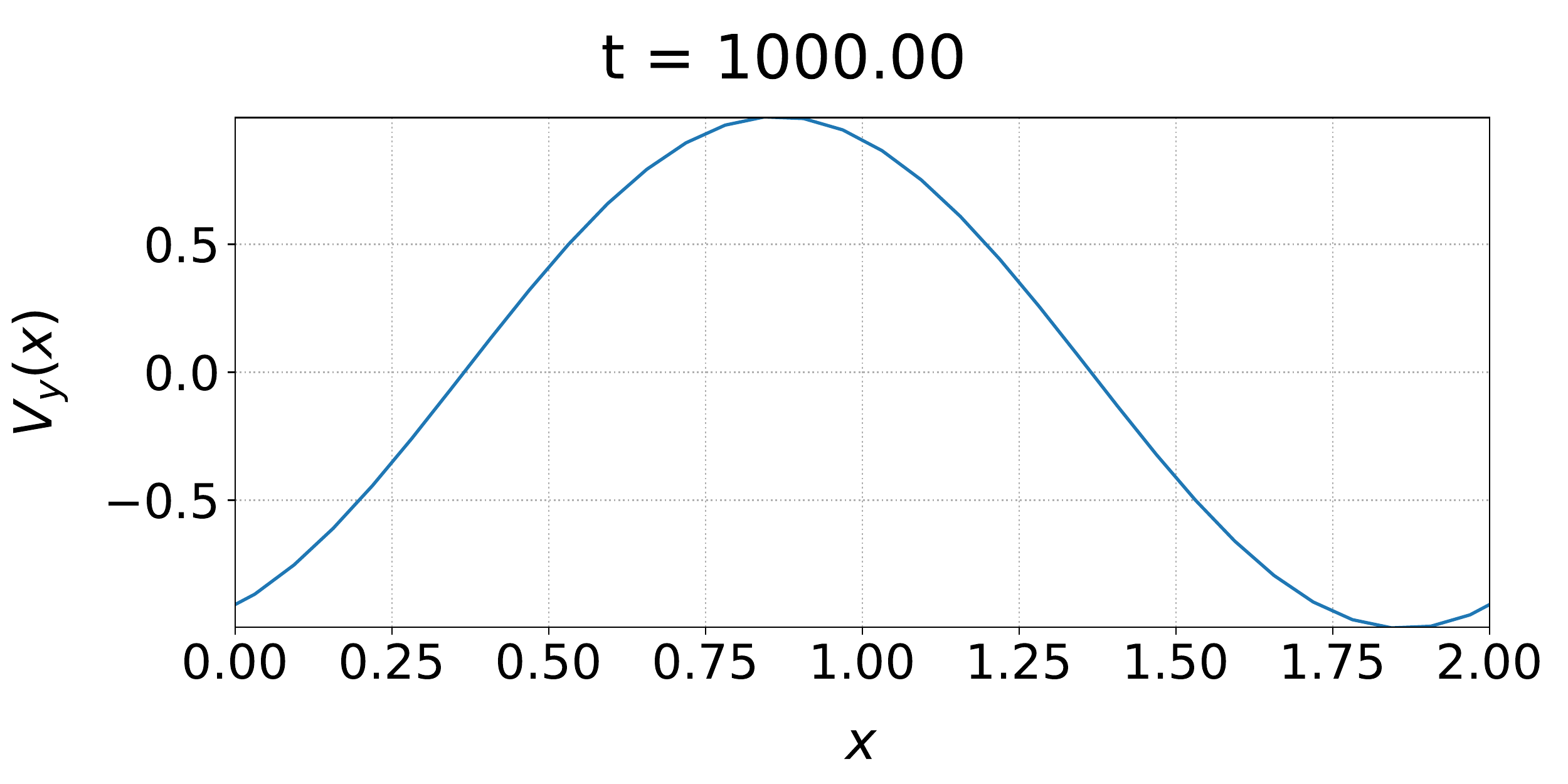}
	}

	\caption{Alfv\'{e}n wave, $x$-profile of the $y$-components of the magnetic field and the velocity at $y = 1$.}
	\label{fig:alfven_wave_travelling_waveform}
\end{figure}

Although this example is rather simple, the results of our variational integrator are already remarkable.
Figure~\ref{fig:mhd_alfven_wave_travelling} shows the time traces of the errors in the total energy, the magnetic helicity and the cross helicity. For most of the simulation, the amplitude of the oscillations is of order $10^{-15}$, i.e., machine precision.
We want to stress that within $500$ passings of the wave there is no change in the energy within the machine accuracy, implying that there is no damping due to numerical effects.
In this respect, it is worth emphasising that this is a fully nonlinear wave, i.e., the amplitudes of the perturbations of the magnetic field as well as the velocity field are of order one.
As can be seen in Figure~\ref{fig:alfven_wave_travelling_waveform}, the shape of the wave is also well preserved, albeit there is some phase error: for $t = 2n$, with integer $n$, the solution should reproduce the initial condition, which is not the case at $t = 1000$.
As can be seen by Fourier analysis, the phase velocity in the simulation is approximately $0.9901$ and thus does not match the theoretical Alfv\'{e}n velocity, which is one.
This is not surprising given the low order and dispersiveness of the proposed method (c.f. the discussion in Reference~\cite{KrausMaj:2015}).

\begin{figure}[p]
\centering

\subfloat{
\includegraphics[width=.42\textwidth]{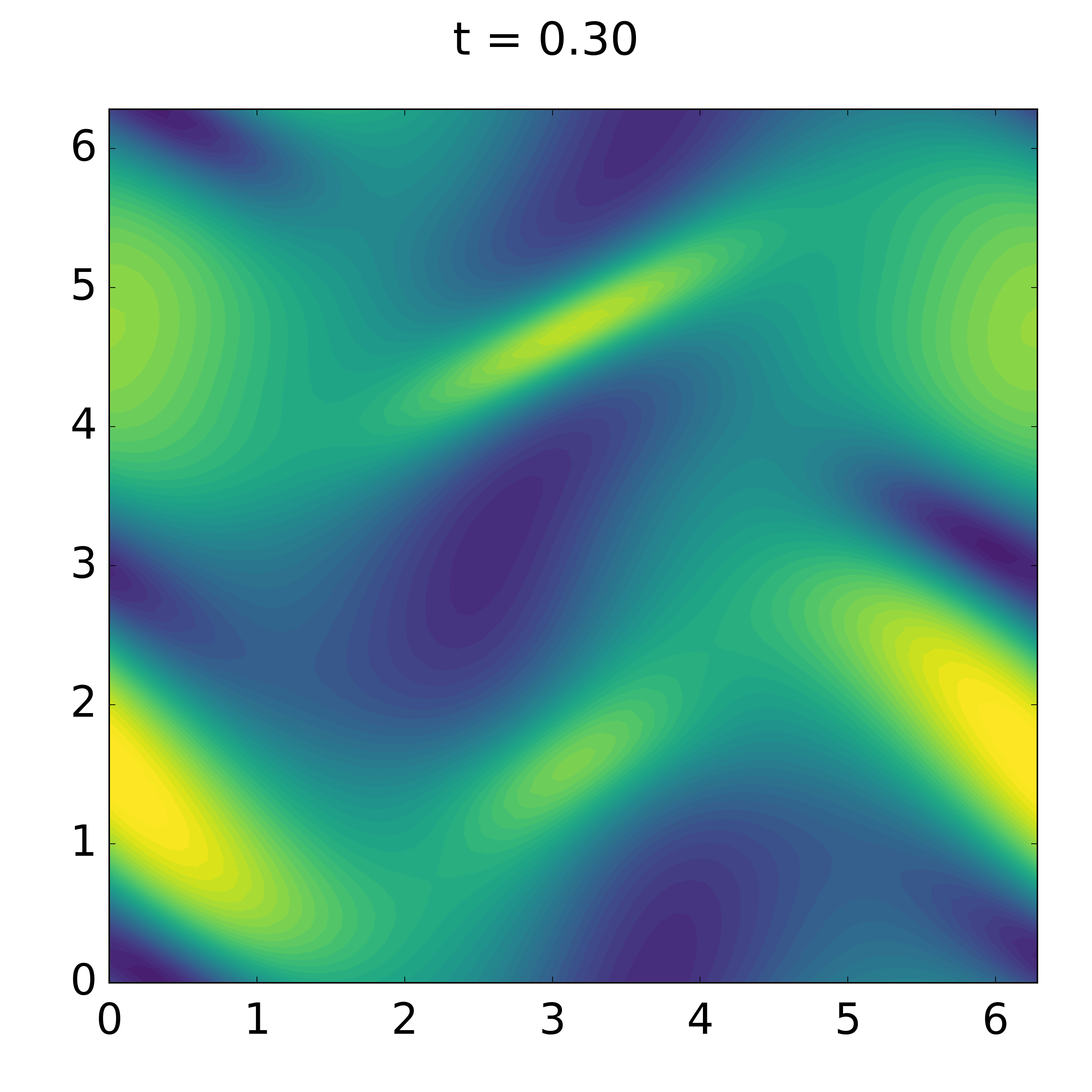}
}
\subfloat{
\includegraphics[width=.42\textwidth]{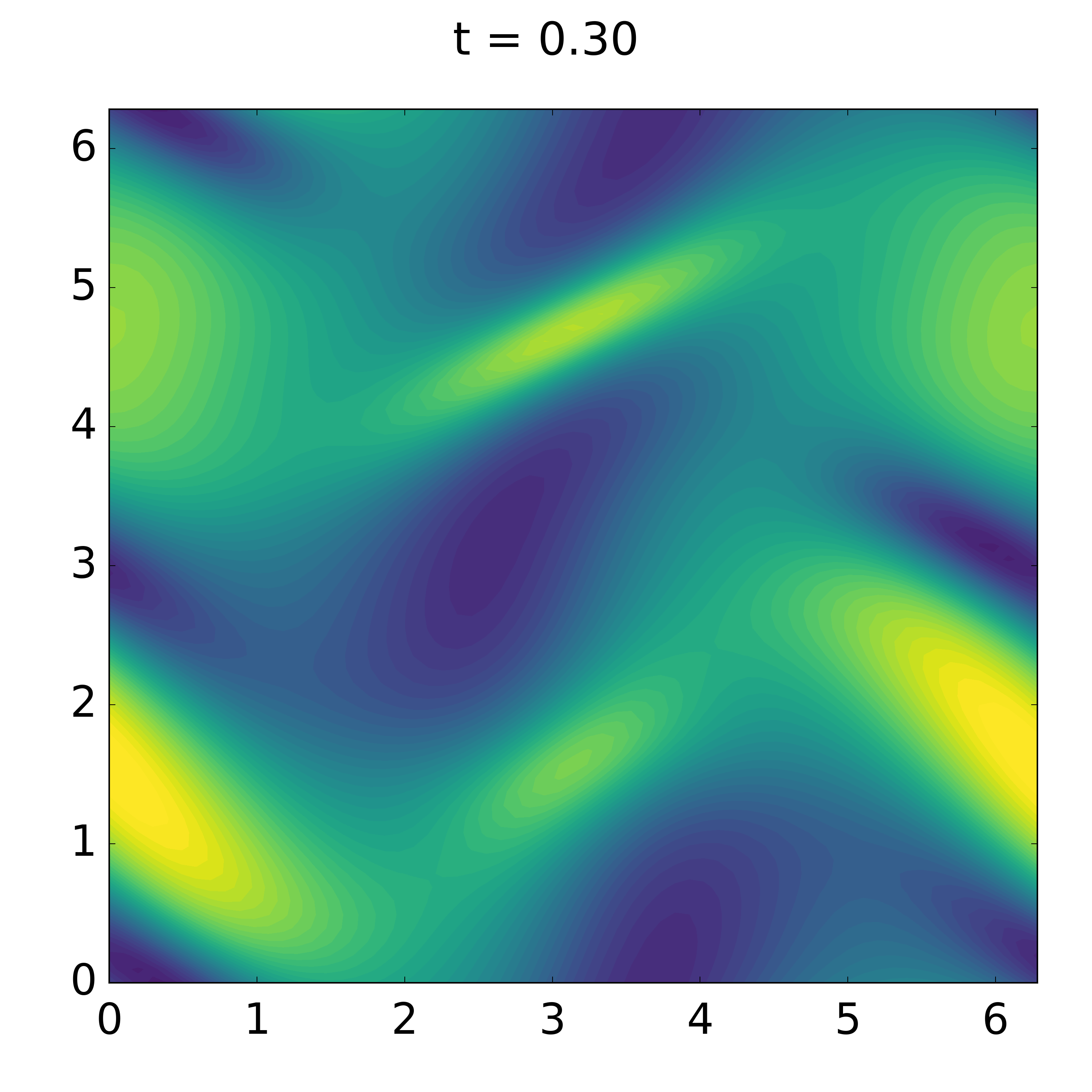}
}

\subfloat{
\includegraphics[width=.42\textwidth]{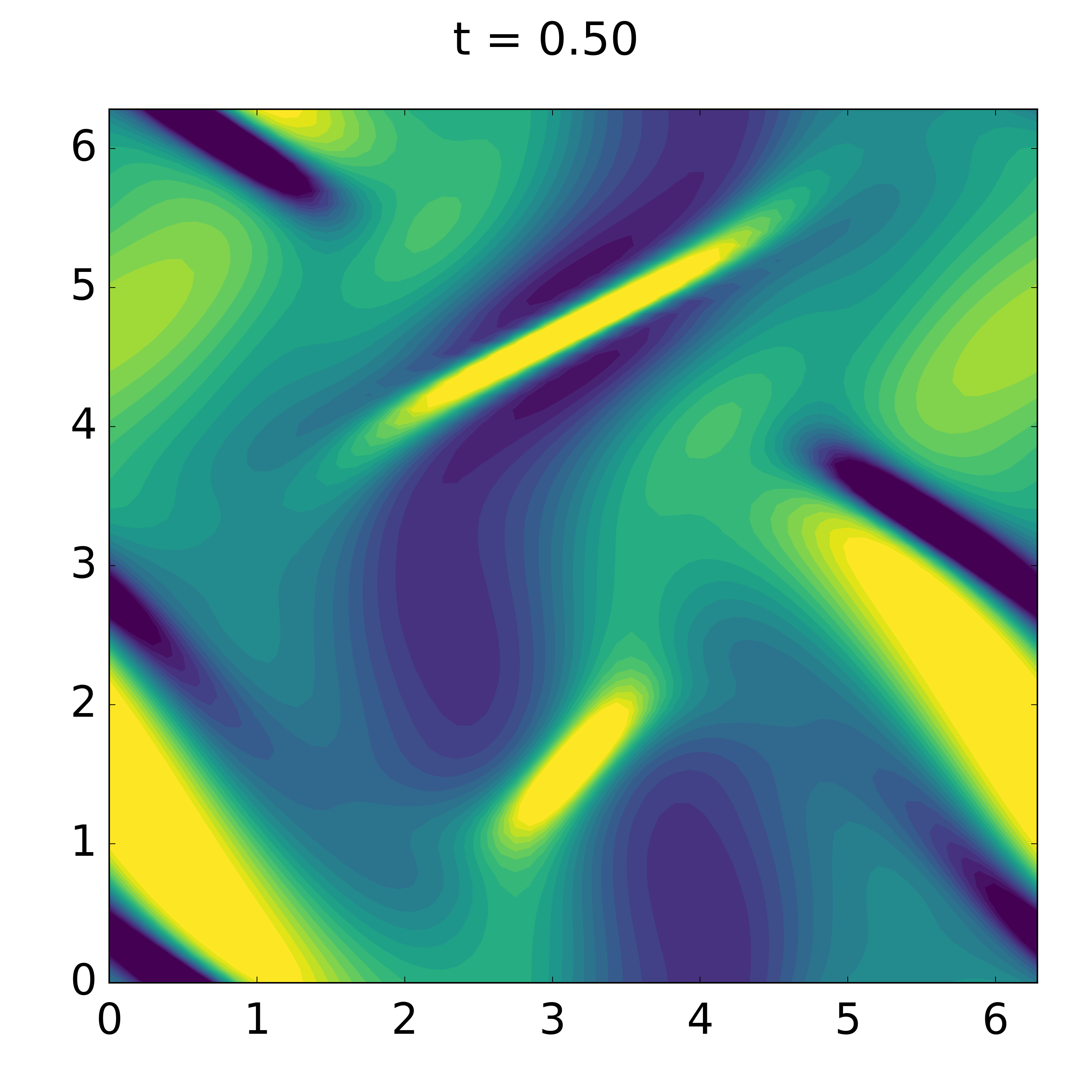}
}
\subfloat{
\includegraphics[width=.42\textwidth]{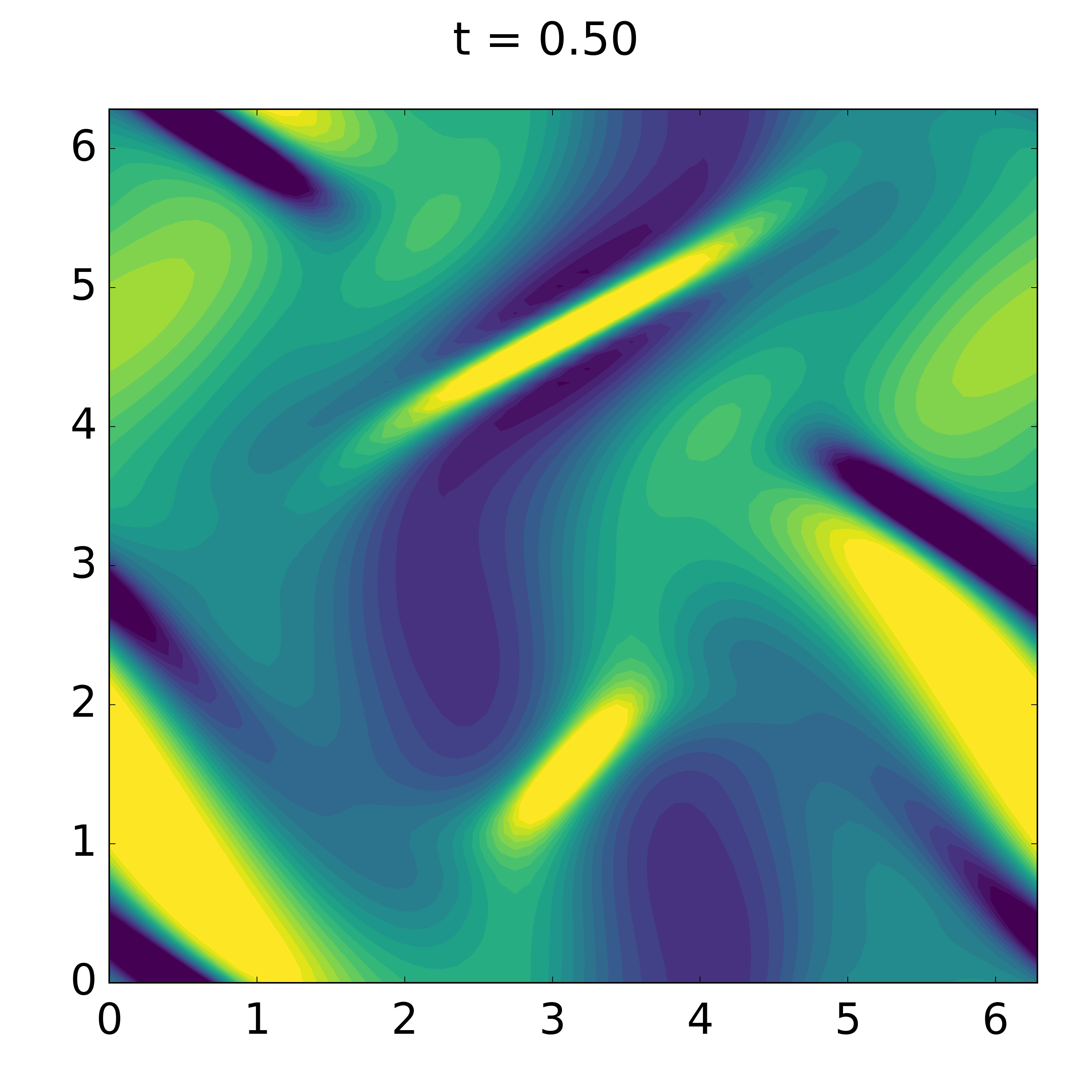}
}

\subfloat{
\includegraphics[width=.42\textwidth]{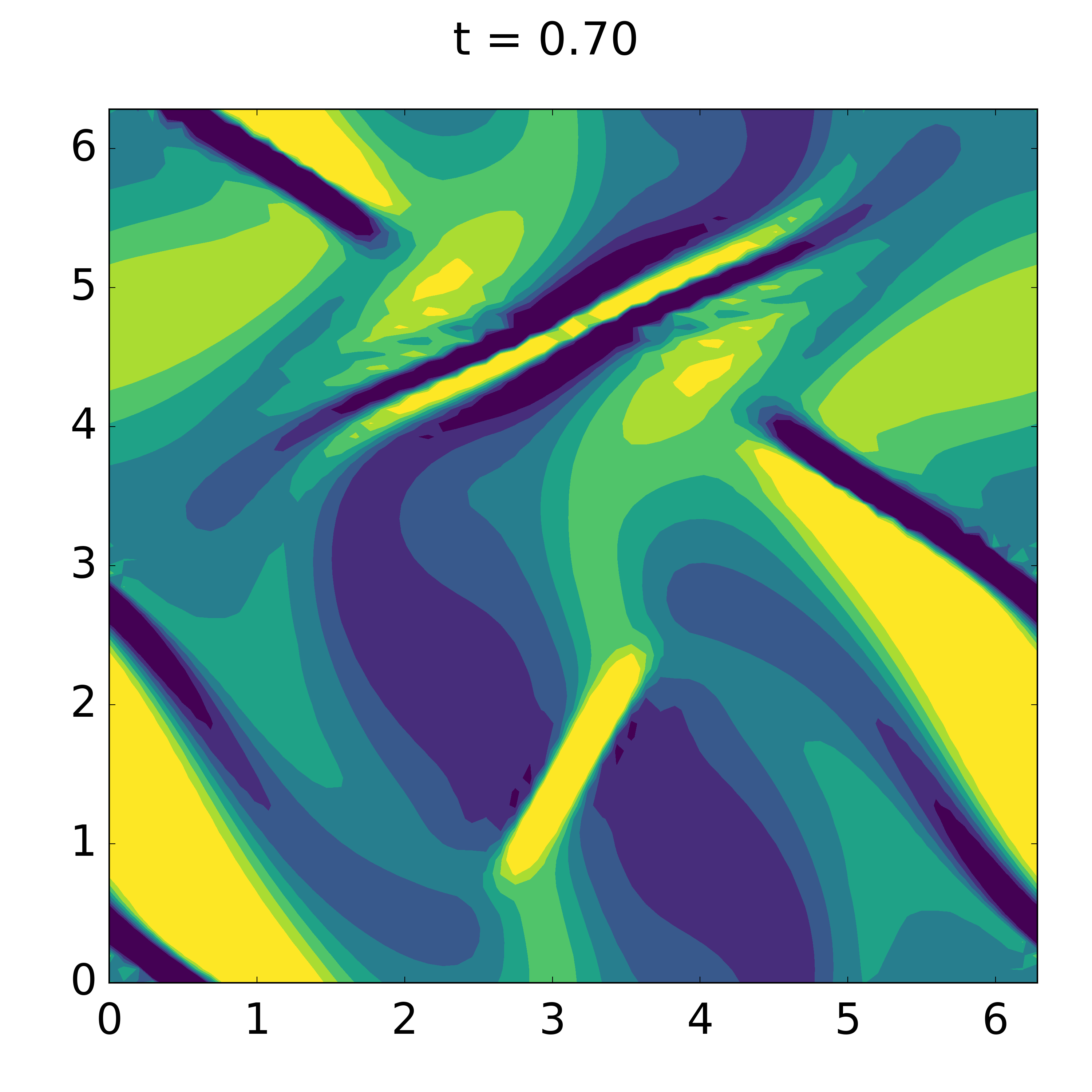}
}
\subfloat{
\includegraphics[width=.42\textwidth]{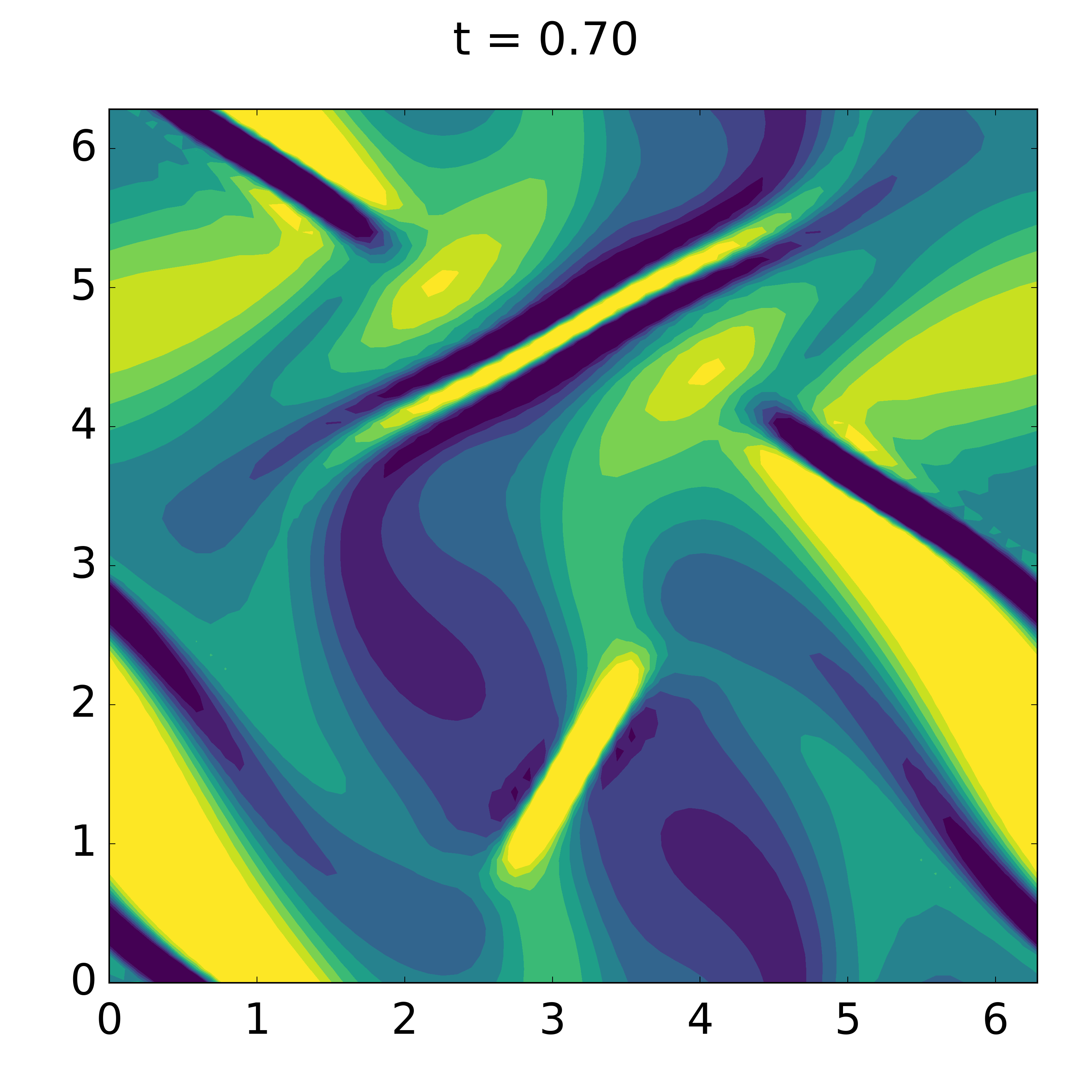}
}

\caption{Orszag Tang Vortex with ideal MHD integrator (left) and reduced MHD integrator (right). Current density $j$. Fixed colour scale.}
\label{fig:orszag_tang_vortex_64x64_current_density}
\end{figure}

\subsection{Orszag-Tang Vortex}\label{sec:ex_orszag_tang_vortex}

Next we consider the evolution of current sheets in an Orszag-Tang vortex, where we use the same initial conditions as \citet{CordobaMarliani:2000}, namely
\begin{align*}
V^{x} &=   \partial_{y} \psi, &
V^{y} &= - \partial_{x} \psi, &
B^{x} &=   \partial_{y} A, &
B^{y} &= - \partial_{x} A , &
P &= 0.1 ,
\end{align*}
with the streaming function $\psi$ and the vector potential $A$ given by
\begin{align*}
\psi &= 2 \sin (y) - 2 \cos (x) , &
A    &= \cos (2y) - 2 \cos (x) .
\end{align*}
The spatial domain is $\Omega = [0,2 \pi] \times [0,2 \pi]$ with periodic boundaries. We consider a spatial resolution of $n_{x} \times n_{y} = 64 \times 64$ grid points and the time step $h_{t} = 0.01$.

The Orszag-Tang vortex constitutes a turbulent setting that leads to the growths of current sheets. These are narrow areas of large current density due to a change of sign in the magnetic field.
In Figure~\ref{fig:orszag_tang_vortex_64x64_current_density}, the current density computed by~(\ref{eq:mhd_diagnostics_current_density}) is plotted.
Current sheets are those regions of the domain where the current density concentrates.
Starting from about $t = 60$, the simulation is under-resolved and subgrid modes start to impair the quality of the solution.
Note that in the original work, \citet{CordobaMarliani:2000} used an adaptive mesh refinement approach with an initial resolution of $1024 \times 1024$ points.
The important observation is, that even with low resolution energy, magnetic helicity and cross helicity are preserved to machine precision (see Figure~\ref{fig:orszag_tang_64x64_errors}). Even though a slight growth in the errors is observed, the amplitudes do not exceed $3 \times 10^{-15}$ throughout the whole simulation.

\subsection{Loop Advection}\label{sec:ex_loop_advection}

We now consider a case with very small magnetic field, such that the momentum and induction equations are almost decoupled and the magnetic field is passively advected by the fluid.
The initial conditions are the ones proposed by~\citet{GardinerStone:2005}, namely
\begin{align*}
V^{x} &= V_{0} \, \cos (\theta) , &
V^{y} &= V_{0} \, \sin (\theta) , &
B^{x} &= \partial_{y} A, &
B^{y} &= - \partial_{x} A, &
P &= 1.0 ,
\end{align*}
with the magnetic potential $A$ given by
\begin{align*}
A &= \begin{cases}
A_{0} \, ( R - \sqrt{x^{2} + y^{2}} ) & \quad \text{for} \; r \leq R , \\
0 & \quad \text{for} \; r > R , \\
\end{cases}
\end{align*}
essentially describing a cone, and the following constants
\begin{align*}
V_{0} &= \sqrt{L_x^2 + L_y^2} = \sqrt{5}, &
\theta &= \tan^{-1} (L_y / L_x) = \tan^{-1} (0.5) , &
A_{0} &= 10^{-3}, &
R &= 0.3 .
\end{align*}
The spatial domain is $\Omega = [-1, +1] \times [-0.5, +0.5]$ with periodic boundaries, so that the lengths of the domain are $L_x = 2$ and $L_y = 1$. We consider a spatial resolution of $n_{x} \times n_{y} = 128 \times 64$ grid points and the time step $h_{t} = 0.01$.

\begin{figure}[bth]
\centering
\subfloat[Loop advection, non-smooth, $128 \times 64$]{\label{fig:loop_advection_128x64_errors}
\begin{minipage}{.48\textwidth}
\includegraphics[width=\textwidth]{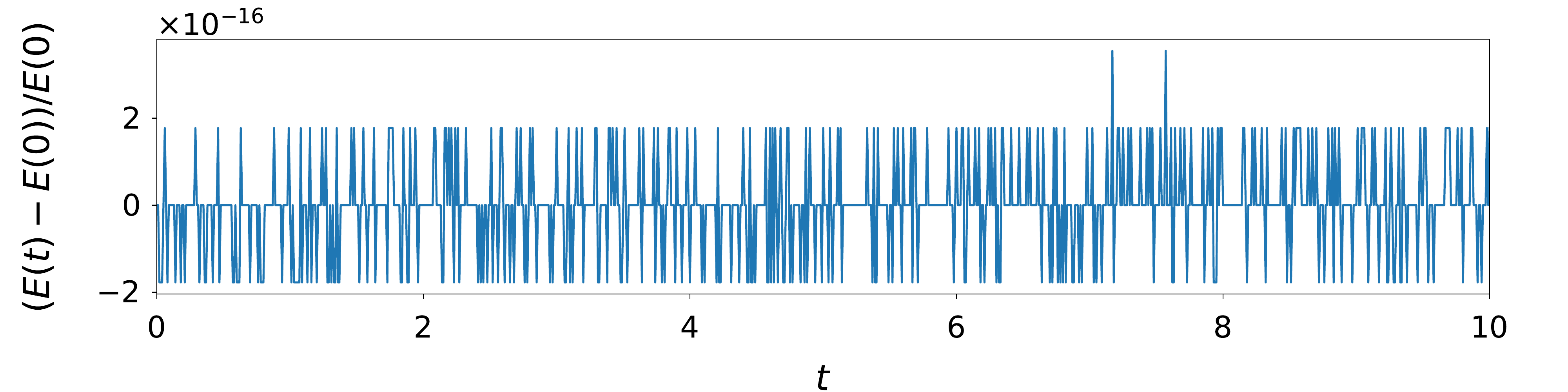}

\includegraphics[width=\textwidth]{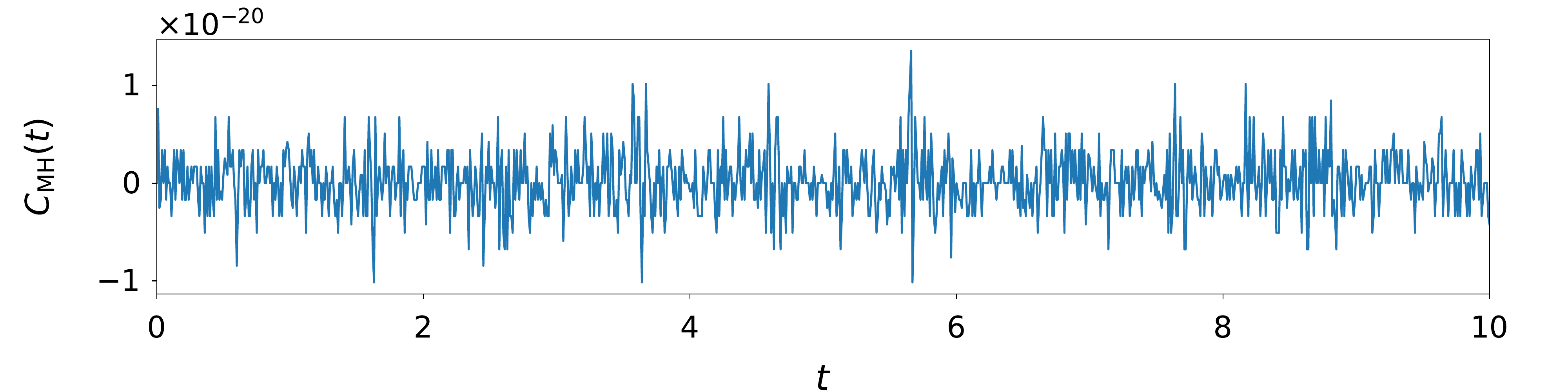}

\includegraphics[width=\textwidth]{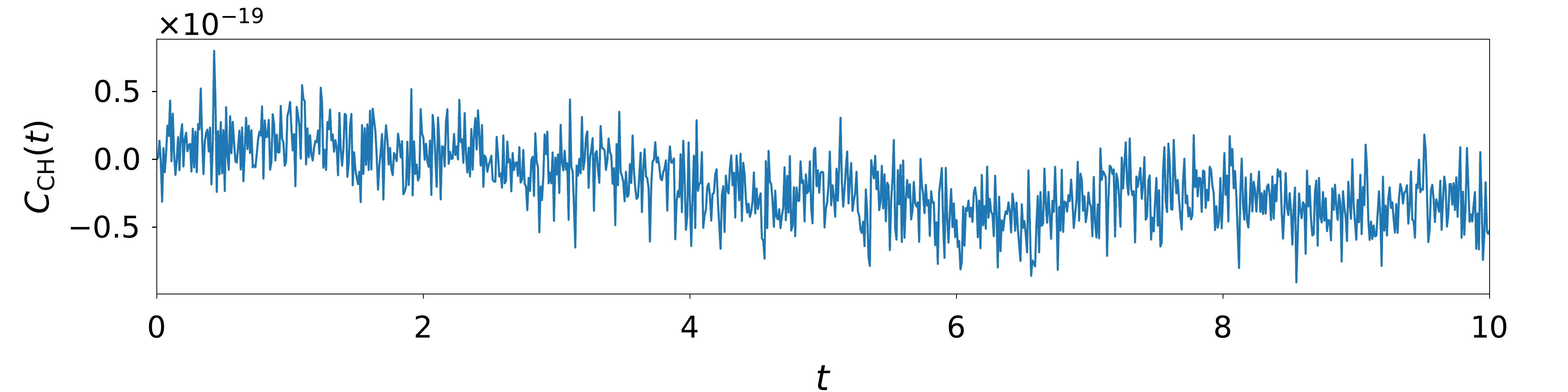}
\end{minipage}
}
\subfloat[Loop advection, smooth, $256 \times 256$]{\label{fig:loop_advection_periodic_errors}
\begin{minipage}{.48\textwidth}
\includegraphics[width=\textwidth]{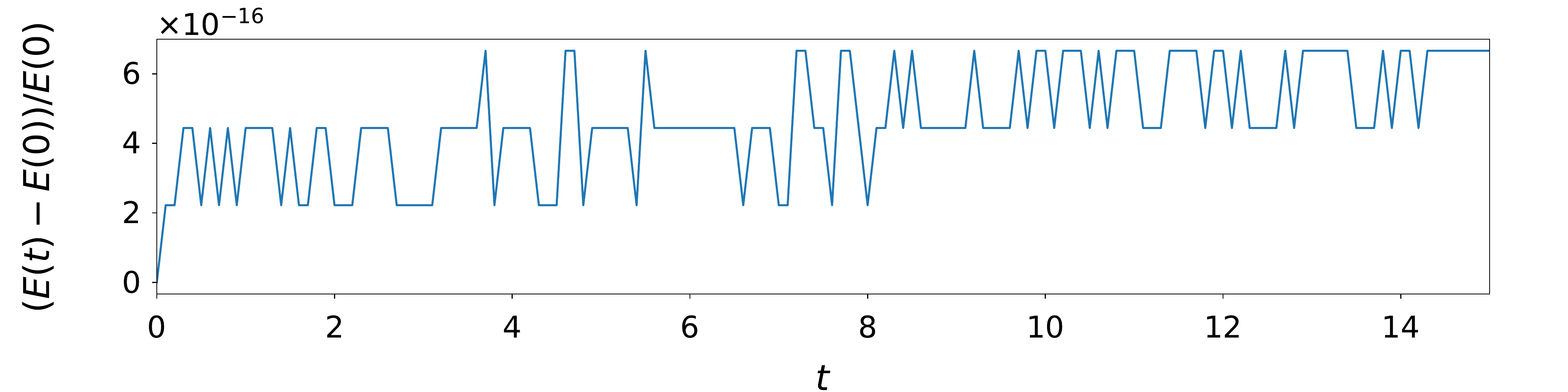}

\includegraphics[width=\textwidth]{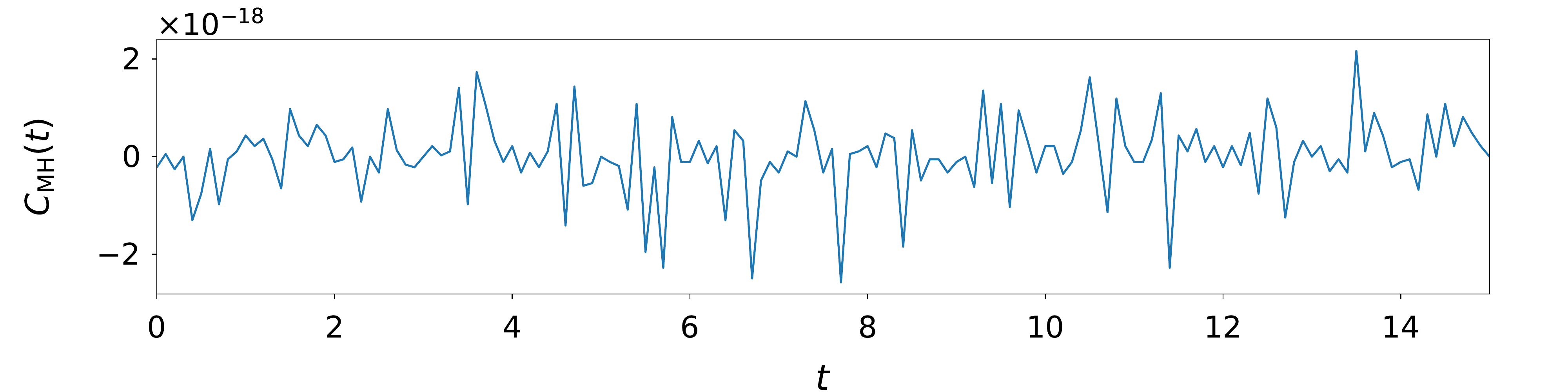}

\includegraphics[width=\textwidth]{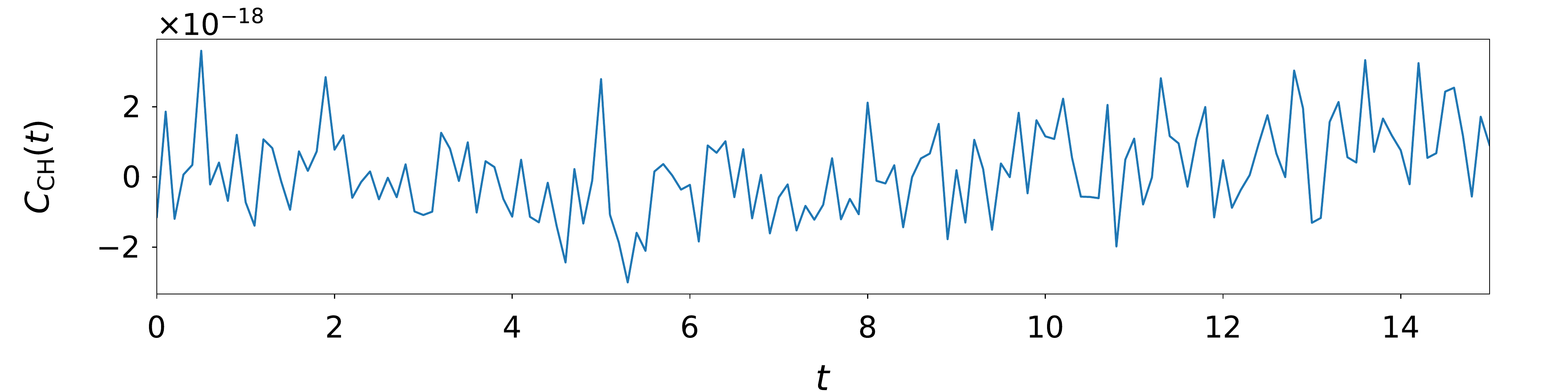}
\end{minipage}
}

\caption{Conservation of energy, magnetic helicity and cross helicity for loop advection, non-smooth (left) and smooth (right).}
\end{figure}

\begin{figure}[p]
\centering
\includegraphics[width=.4\textwidth]{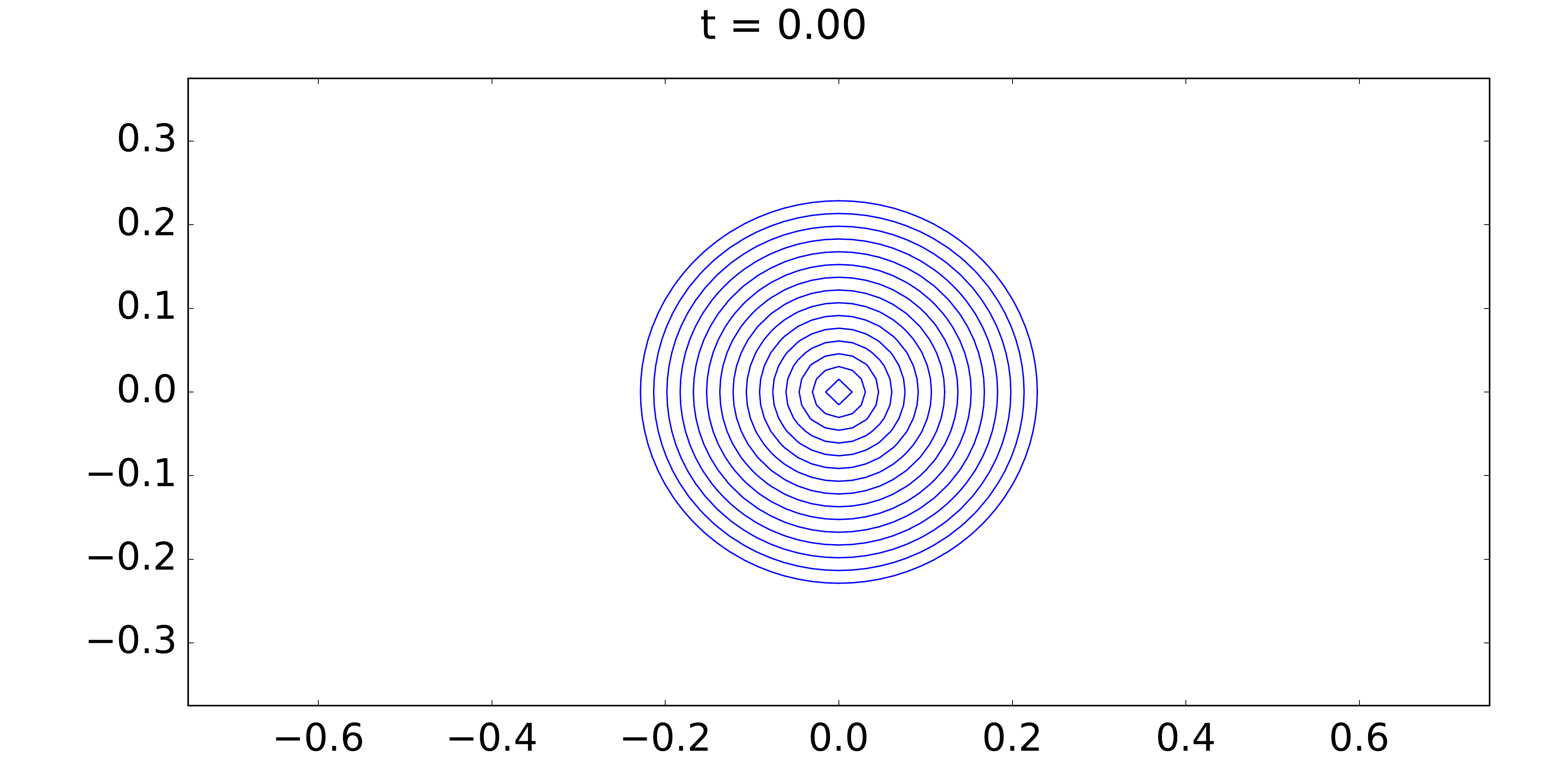}
\includegraphics[width=.4\textwidth]{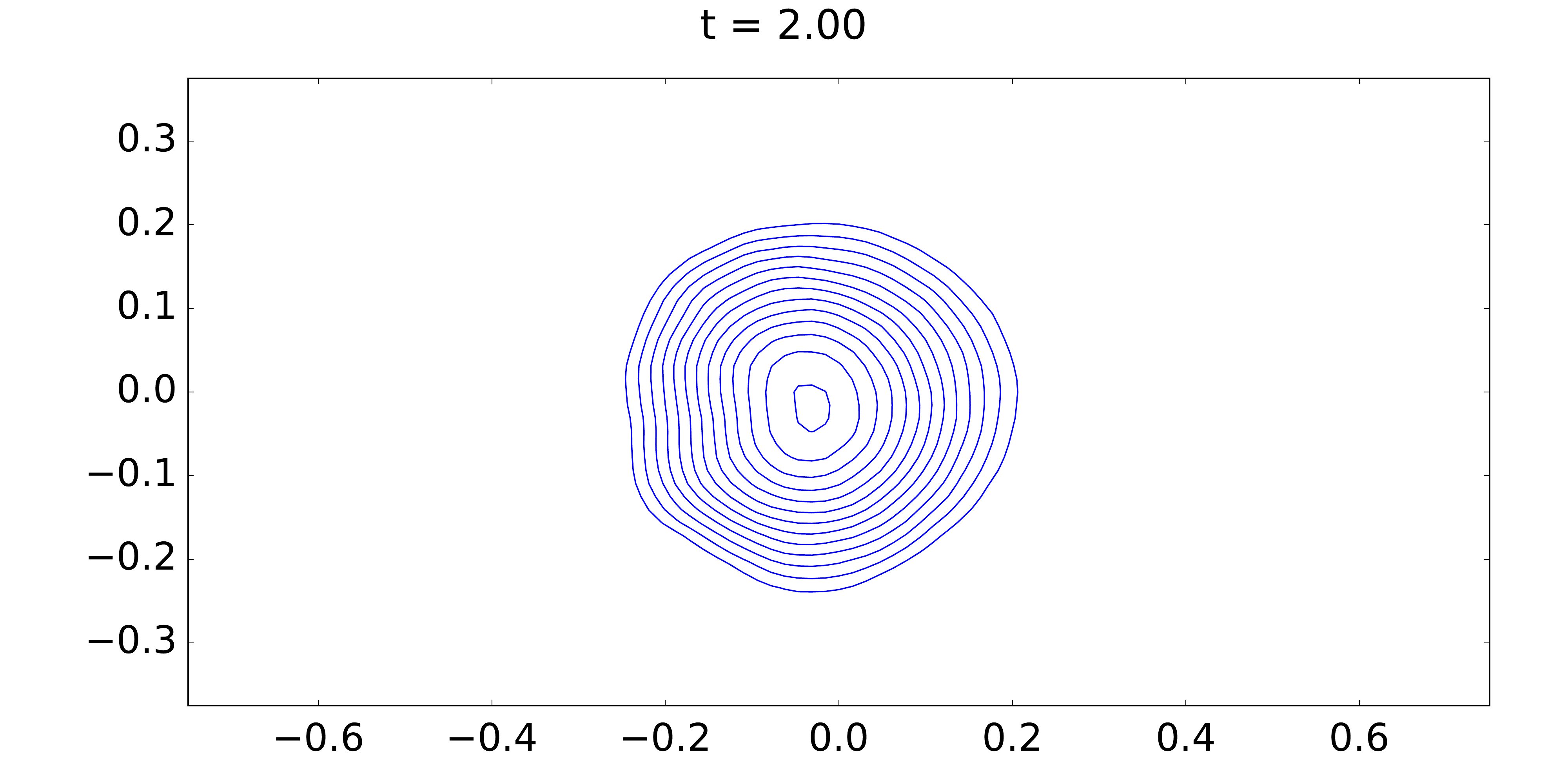}

\includegraphics[width=.4\textwidth]{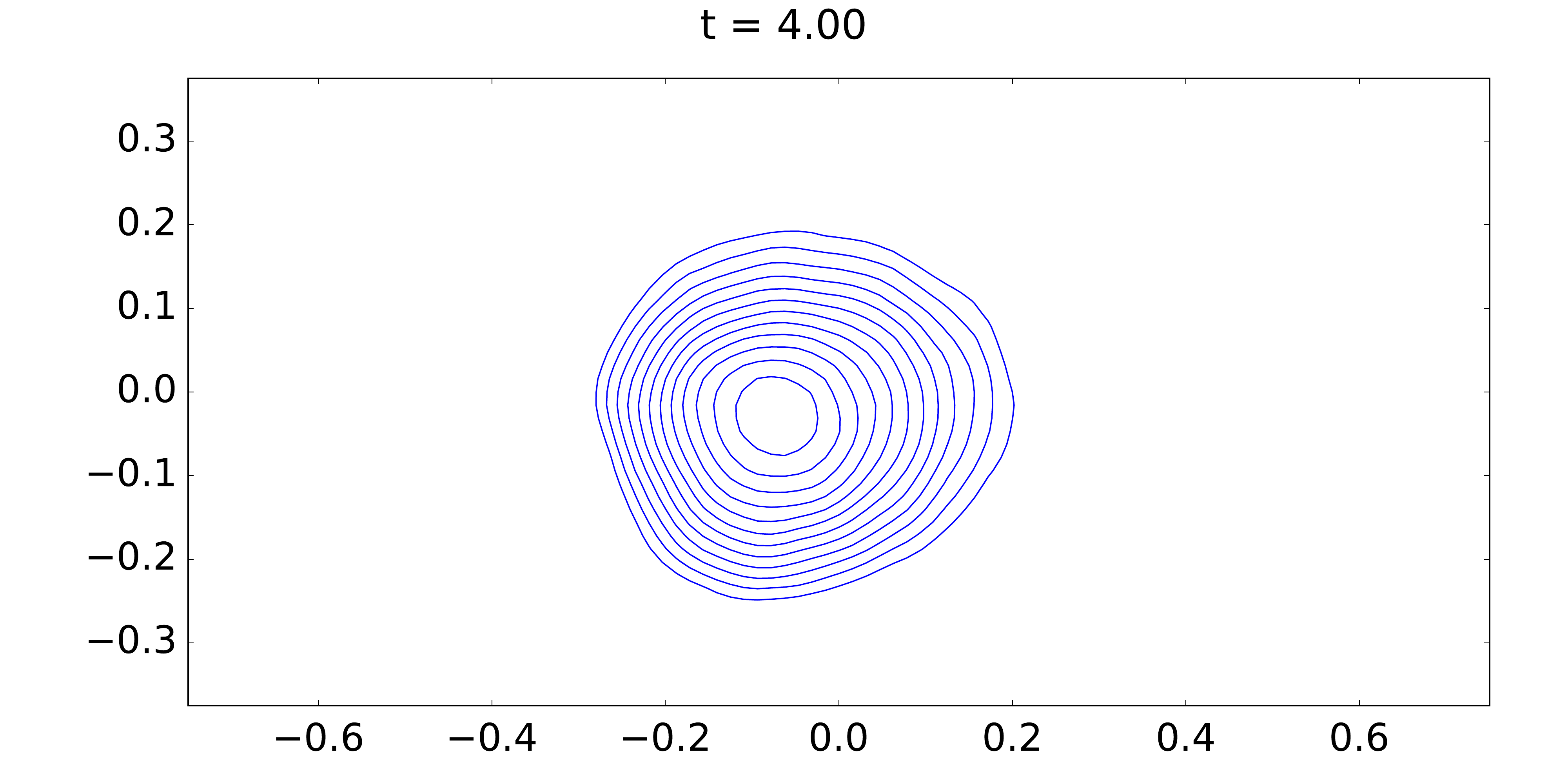}
\includegraphics[width=.4\textwidth]{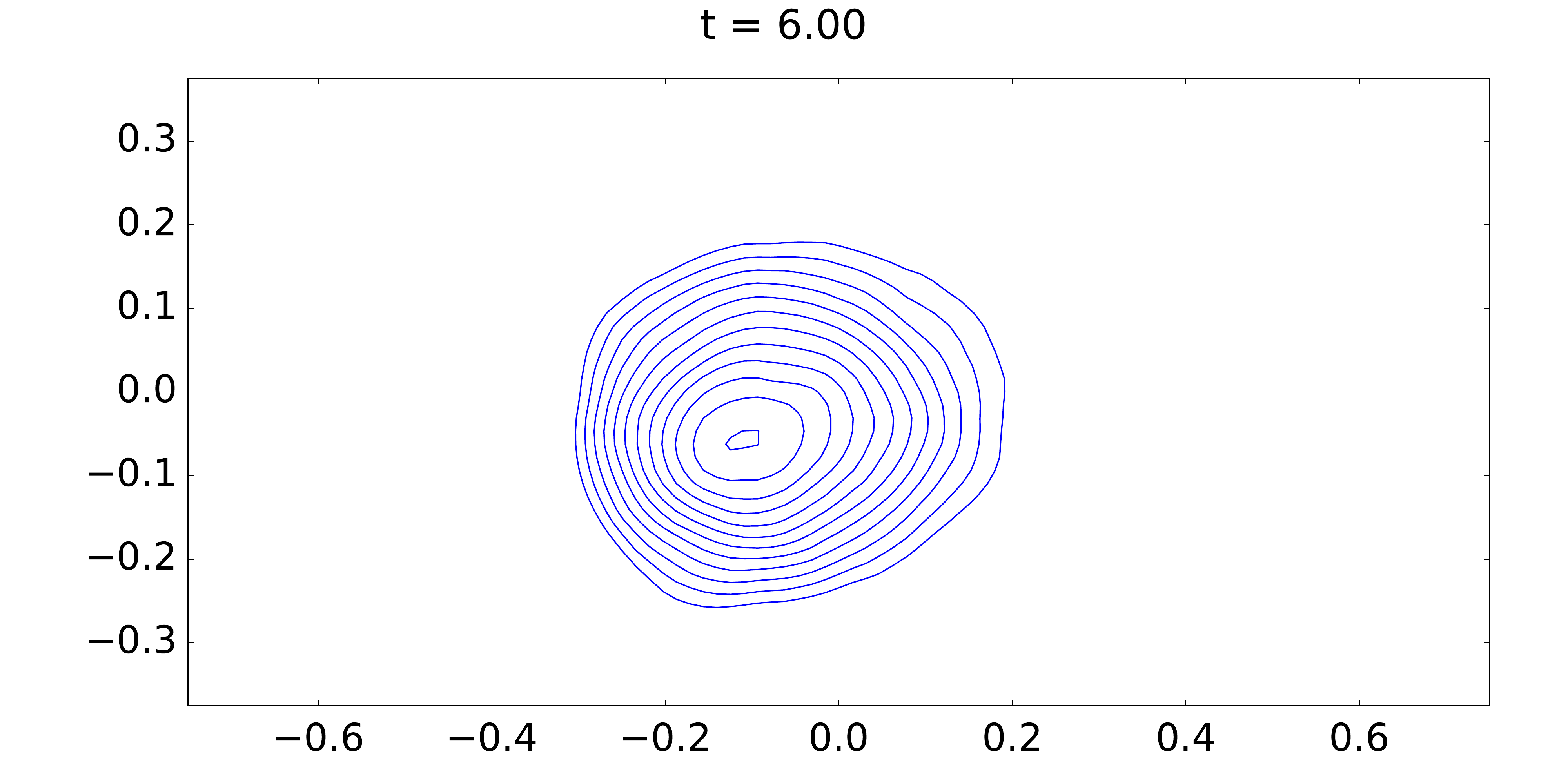}

\includegraphics[width=.4\textwidth]{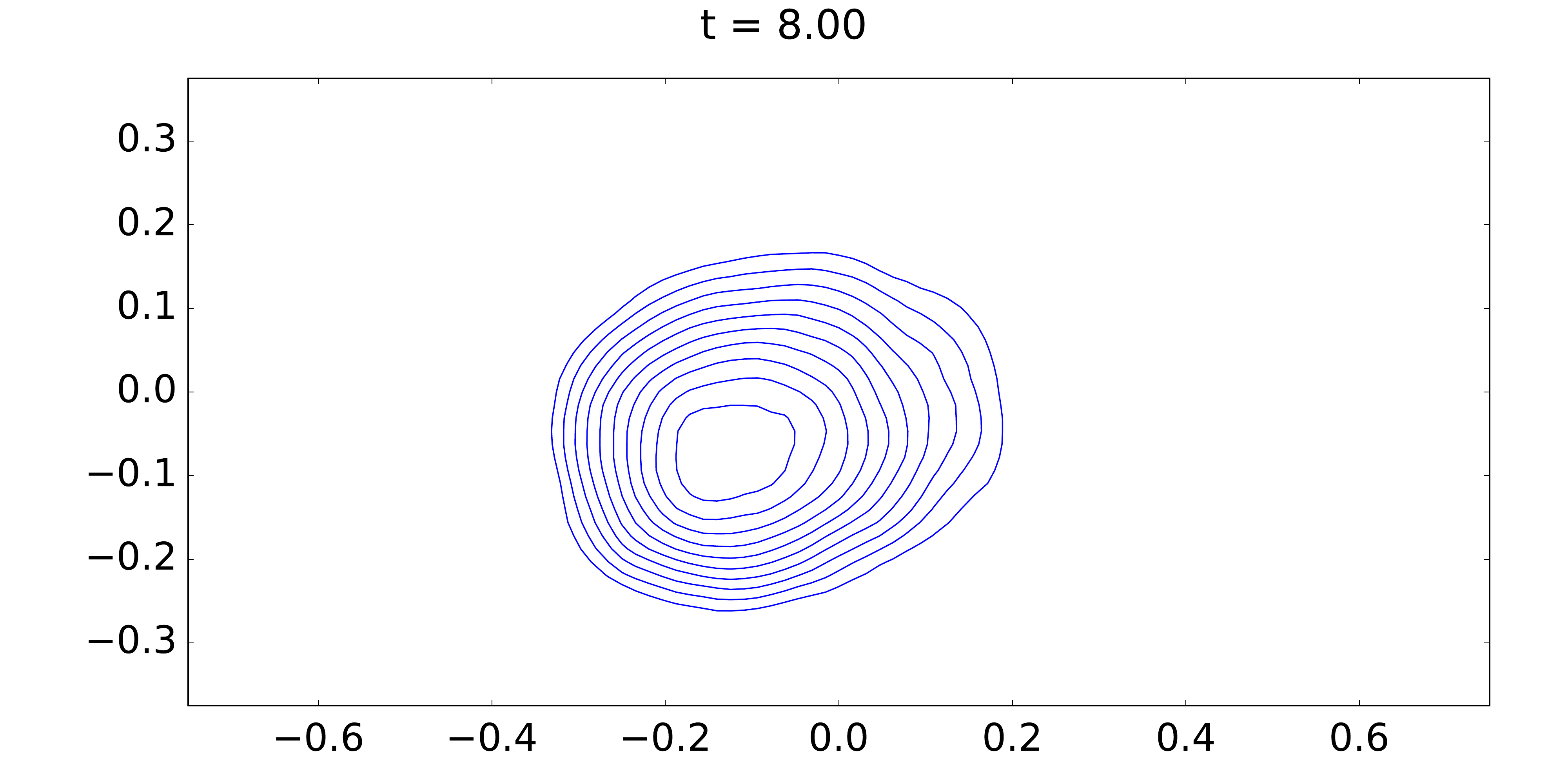}
\includegraphics[width=.4\textwidth]{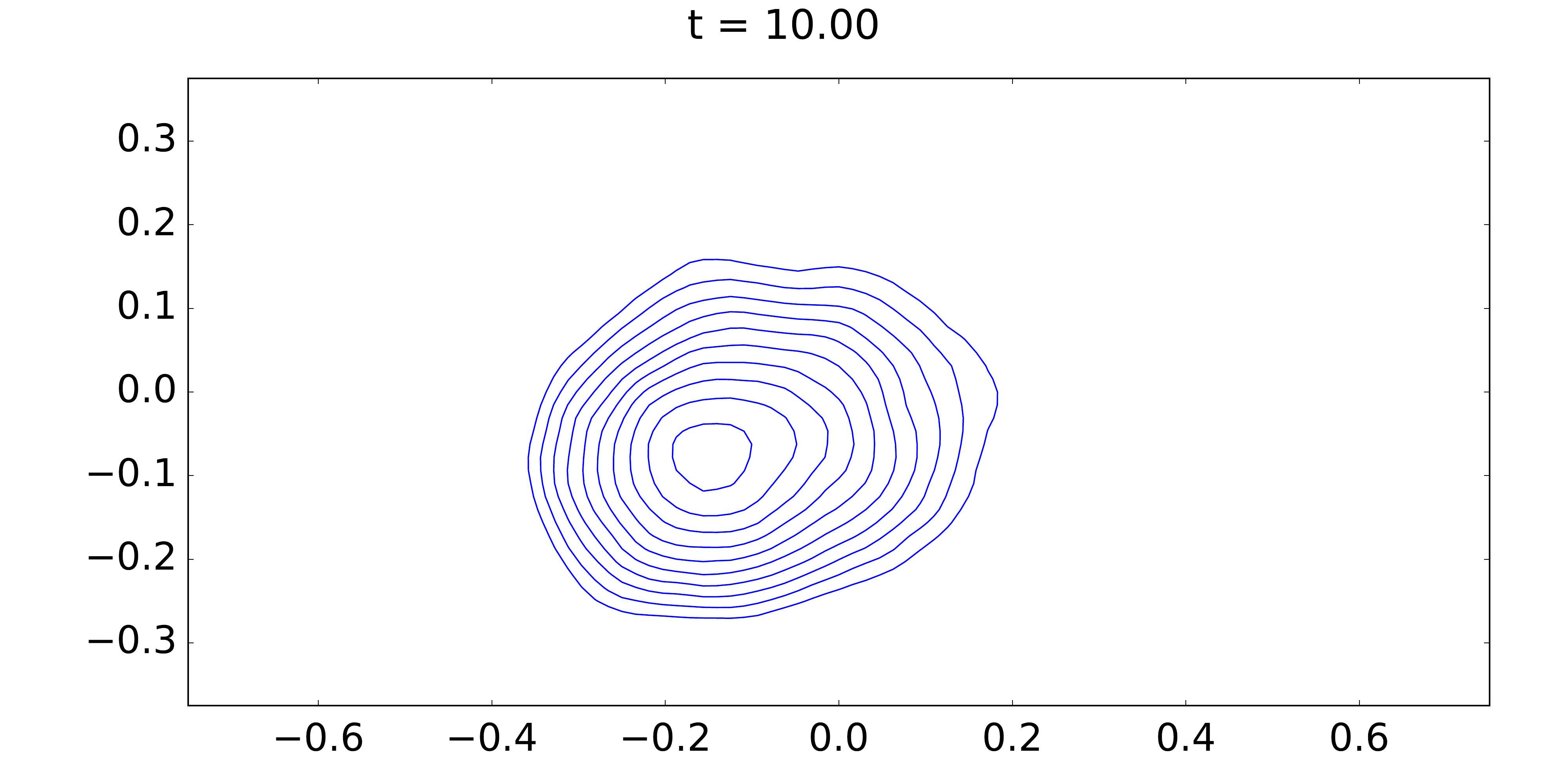}

\caption{Loop advection, non-smooth. Evolution of the magnetic field line contours. Only contours in the interval $[0.5 \times 10^{-4} , 2.8 \times 10^{-4}]$ are plotted.}
\label{fig:loop_advection_128x64_field_lines}
\end{figure}

The problem is set up so that the loop should return to its initial position after integer times $t = 1, 2, 3, ...$.
In Figure~\ref{fig:loop_advection_128x64_field_lines} it can be seen that this is initially the case, but after some time, the loop gets deformed, such that its centre is slightly displaced from its initial position at integer times.
One possible reason for this behaviour is that the proposed scheme is of low order and can be affected by phase errors as already observed for the example of Alfv\'{e}n Waves in Section~\ref{sec:ex_alfven_waves}.
Another possible reason is that we are using a finite difference discretisation, which assumes a sufficient degree of smoothness of the solution. Here, however, the magnetic field forms a cone and thus is only continuous.
Despite the quantitative deficits of the numerical solution, energy, magnetic helicity and cross helicity are preserved to machine accuracy throughout the whole simulation (c.f. Figure~\ref{fig:loop_advection_128x64_errors}). The magnetic pressure decreases slightly, but after $10$ passings the difference to the initial pressure is still smaller than $10^{-10}$.

\begin{figure}[p]
	\centering
	\subfloat[Initial Condition]{
		\includegraphics[width=.3\textwidth]{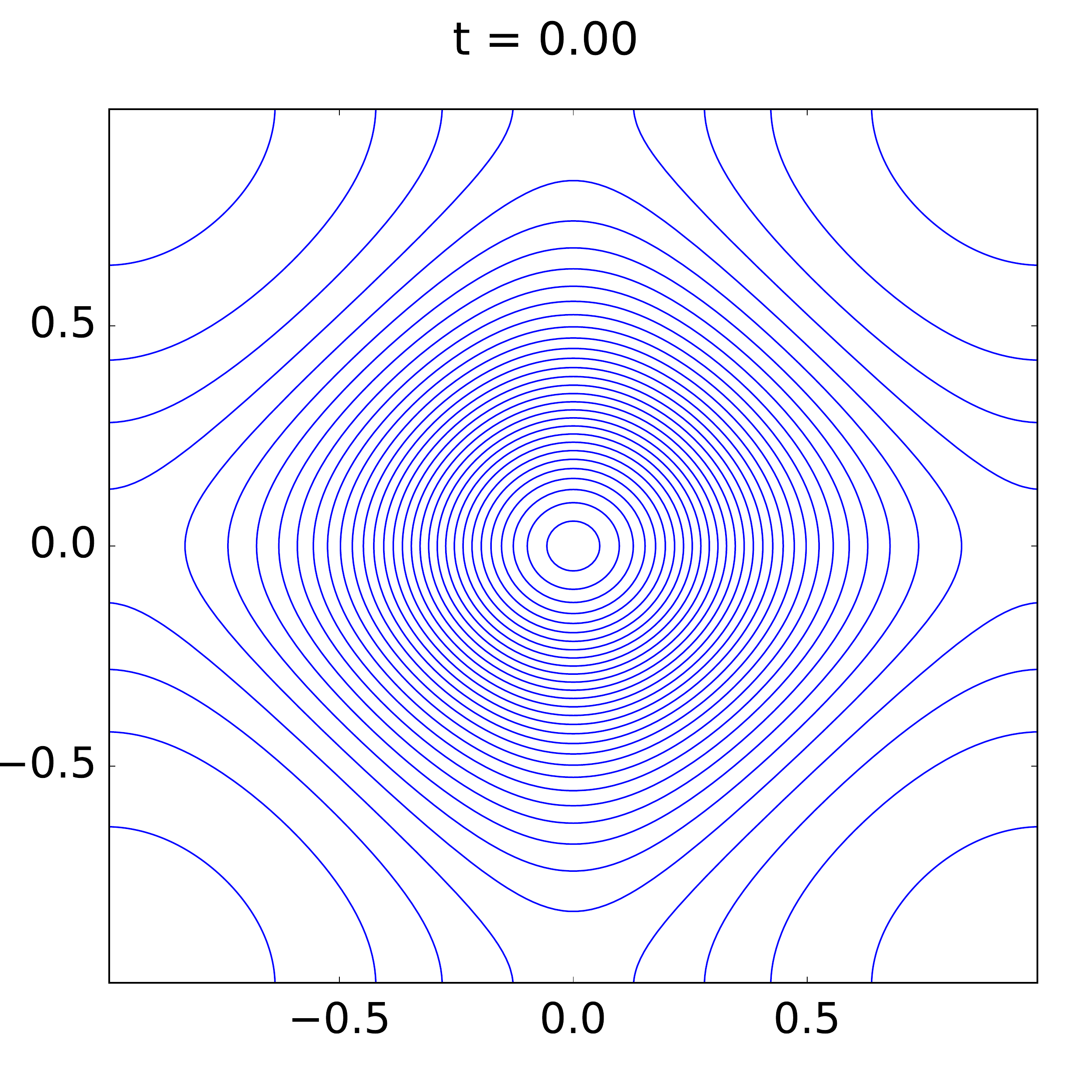}
	}
	\subfloat[$64 \times 64$ Points]{
		\includegraphics[width=.3\textwidth]{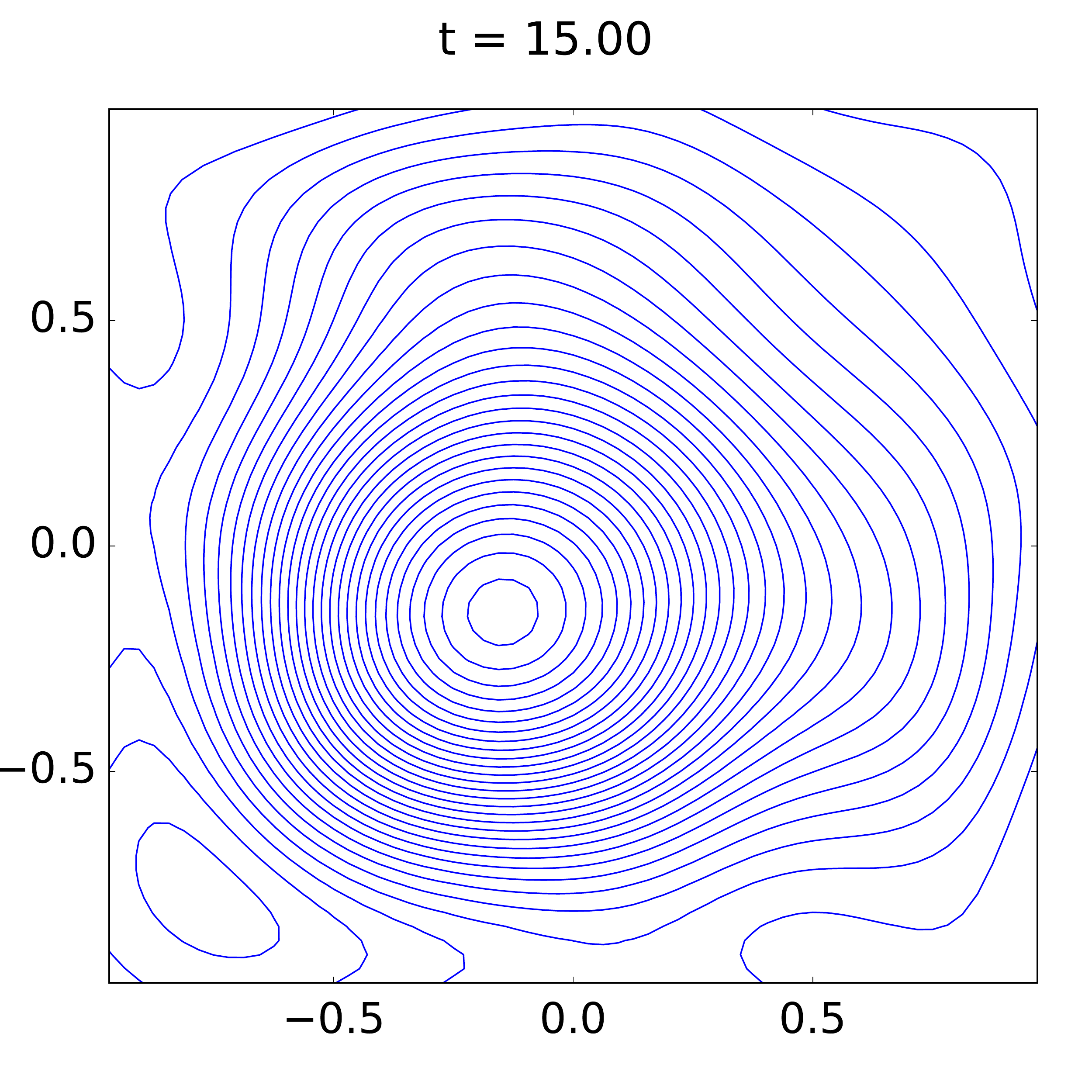}
	}
	
	\subfloat[$128 \times 128$ Points]{
		\includegraphics[width=.3\textwidth]{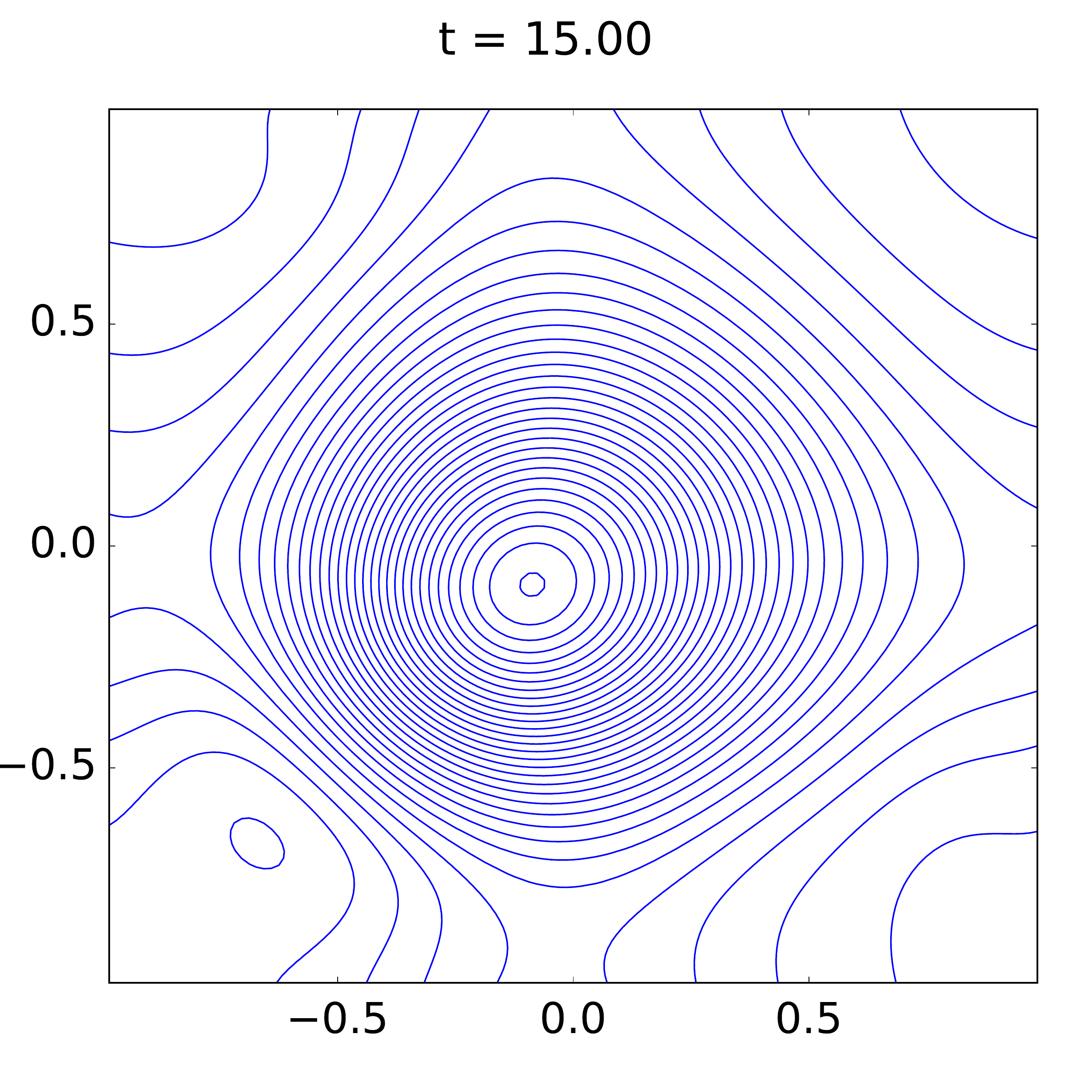}
	}
	\subfloat[$256 \times 256$ Points]{
		\includegraphics[width=.3\textwidth]{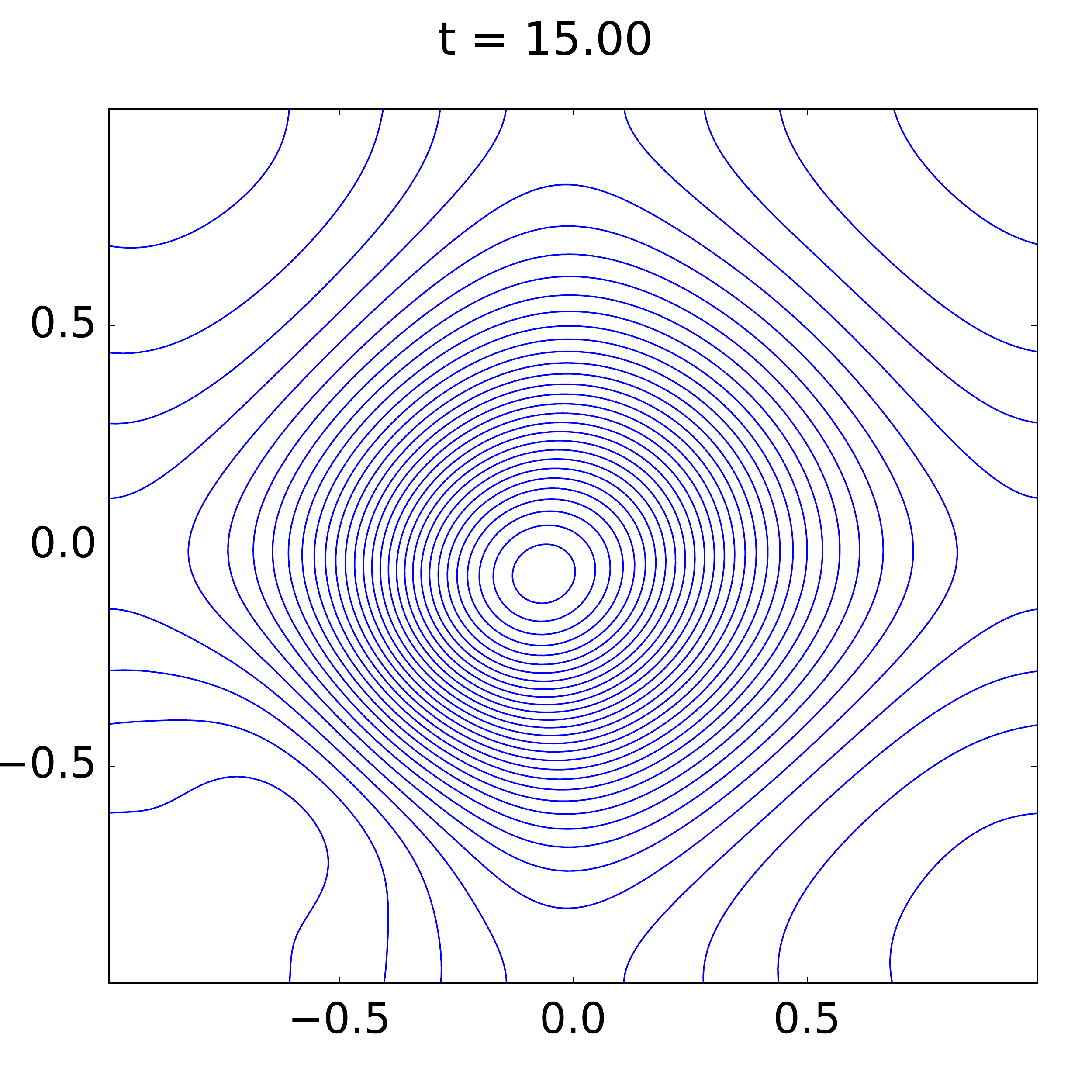}
	}

\caption{Loop advection, smooth, $128 \times 64$. Evolution of the magnetic field line contours.}
\label{fig:loop_advection_periodic_field_lines}
\end{figure}

In order to gain a better understanding of the origin of the deformation of the loop, we repeat the loop advection experiment with a different, smooth magnetic potential, namely
\begin{align*}
A &= A_{0} \exp \{ \cos (\pi x) + \cos (\pi y) \} .
\end{align*}
leading to a magnetic loop with smooth magnetic field
\begin{align*}
B^{x} &= \hphantom{-} \partial_{y} A = - A_{0} \pi \sin (\pi y) \exp \{ \cos (\pi x) + \cos (\pi y) \} , \\
B^{y} &= - \partial_{x} A = \hphantom{-} A_{0} \pi \sin (\pi x) \exp \{ \cos (\pi x) + \cos (\pi y) \} .
\end{align*}
The velocity and pressure are kept the same, but the parameters $V_{0}$ and $\theta$ are set to $V_{0} = \sqrt{8}$ and $\theta = \tan^{-1} (1)$ as the spatial domain is now $\Omega = [-1, +1] \times [-1, +1]$, still with periodic boundaries, but the lengths of the domain given by $L_x = 2$ and $L_y = 2$. We consider several spatial resolutions, namely $n_{x} \times n_{y} \in \{ 64 \times 64 , \, 128 \times 128 , \, 256 \times 256 \}$ grid points and the time step is $h_{t} = 0.01$.

The results are plotted in Figure~\ref{fig:loop_advection_periodic_field_lines}. We can see that the deformations at low resolution ($64 \times 64$) are comparable to those observed for the non-smooth loop, but become less pronounced with increasing resolution ($128 \times 128$ and $256 \times 256$). This hints at the low order of the scheme as the main reason for the inaccuracies.
However, in the smooth case, the deformations are mirror-symmetric with respect to the axis of advection, which is not the case for the non-smooth case. This can be attributed to the low regularity of the initial conditions of the non-smooth loop, which leads to more pronounced effects of dispersion than with smooth initial conditions.

\subsection{Current Sheet}\label{sec:ex_current_sheet}

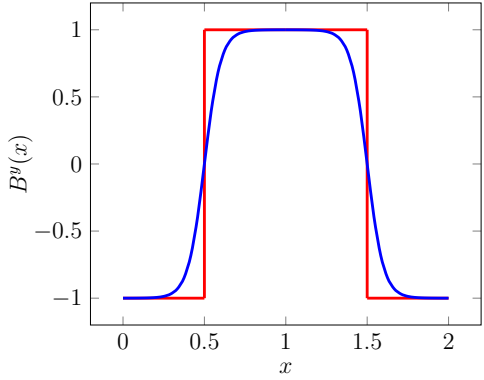
\begin{figure}[tb]
	\centering
	\begin{tikzpicture}[scale=.75]%
		\begin{axis}[domain=0:2,
					 xlabel=$x$,
					 ylabel=$B^{y} (x)$]
			\addplot [red, line width=.5mm, no markers] coordinates {(0.0,-1) (0.5,-1)};
			\addplot [red, line width=.5mm, no markers] coordinates {(0.5,-1) (0.5,+1)};
			\addplot [red, line width=.5mm, no markers] coordinates {(0.5,+1) (1.5,+1)};
			\addplot [red, line width=.5mm, no markers] coordinates {(1.5,+1) (1.5,-1)};
			\addplot [red, line width=.5mm, no markers] coordinates {(1.5,-1) (2.0,-1)};
			\addplot[domain=0:1, color=blue, line width=.5mm, smooth] {+tanh(10*(x-0.5))};
			\addplot[domain=1:2, color=blue, line width=.5mm, smooth] {-tanh(10*(x-1.5))};
		\end{axis}
	\end{tikzpicture}
	\caption{Current sheet models: discontinuous magnetic field (red), smooth field (blue).}
	\label{fig:current_sheet_models}
\end{figure}

In the following, we consider as initial conditions for the magnetic field two different current sheet models that appear in reconnection studies (see Figure~\ref{fig:current_sheet_models}). That is a discontinuous magnetic field with a sharp jump~\cite{GardinerStone:2005}, caused by two singular current sheets,
\begin{align*}
B^{y}_{\text{sharp}} &=
\begin{cases}
+ 1 & x < x_{1} \\
- 1 & x_{1} \leq x \leq x_{2} \\
+ 1 & x > x_{2}
\end{cases} &
& \text{with} &
x_{1} &= 0.5 , &
x_{2} &= 1.5 , 
\end{align*}
and a $\tanh$ profile similar to \cite{GrassoCalifano:2001}, caused by current sheets with finite thickness,
\begin{align*}
B^{y}_{\tanh} &= \begin{cases}
+ \tanh(10 (x-x_{1})) & x < 1 \\
- \tanh(10 (x-x_{2})) & 1 \geq x
\end{cases} &
& \text{with} &
x_{1} &= 0.5 , &
x_{2} &= 1.5 .
\end{align*}
In both cases we have $B^{x} = 0$ and the following initial conditions for the fluid,
\begin{align*}
V^{x} &= V_{0} \, \sin (\pi y), &
V^{y} &= 0, &
P &= 0.1 .
\end{align*}
The spatial domain is $\Omega = [0,2] \times [0,2]$, which is discretised by $n_{x} \times n_{y} = 32 \times 32$ grid points. We use periodic boundaries and a time step of $h_{t} = 0.1$.

\begin{figure}[bt]
\centering
\subfloat[Singular Current Sheet]{\label{fig:current_sheet_32x32_errors}
\begin{minipage}{.48\textwidth}
\includegraphics[width=\textwidth]{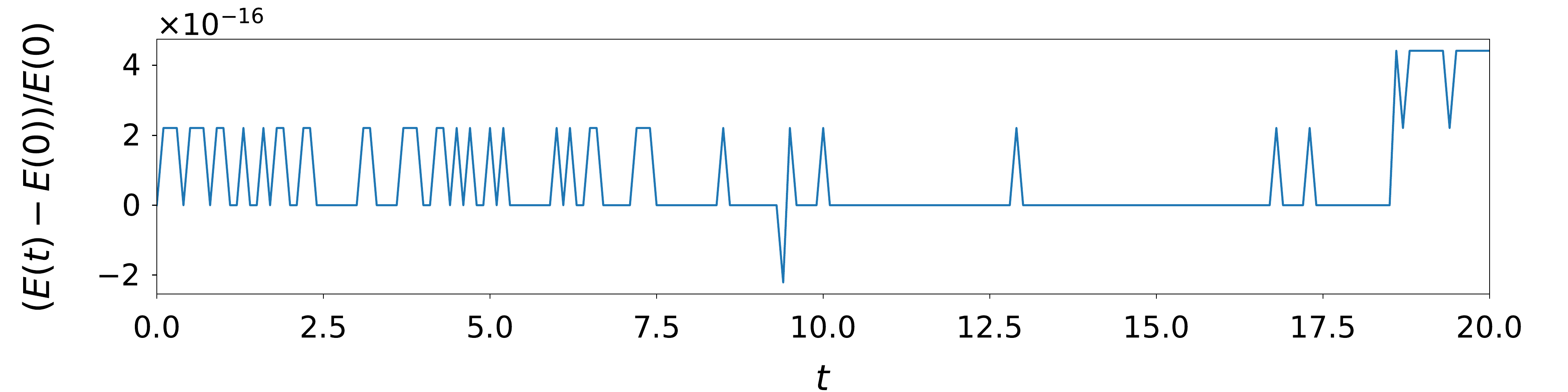}

\includegraphics[width=\textwidth]{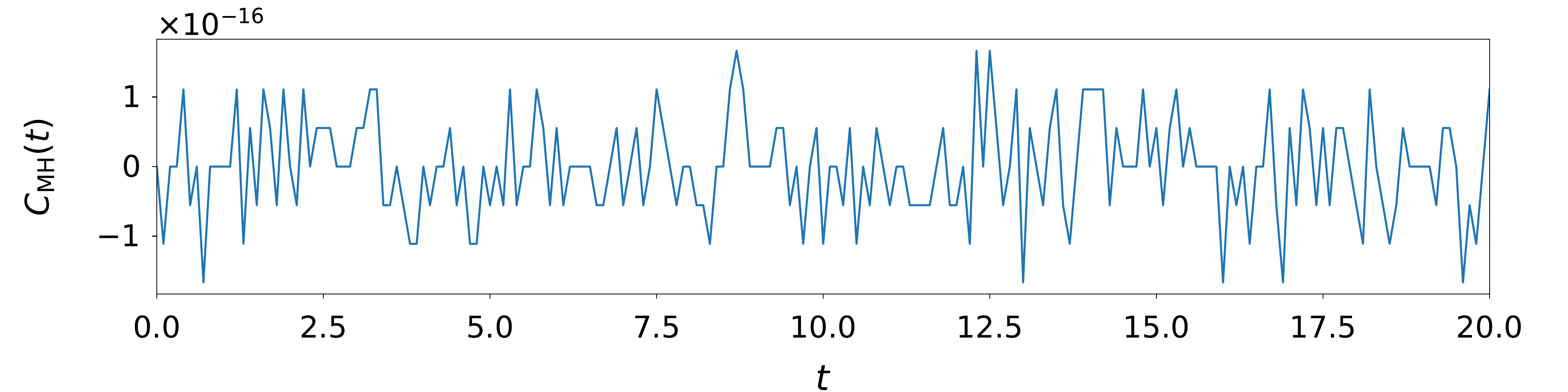}

\includegraphics[width=\textwidth]{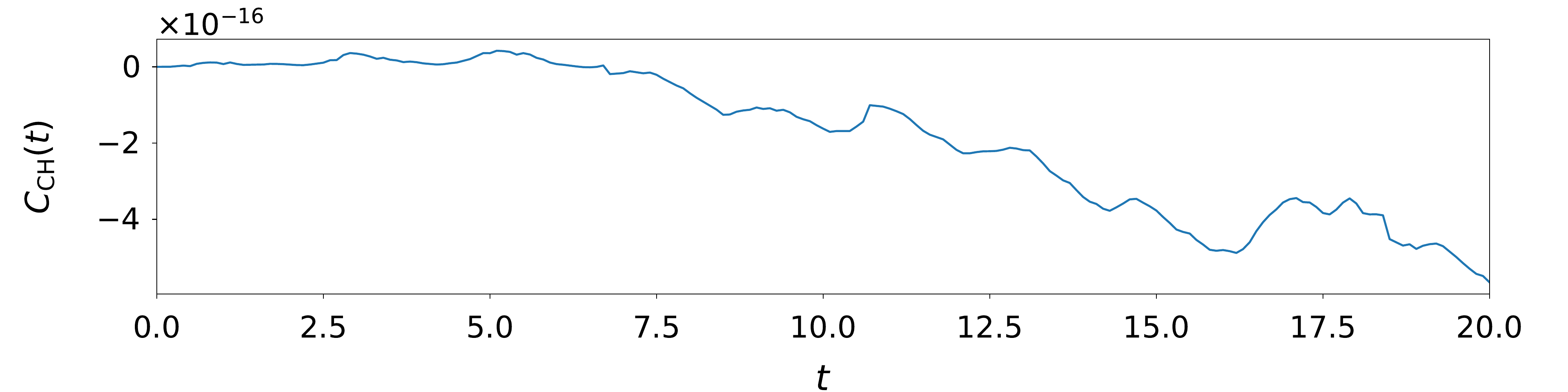}
\end{minipage}
}
\subfloat[Smooth Current Sheet]{\label{fig:current_sheet_tanh_errors}
\begin{minipage}{.48\textwidth}
\includegraphics[width=\textwidth]{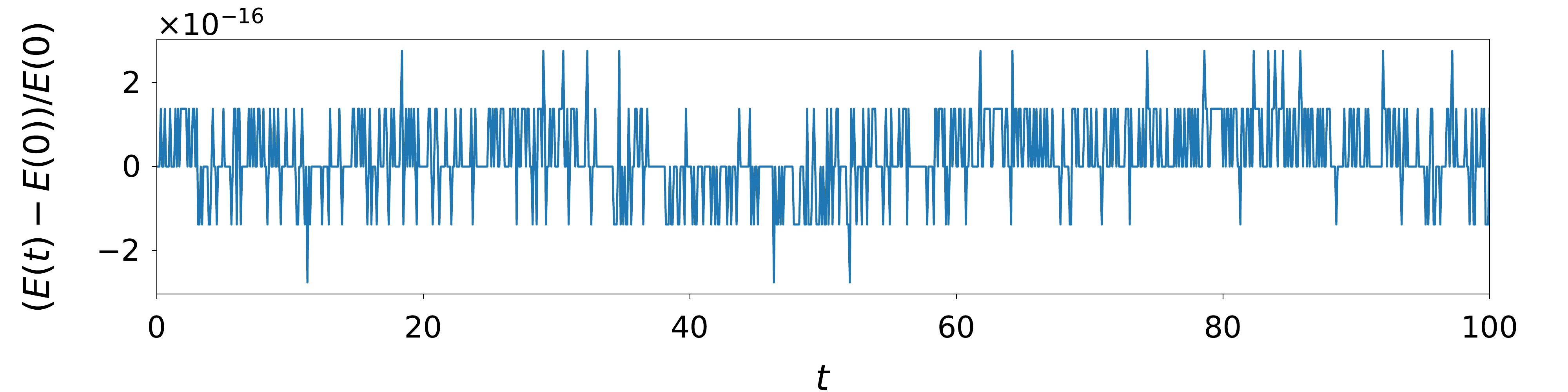}

\includegraphics[width=\textwidth]{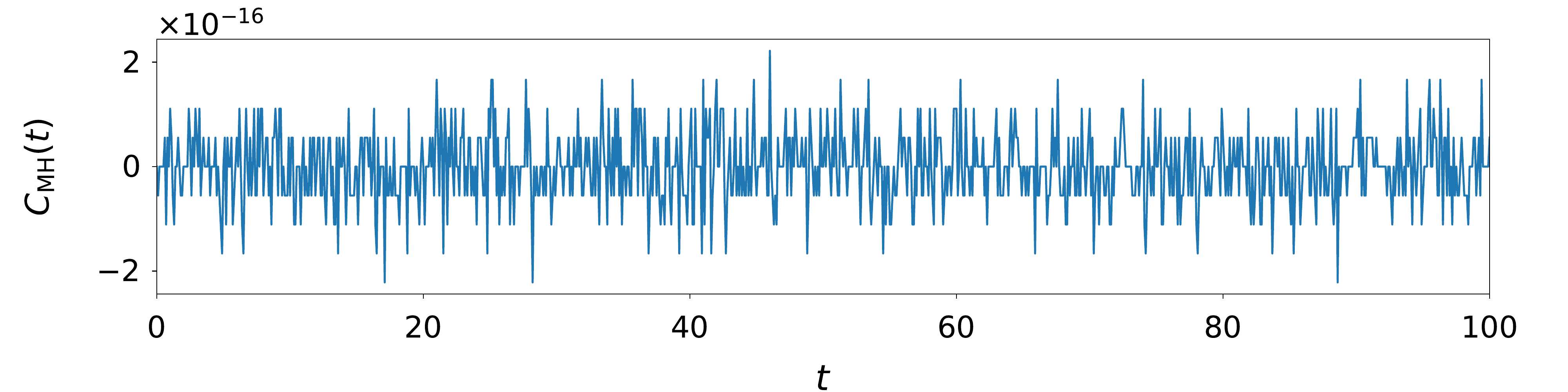}

\includegraphics[width=\textwidth]{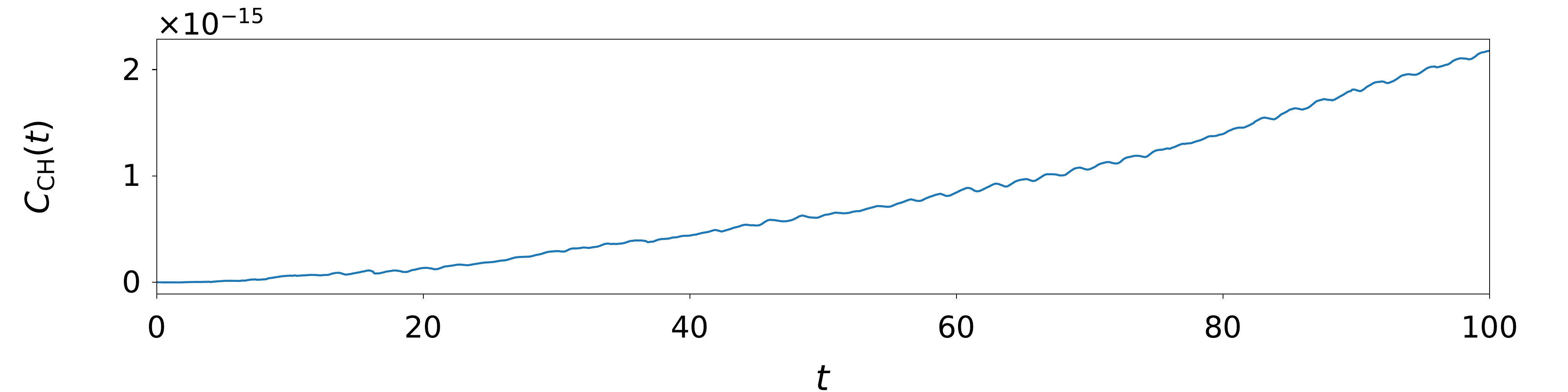}
\end{minipage}
}
\caption{Conservation of energy, magnetic and cross helicity for singular current sheet (left) and smooth current sheet (right).}
\end{figure}

In both cases, energy and magnetic helicity are exactly preserved, i.e., up to machine accuracy (see Figures~\ref{fig:current_sheet_32x32_errors}~-~\ref{fig:current_sheet_tanh_errors}).
For cross helicity, we see a slight drift (Figure~\ref{fig:current_sheet_tanh_errors}), however, after 100 characteristic times the error is still of order $10^{-15}$.
In the following, we want to focus on the conservation of field line topology.

\begin{figure}[p]
\centering
\subfloat{
\includegraphics[width=.32\textwidth]{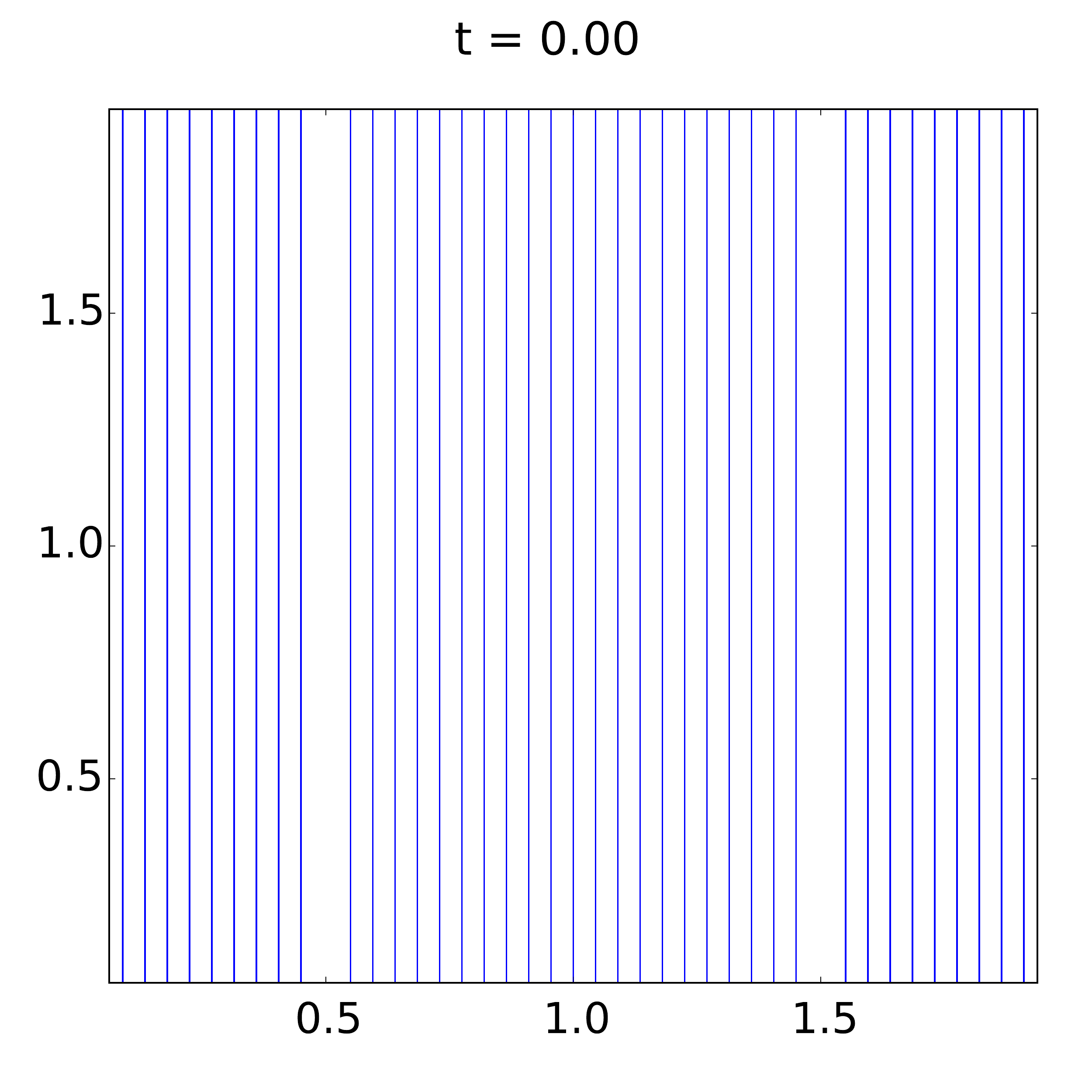}
}
\subfloat{
\includegraphics[width=.32\textwidth]{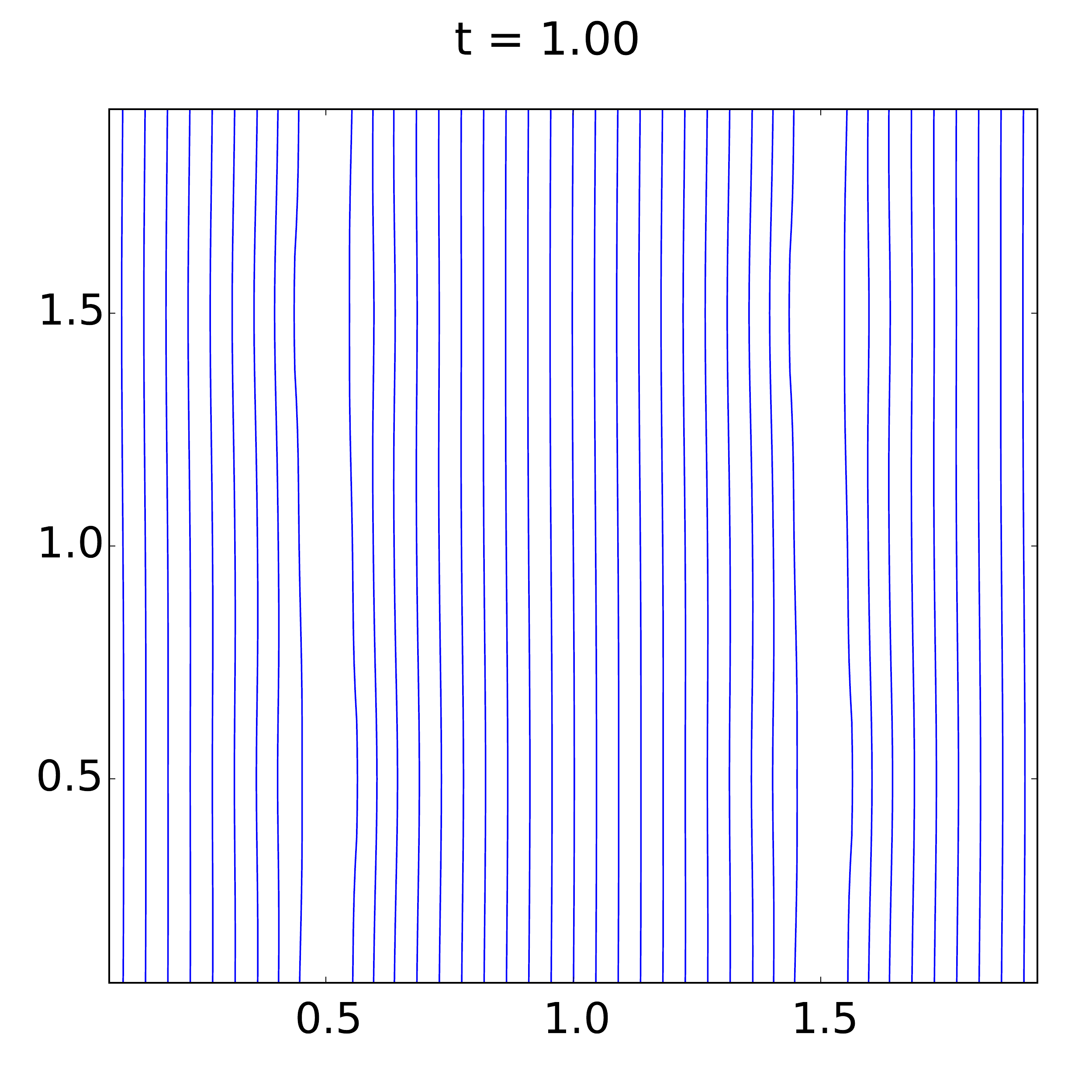}
}
\subfloat{
\includegraphics[width=.32\textwidth]{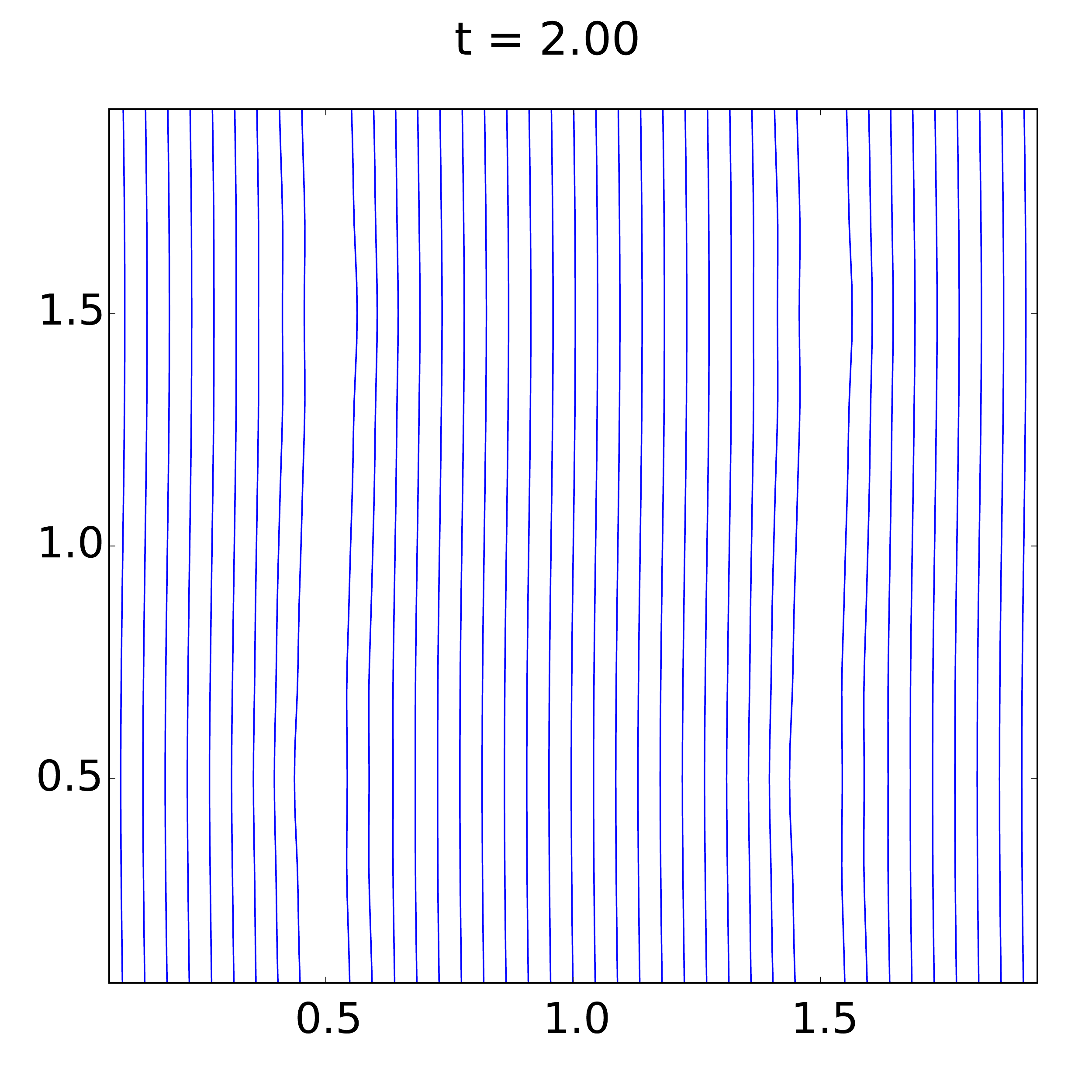}
}

\subfloat{
\includegraphics[width=.32\textwidth]{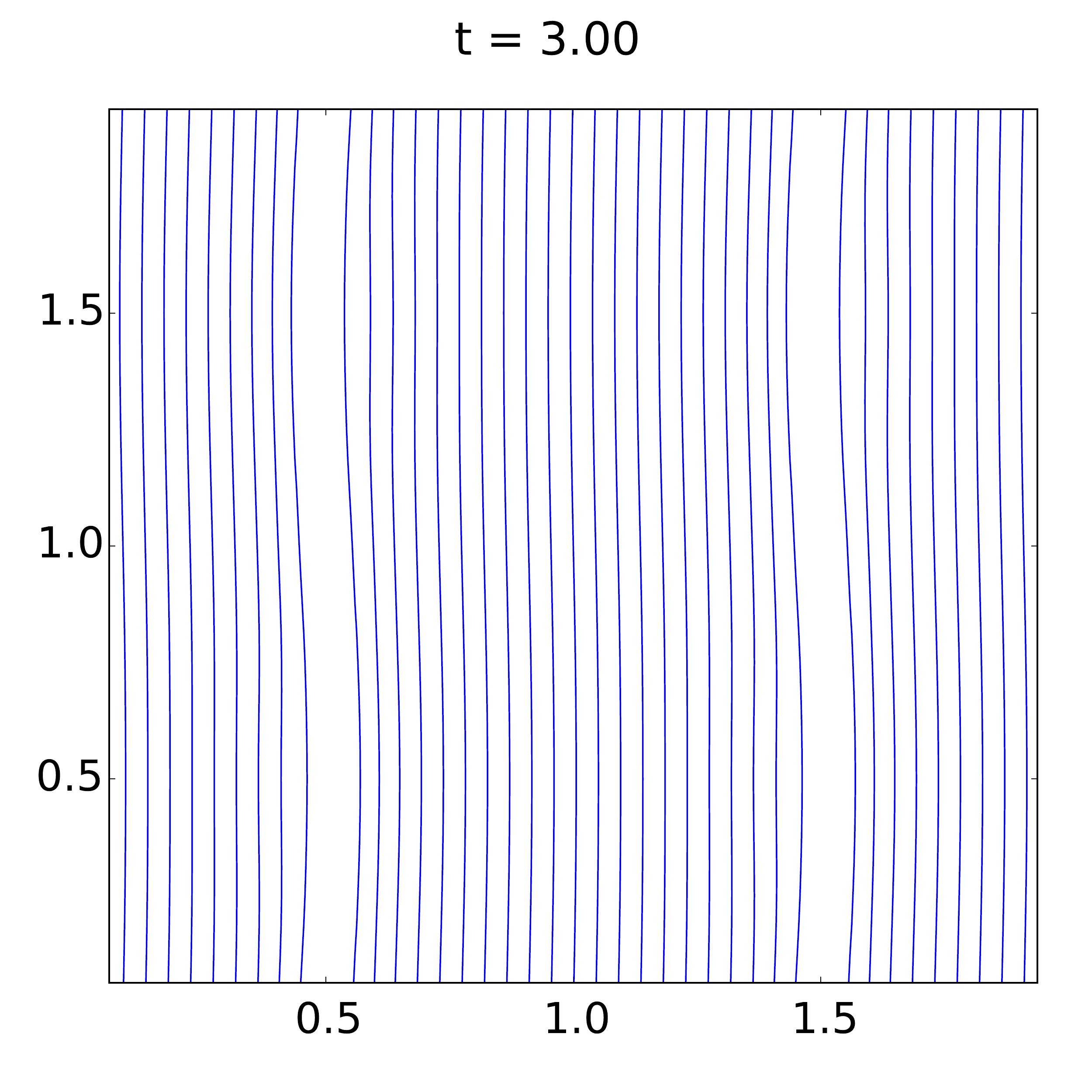}
}
\subfloat{
\includegraphics[width=.32\textwidth]{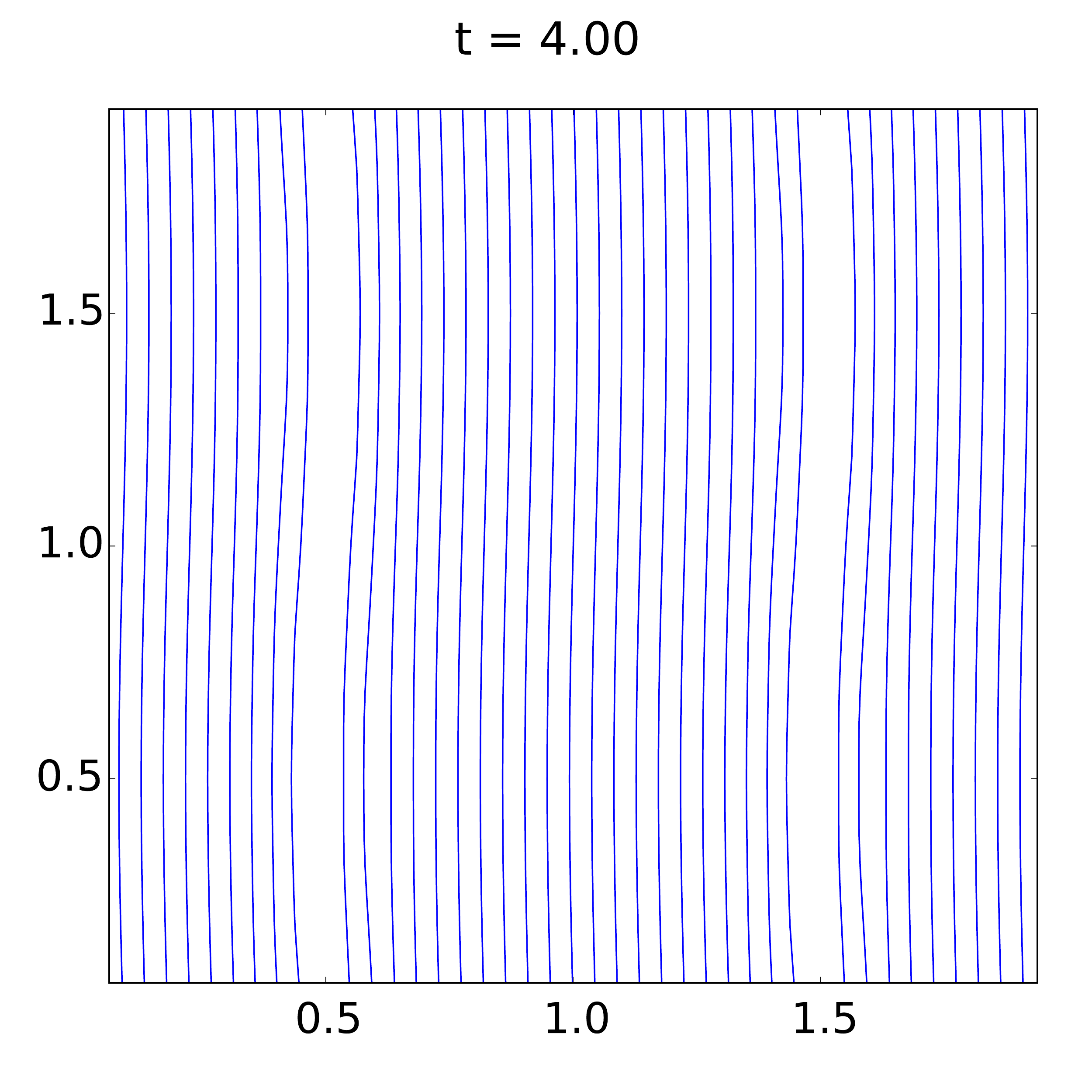}
}
\subfloat{
\includegraphics[width=.32\textwidth]{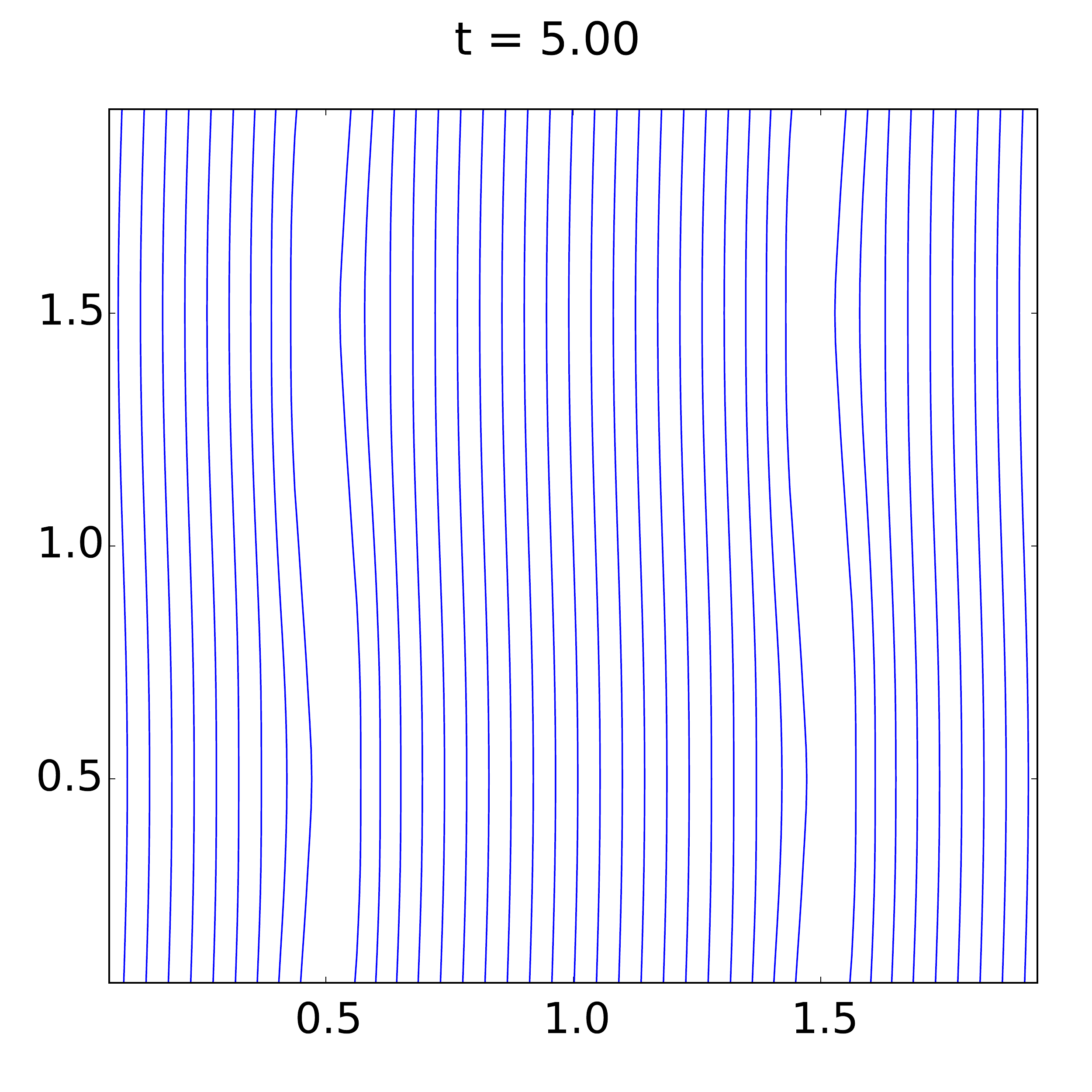}
}

\subfloat{
\includegraphics[width=.32\textwidth]{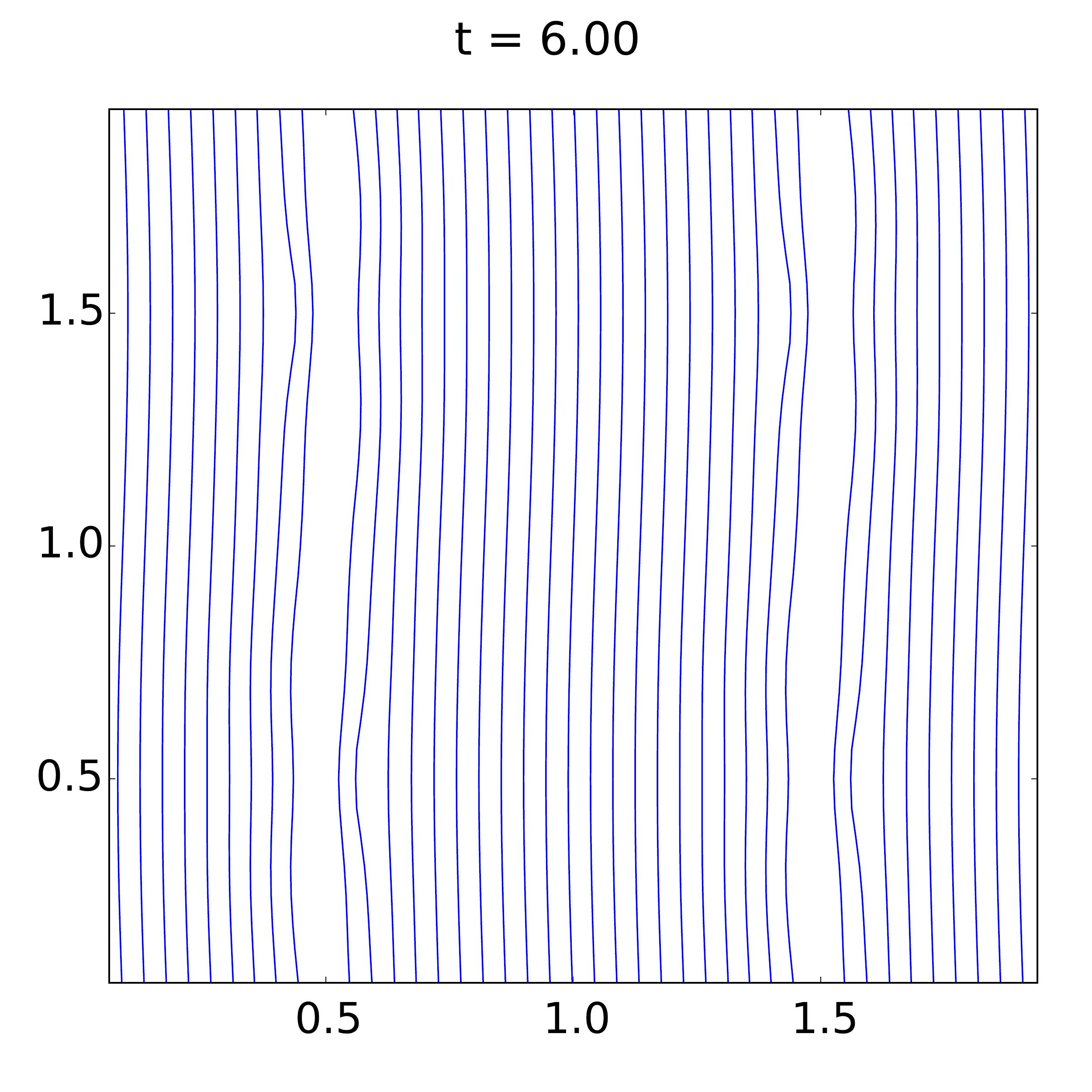}
}
\subfloat{
\includegraphics[width=.32\textwidth]{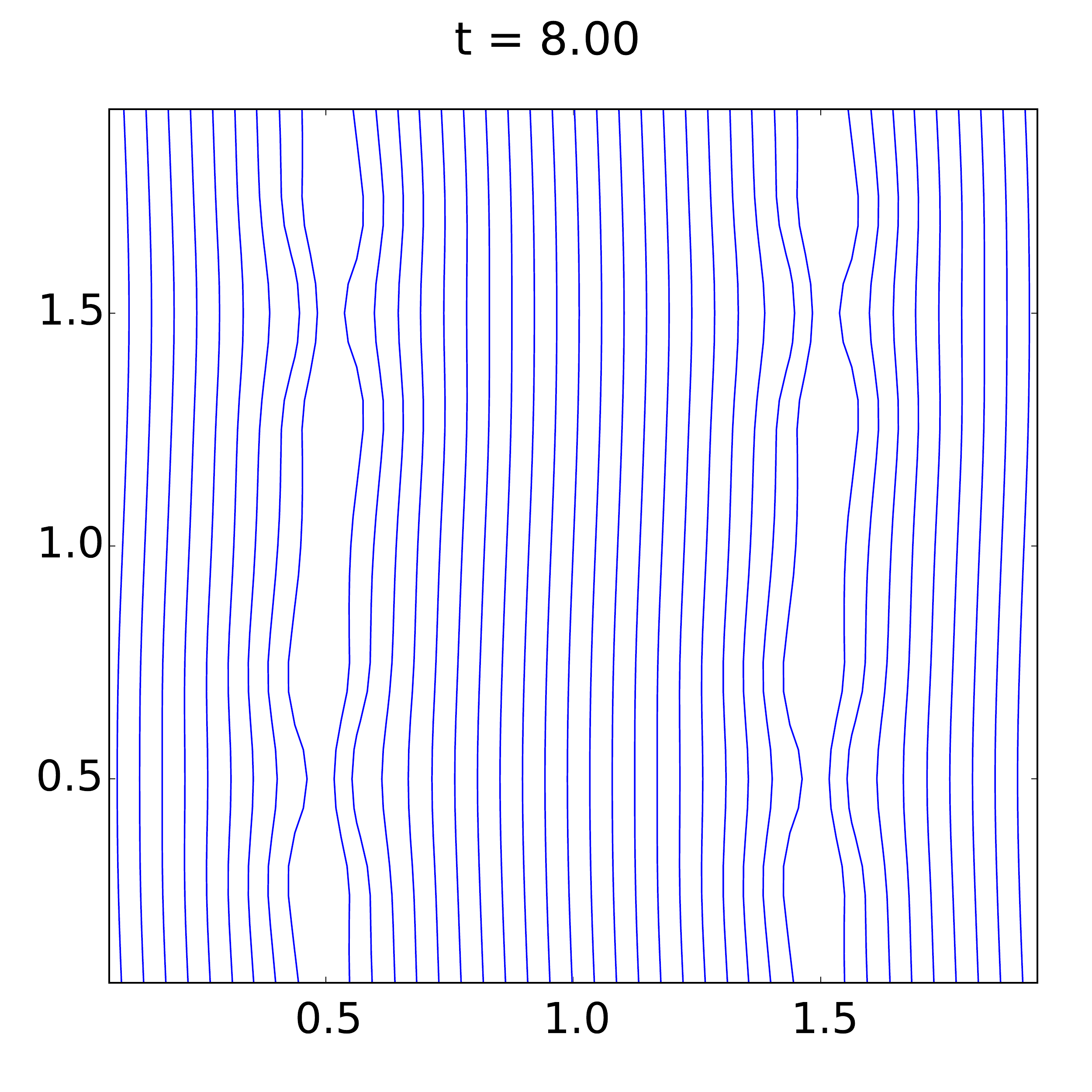}
}
\subfloat{
\includegraphics[width=.32\textwidth]{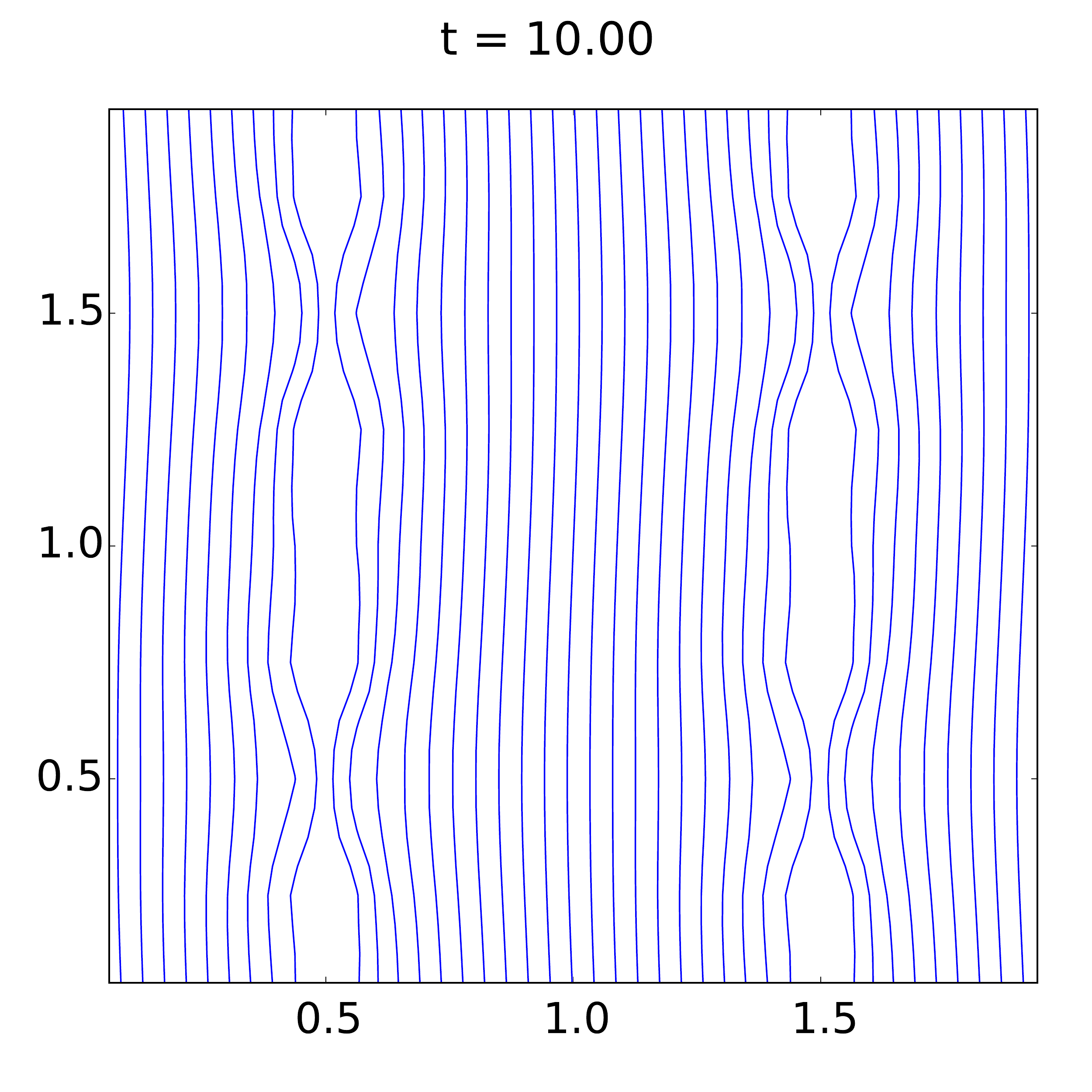}
}

\subfloat{
\includegraphics[width=.32\textwidth]{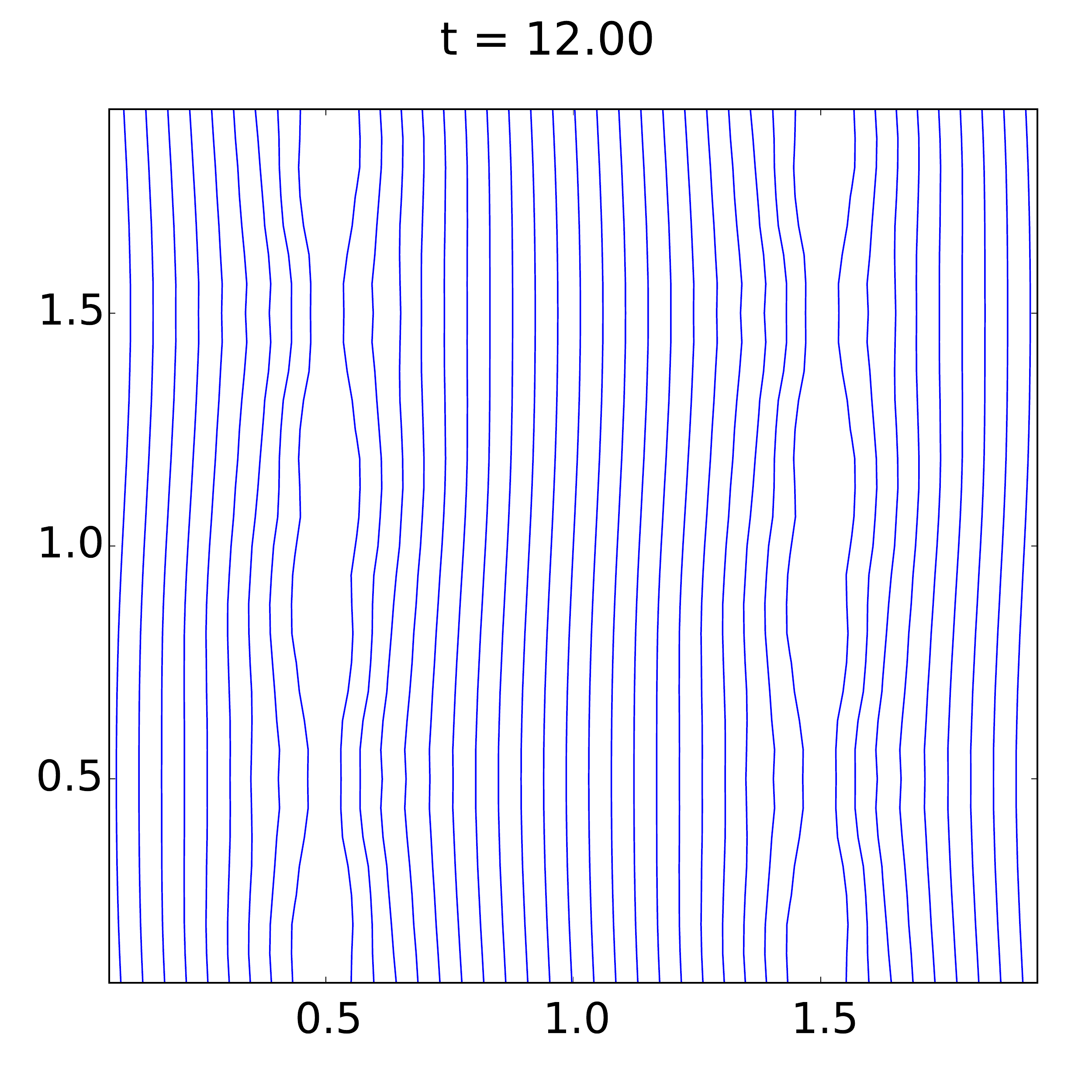}
}
\subfloat{
\includegraphics[width=.32\textwidth]{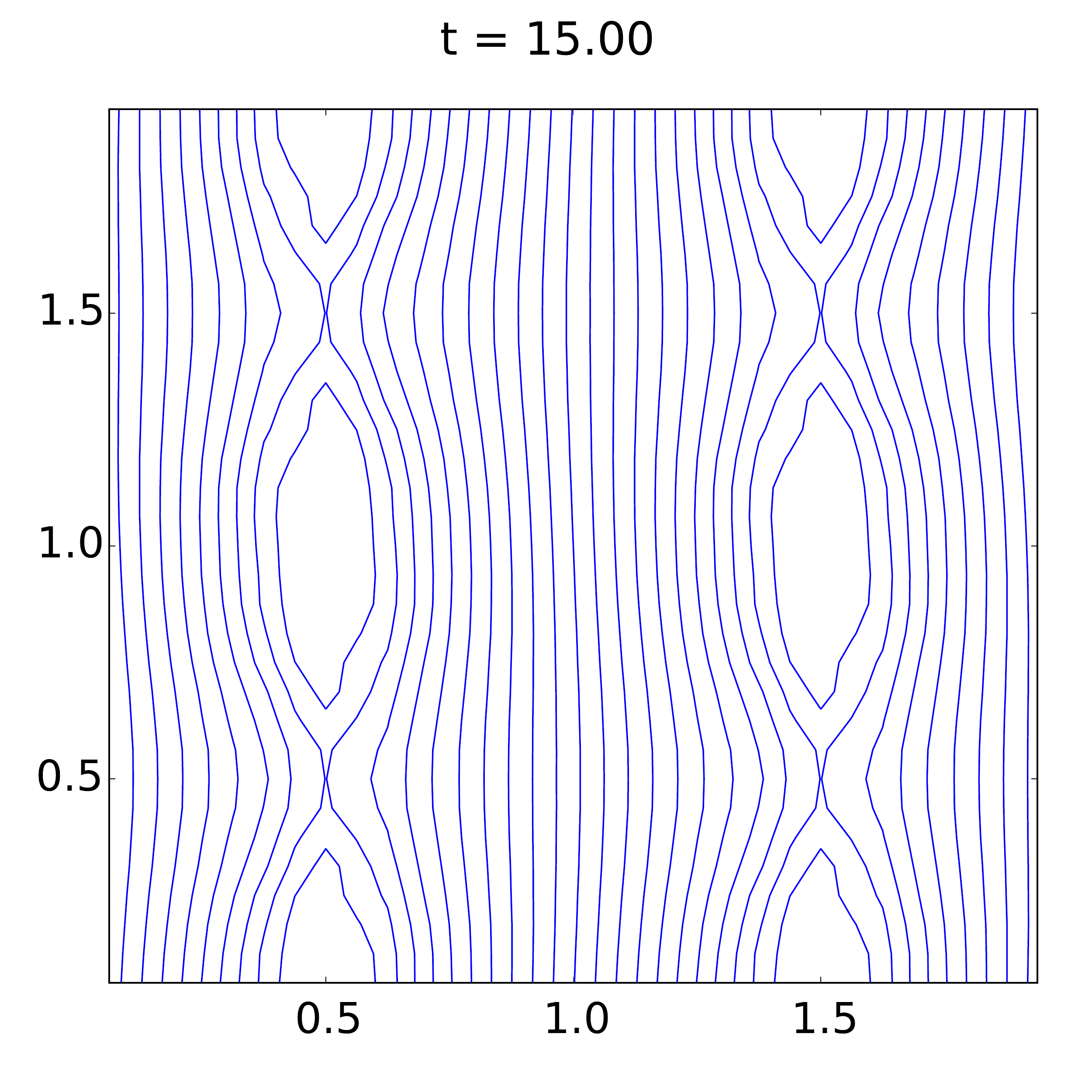}
}
\subfloat{
\includegraphics[width=.32\textwidth]{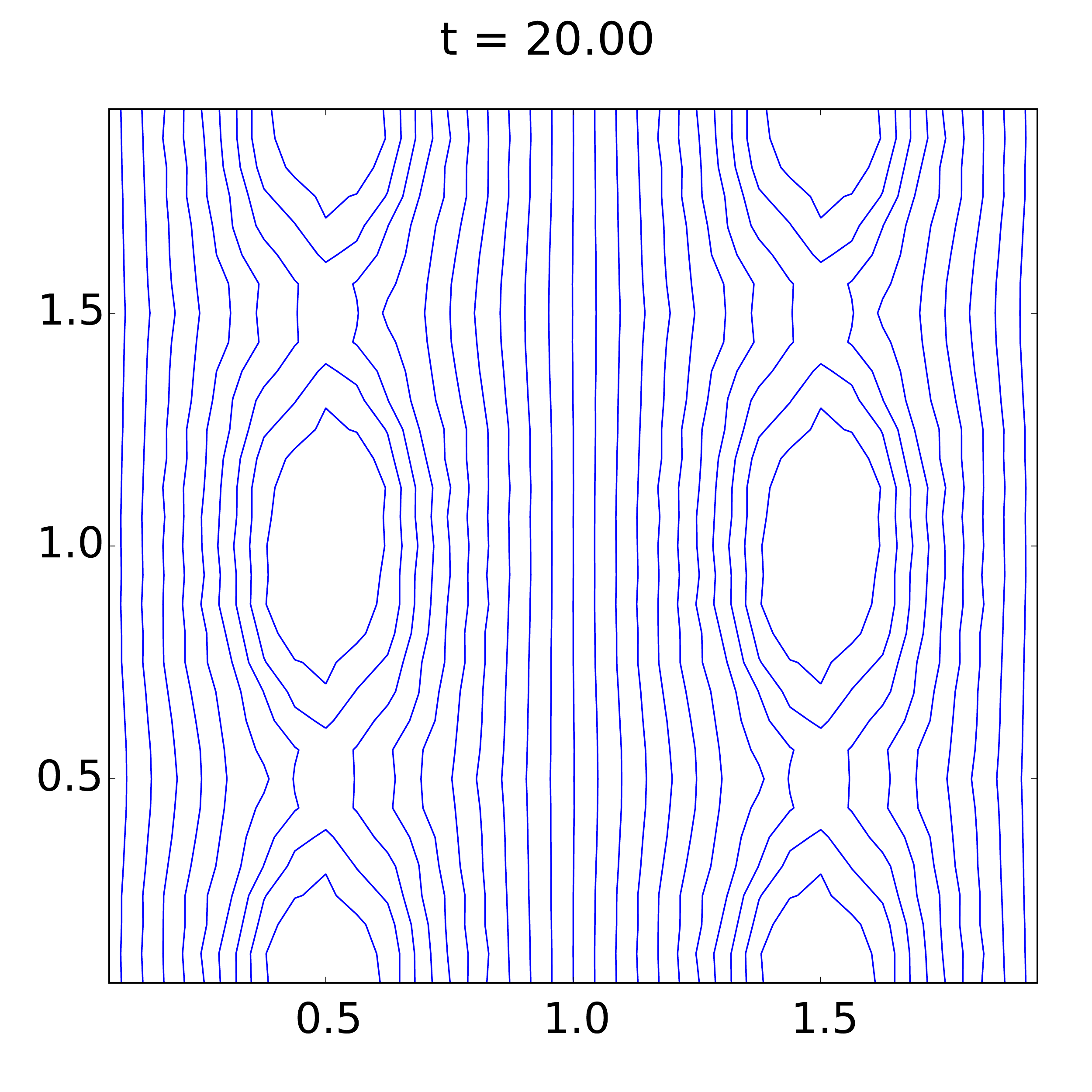}
}

\caption{Singular current sheet. Magnetic field lines.}
\label{fig:current_sheet_field_lines}
\end{figure}

In Figure~\ref{fig:current_sheet_field_lines}, the field line evolution for the discontinuous magnetic field is plotted. Initially all field lines are parallel. Due to the perturbation in the velocity field, the magnetic field lines get bend, but for more then 10 characteristic times, they do not break up and reconnect.
After $t = 12$, however, magnetic islands start to form and consecutively grow as can be seen at $t = 20$. At this point, the solution can not be regarded as physical anymore.
We have to stress here, that this set of initial conditions is quite challenging for most numerical schemes due to the discontinuity, and that with other methods reconnections sets in much earlier, e.g., in the range $t = 0.5 \hdots 1.0$ for the Gudonov scheme of~\citet{GardinerStone:2005}.

\begin{figure}[p]
\centering
\subfloat{
\includegraphics[width=.32\textwidth]{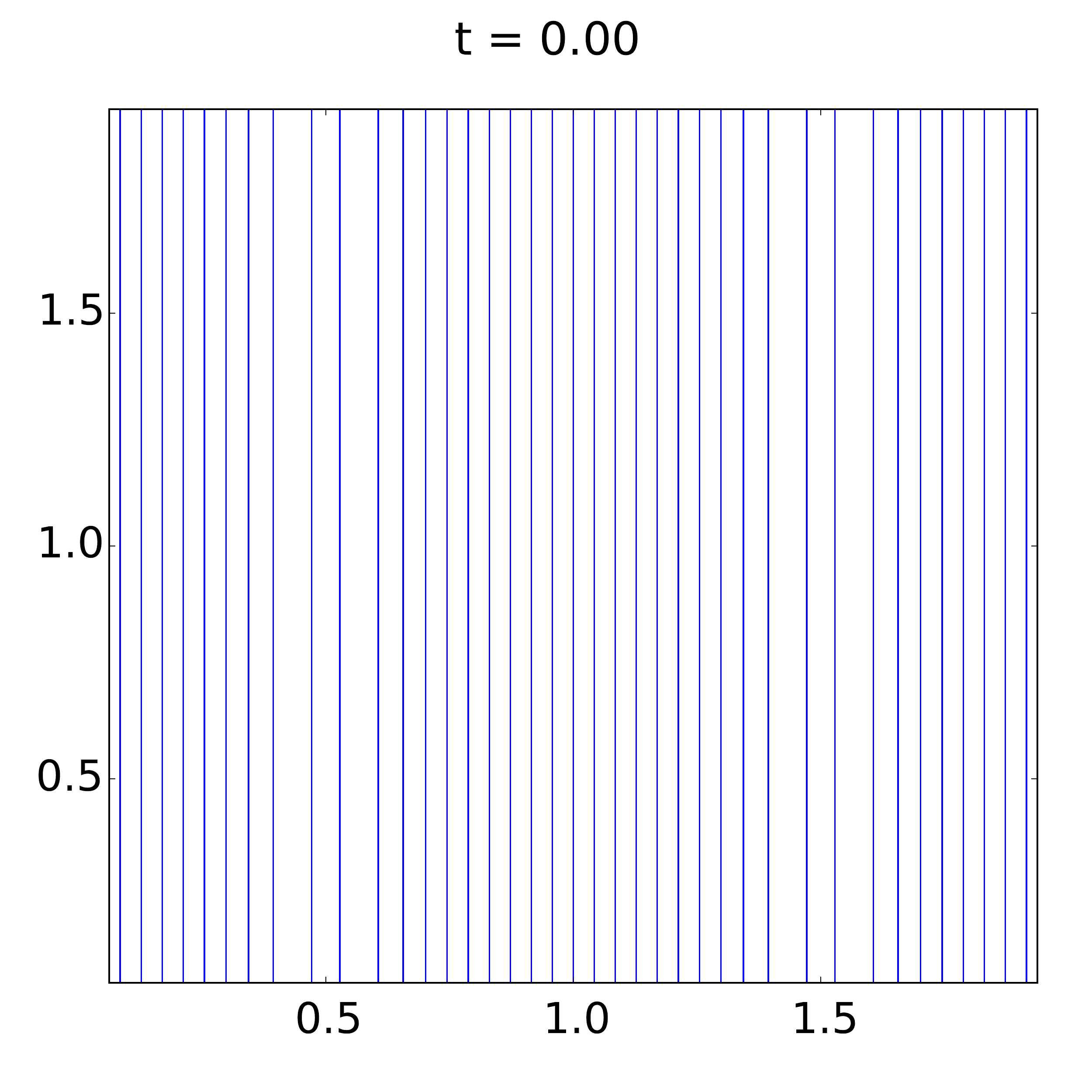}
}
\subfloat{
\includegraphics[width=.32\textwidth]{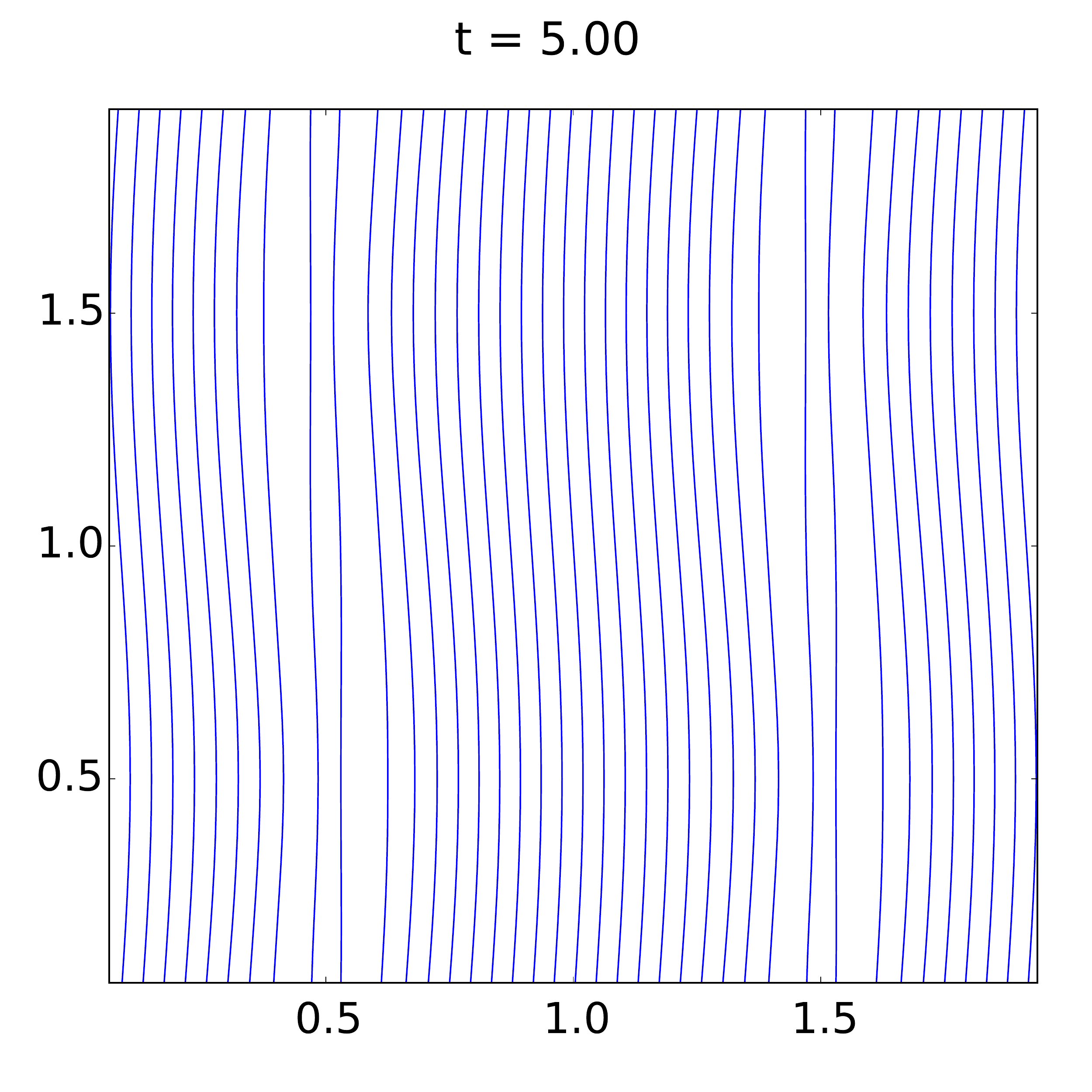}
}
\subfloat{
\includegraphics[width=.32\textwidth]{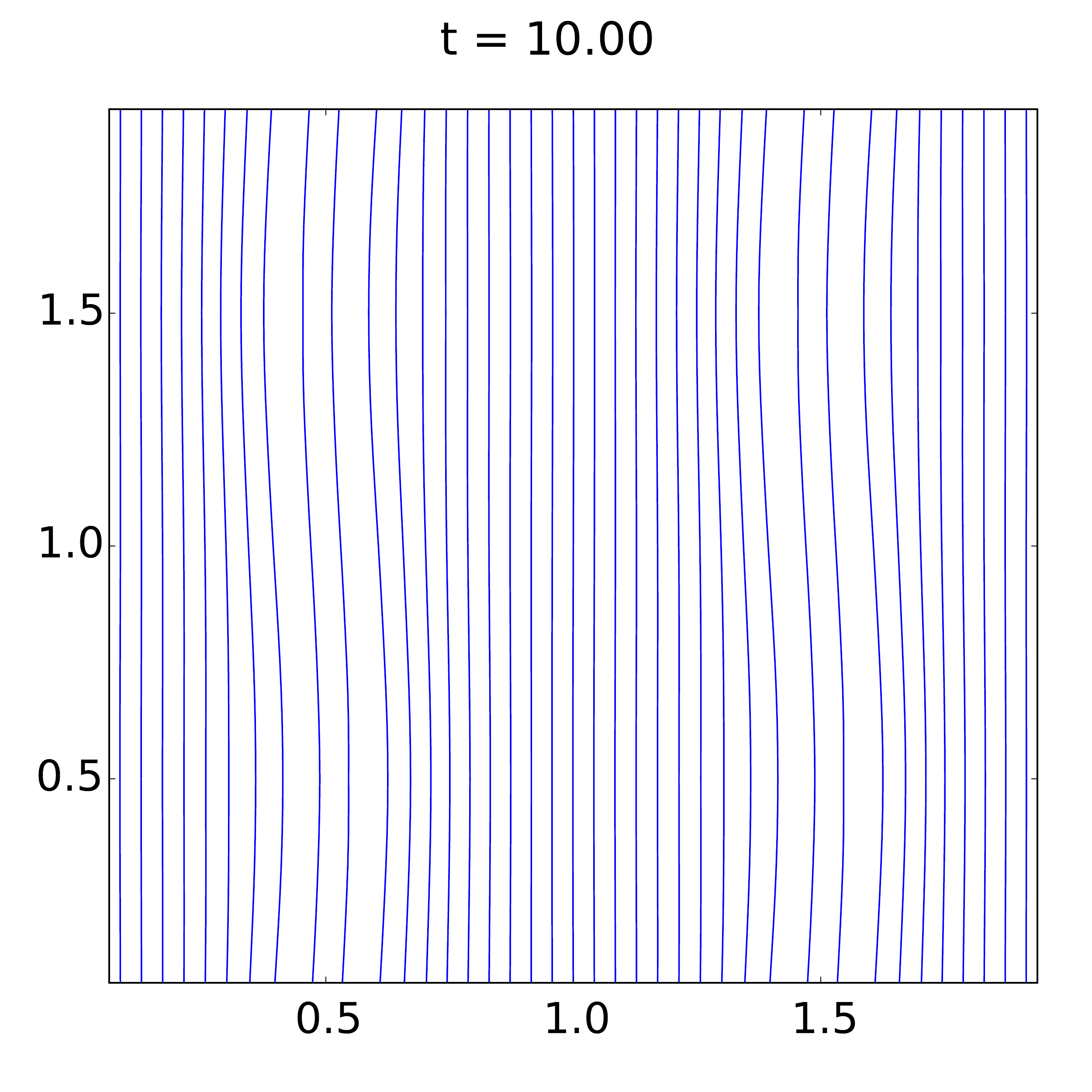}
}

\subfloat{
\includegraphics[width=.32\textwidth]{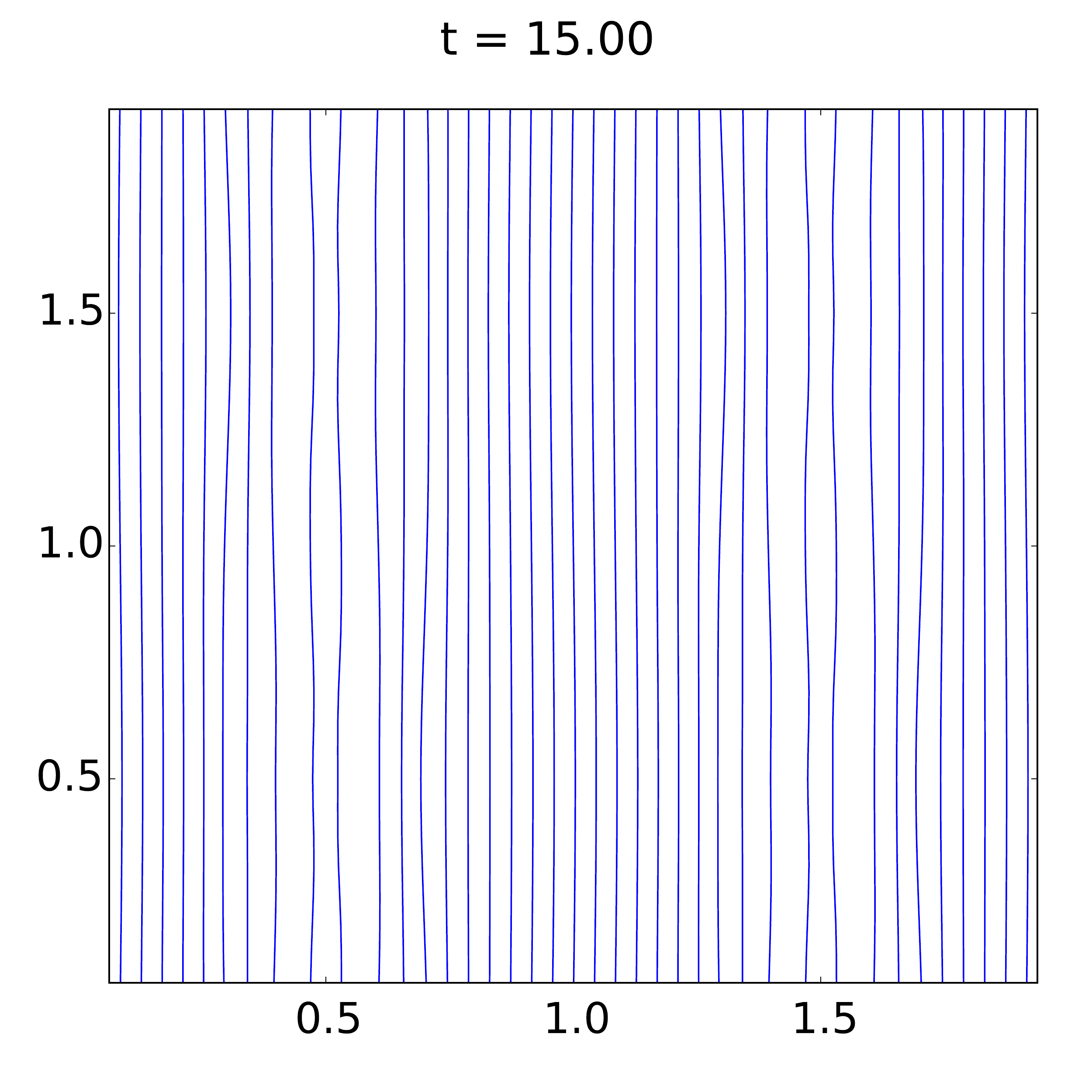}
}
\subfloat{
\includegraphics[width=.32\textwidth]{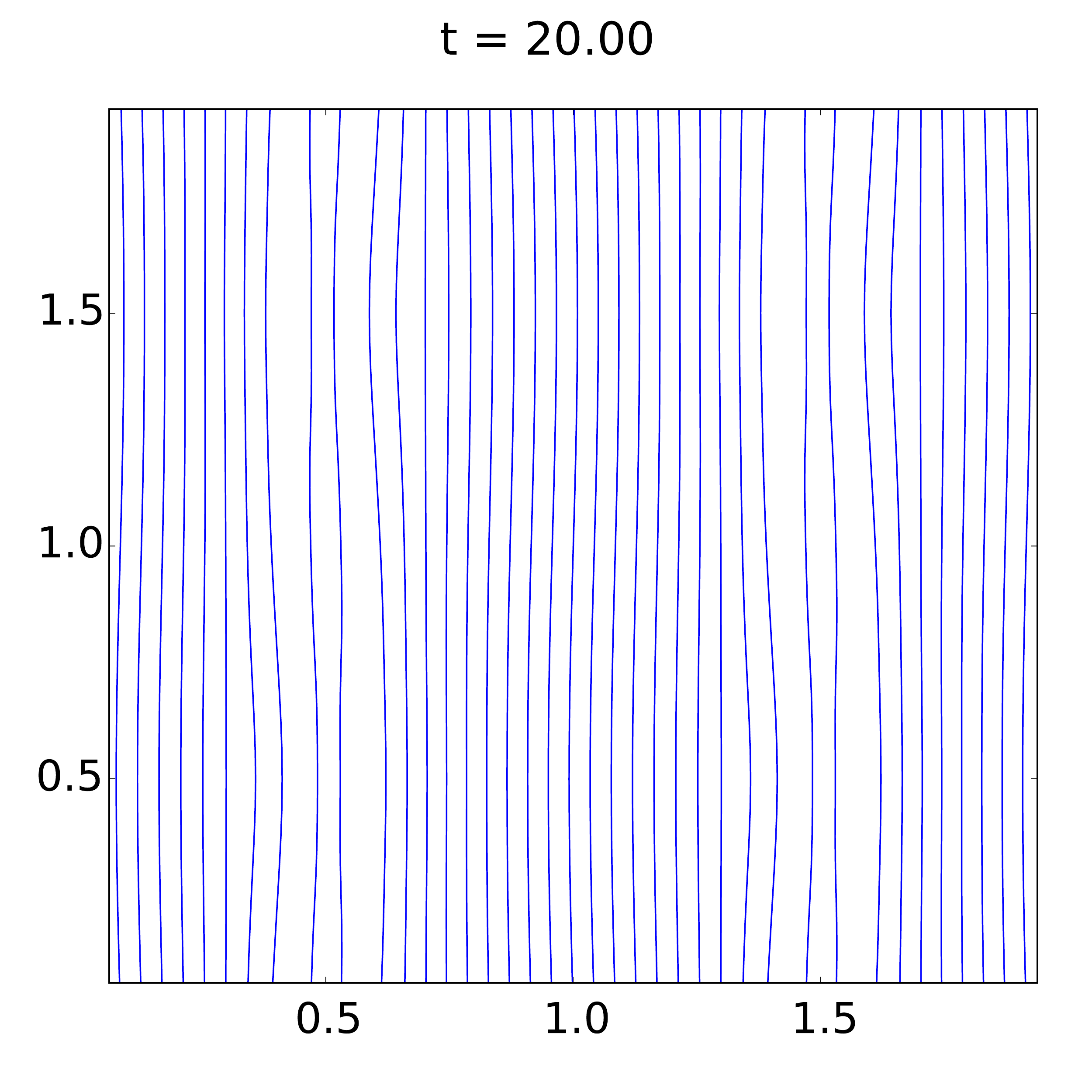}
}
\subfloat{
\includegraphics[width=.32\textwidth]{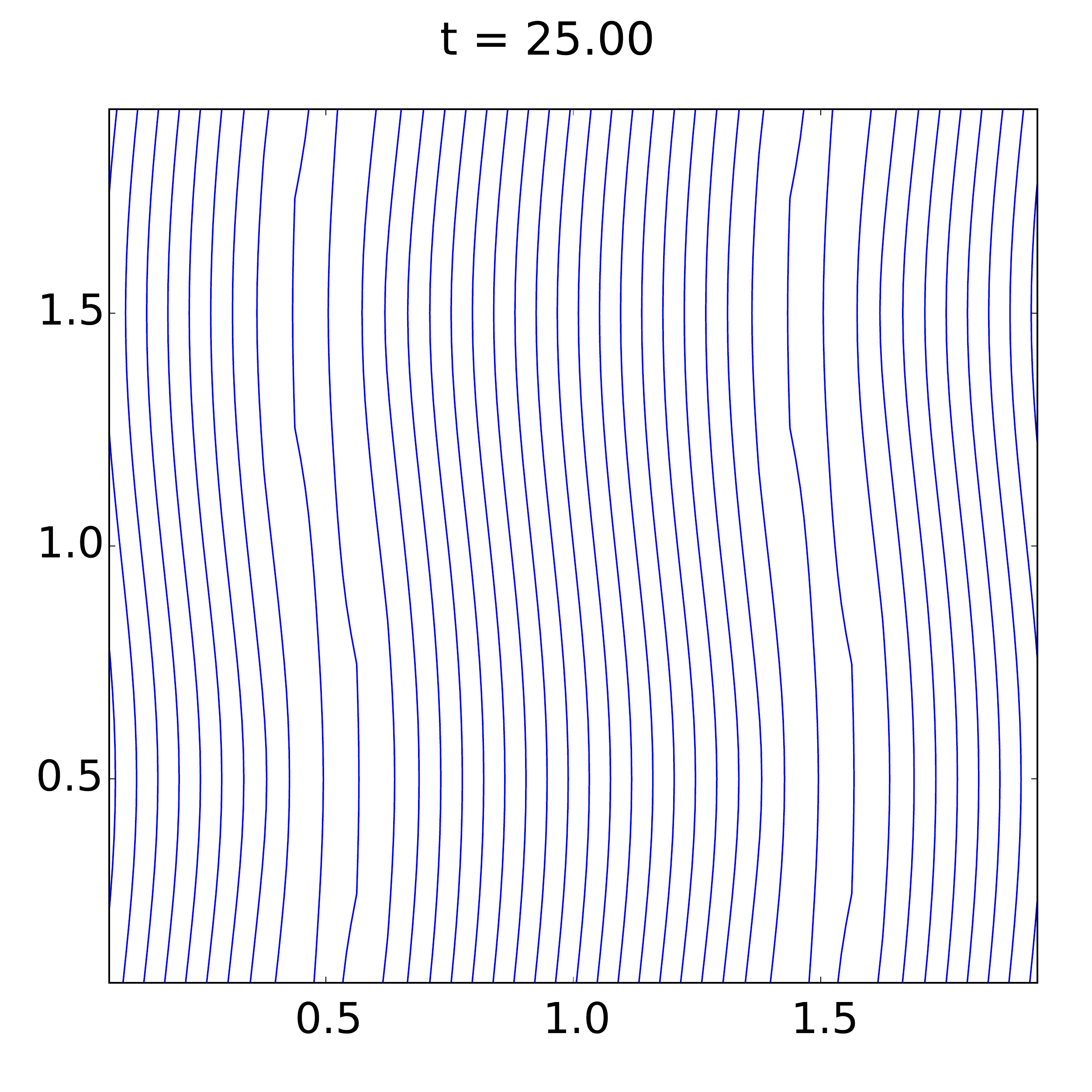}
}

\subfloat{
\includegraphics[width=.32\textwidth]{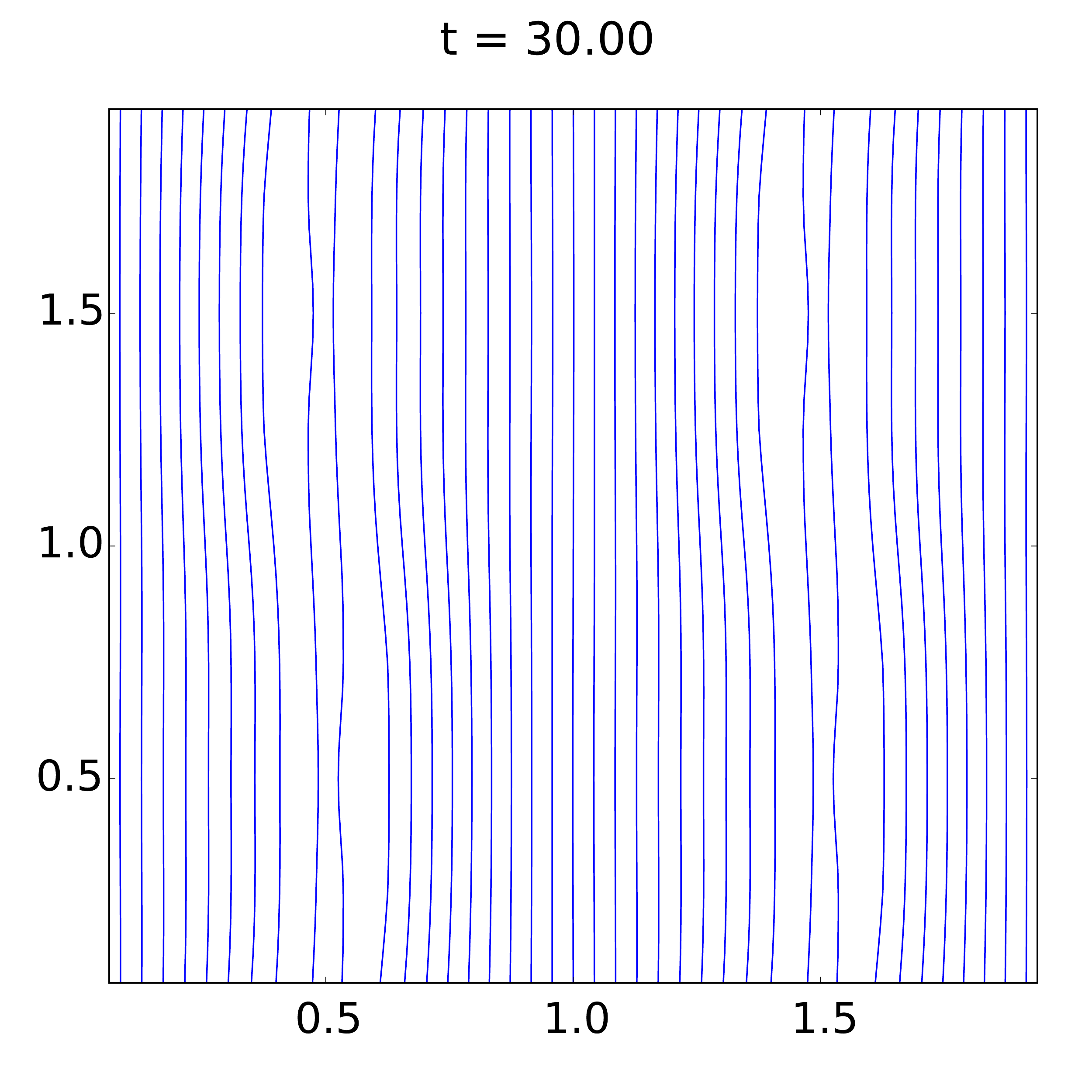}
}
\subfloat{
\includegraphics[width=.32\textwidth]{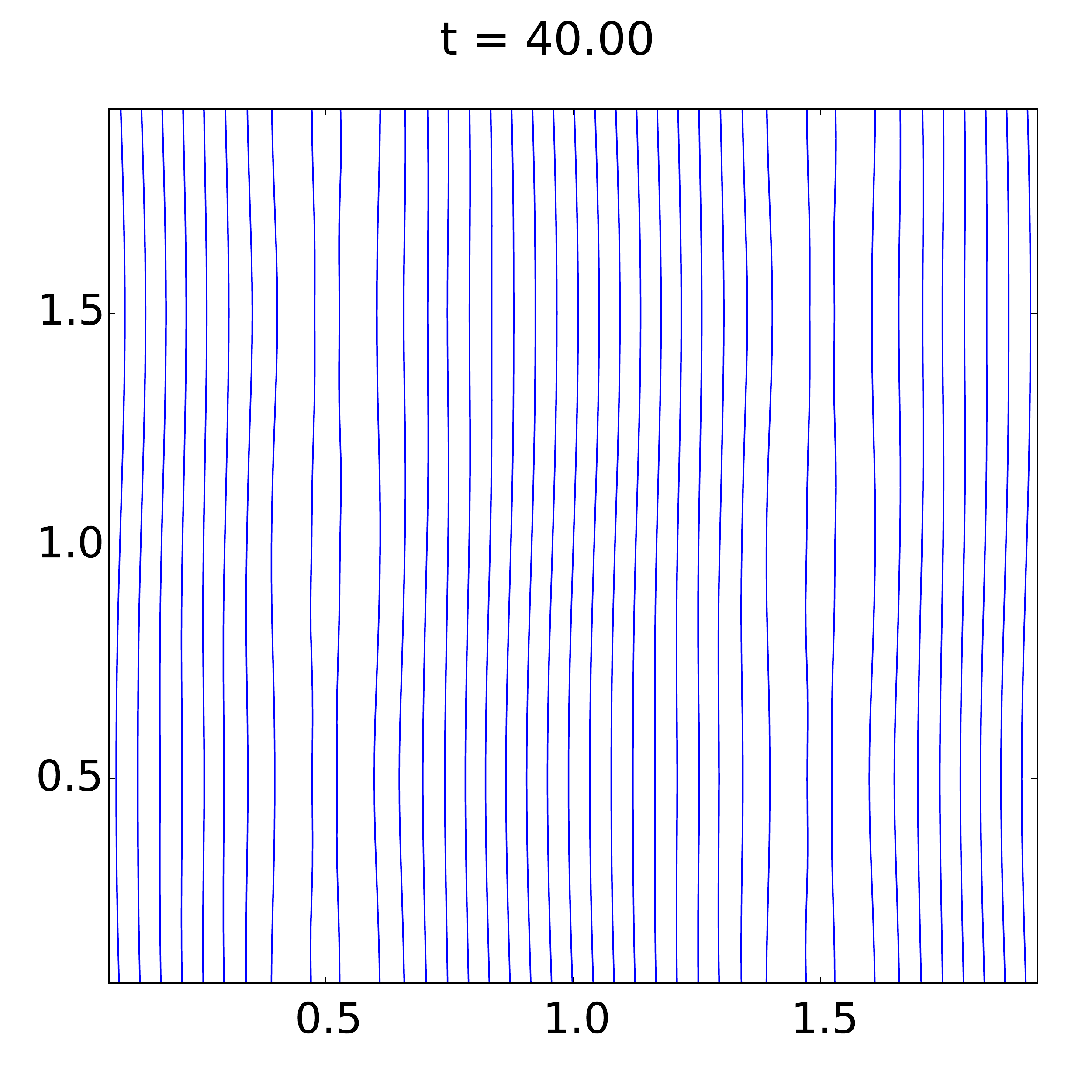}
}
\subfloat{
\includegraphics[width=.32\textwidth]{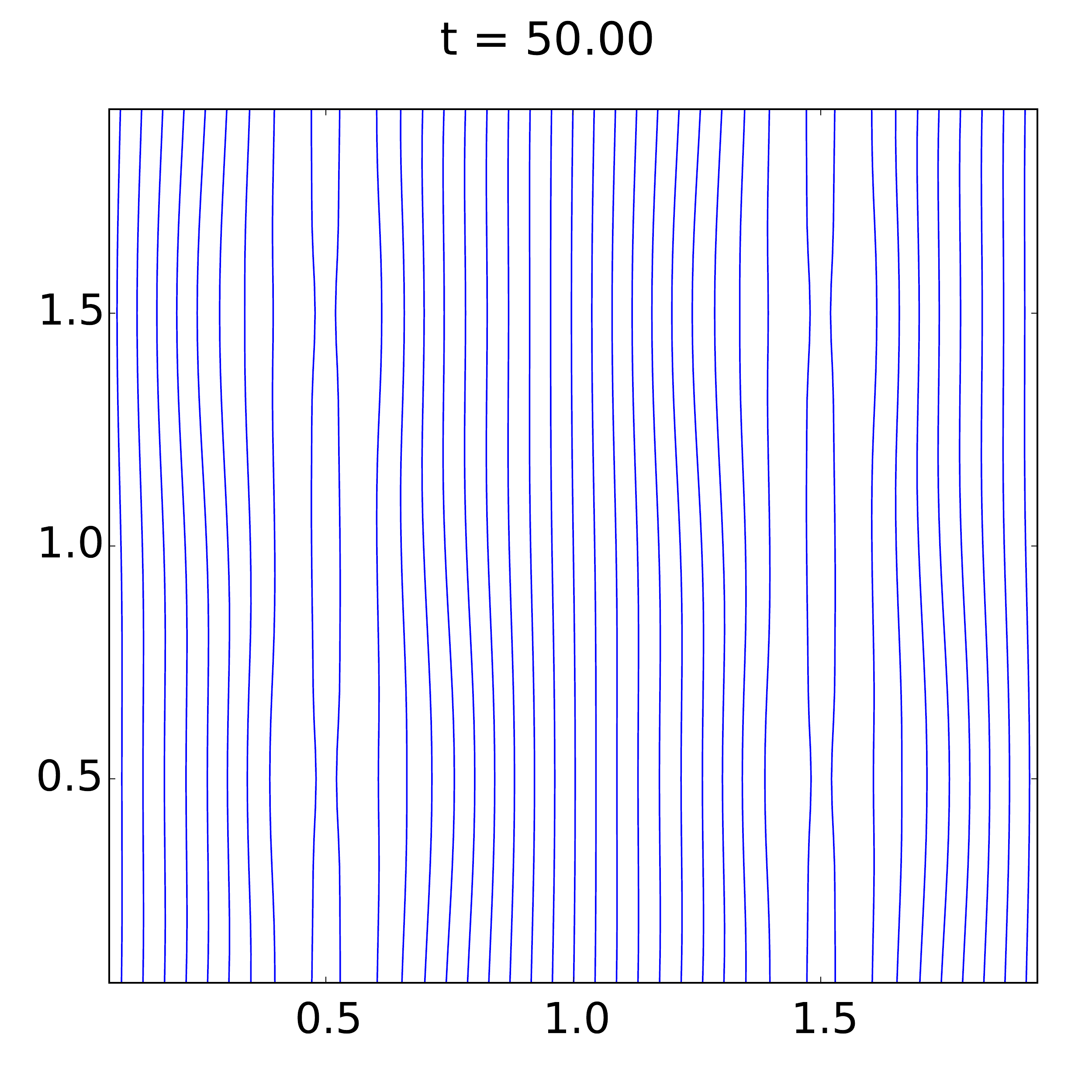}
}

\subfloat{
\includegraphics[width=.32\textwidth]{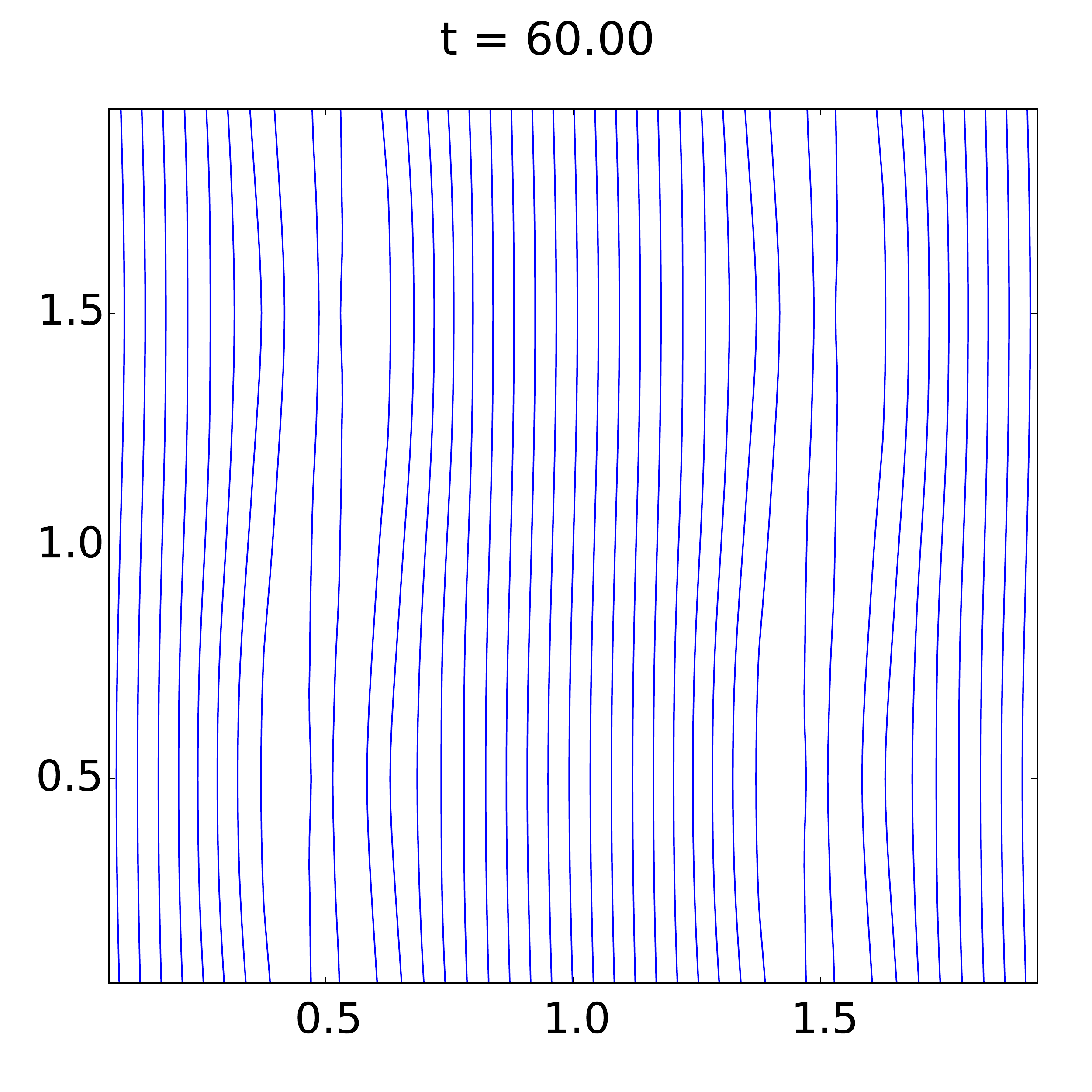}
}
\subfloat{
\includegraphics[width=.32\textwidth]{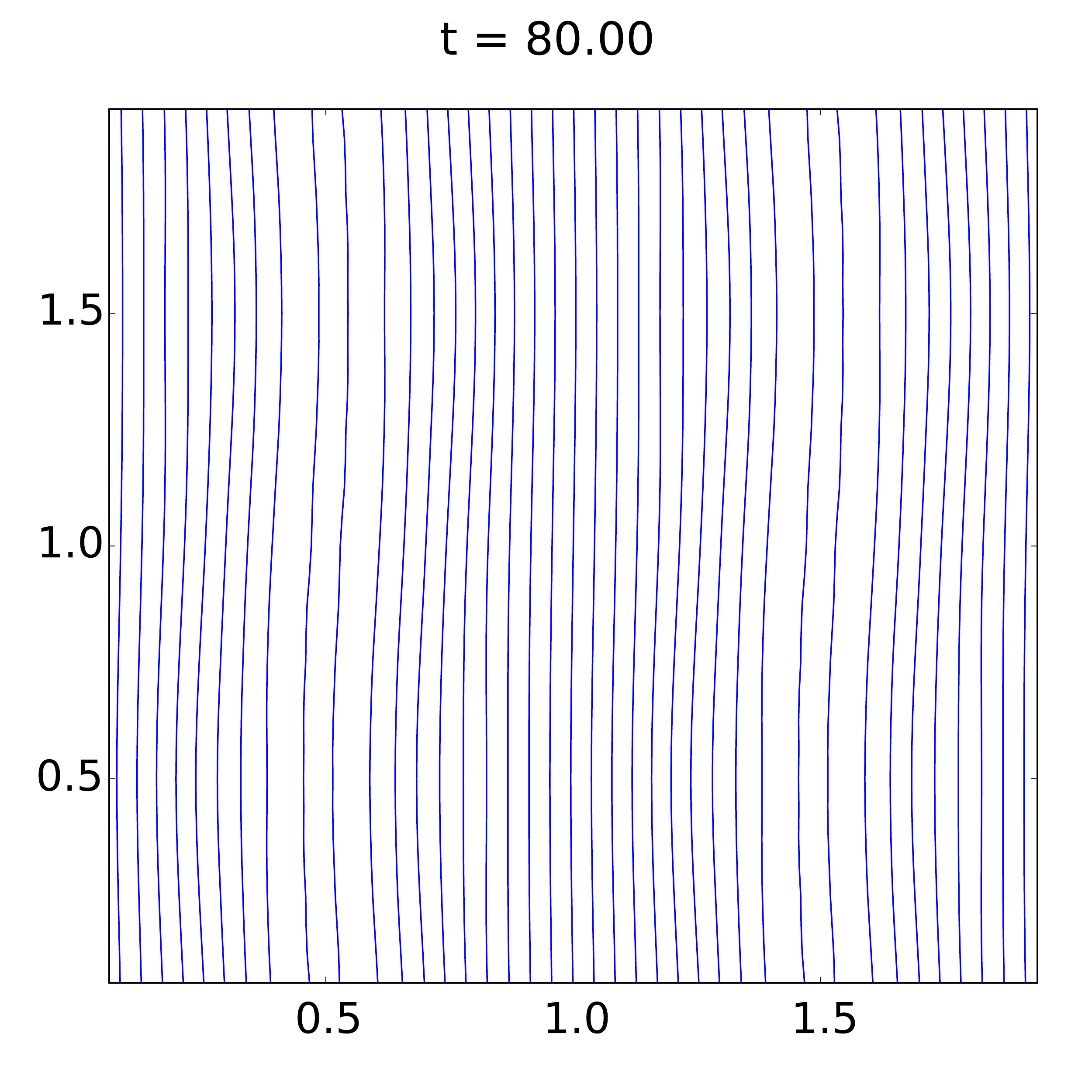}
}
\subfloat{
\includegraphics[width=.32\textwidth]{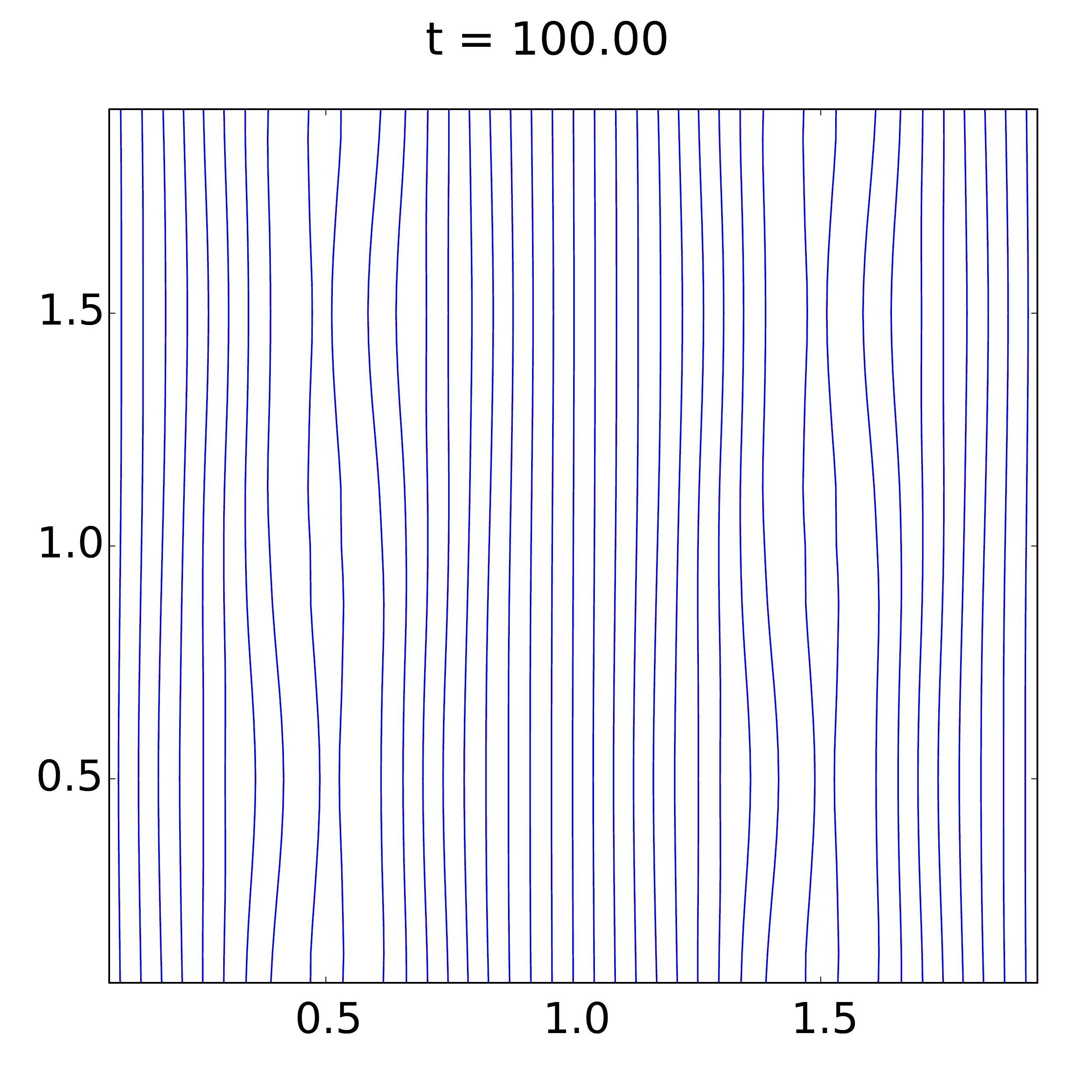}
}

\caption{Smooth current sheet. Magnetic field lines.}
\label{fig:current_sheet_tanh_field_lines}
\end{figure}

To investigate the preservation of the magnetic field line topology on longer time scales, we consider therefore also a less severe current sheet example, following a $\tanh$ profile as it is used in reconnection studies~\cite{GrassoCalifano:2001}. In the $\tanh$ case, the magnetic field changes sign not suddenly but smoothly.
Under this condition, we can run the simulation much longer. Figure~\ref{fig:current_sheet_tanh_field_lines} shows the field line evolution for the case of the smooth magnetic field up to $t = 100$. We observe that the field lines bend but do not reconnect, as is expected from the theory but rarely observed in numerical simulations, especially on the time scales we are considering here.
Most numerical schemes do feature a certain amount of numerical resistivity, leading to unphysical reconnection. In the variational integrator, such spurious resistivity appears to be completely absent, at least in the case of a continuous magnetic field.

\section{Summary}

We have shown how the formal Lagrangian approach to variational integrators described in~\cite{KrausMaj:2015} can be used to derive geometric integration schemes for systems as complicated as magnetohydrodynamics featuring strong nonlinearities.
The ideas of \cite{KrausMaj:2015} have been combined with a staggered grid motivated by discrete differential forms.
While the variational approach guarantees exact conservation of energy, magnetic helicity and cross helicity, the discrete exterior calculus on the staggered grid ensures conservation of the divergence of the magnetic field.
The spatial discretisation thus obtained has already been described in~\cite{LiuWang:2001} and~\cite{Gawlik:2011}. Here, however, it is combined with a symmetric time integrator, exhibiting better energy conservation than previous methods.
The excellent conservation properties with respect to energy, magnetic helicity and cross helicity have been demonstrated in various numerical examples drawn from the literature.
Particularly remarkable is the absence of artificial magnetic reconnection provided that the magnetic field is sufficiently regular.

A limitation of the proposed method is the finite-difference staggered grid approach which is not easily generalised to higher-order methods.
This work, however, should rather be understood as a proof-of-principle of the applicability of the formal Lagrangian approach of~\cite{KrausMaj:2015} to magnetohydrodynamics and that exact conservation properties can be achieved even with very low-order discretisations. In contrast to the Euler-Poincar\'{e} approach of~\cite{Gawlik:2011}, formal Lagrangians can easily be discretised using finite element exterior calculus~\cite{Arnold:2006, Arnold:2010, Christiansen:2011}, mimetic spectral elements~\cite{Gerritsma:2012, Kreeft:2011, Palha:2014} or spline differential forms~\cite{Buffa:2010, Buffa:2011, Ratnani:2012}. This allows for the derivation of numerical schemes of arbitrary order and on general meshes in a straight-forward way (see~\cite{Kraus:2016:Evolution} and~\cite{Kraus:2016:Isogeometric} for developments in this direction).
Moreover, our approach is also applicable to extended magnetohydrodynamics models like inertial MHD~\cite{Kraus:2018:InertialMHD}.

\section*{Acknowledgements}

The first author would like to thank Jonathan Squire and Yao Zhou for helpful discussions, as well as Akihiro Ishizawa for suggesting the application of the formal Lagrangian variational integrator method of~\cite{KrausMaj:2015} to magnetohydrodynamics in order to obtain geometric integrators which preserve important invariants of the system.
Further we thank Yaman G\"{u}{\c c}l\"{u} for valuable comments and reading a draft of the paper.
The first author has received funding from the European Union's Horizon 2020 research and innovation programme under the Marie Sklodowska-Curie grant agreement No 708124. The views and opinions expressed herein do not necessarily reflect those of the European Commission.

\phantomsection
\addcontentsline{toc}{section}{References}
\bibliographystyle{plainnat}
\bibliography{vi_ideal_mhd2d}

\end{document}

%% file: definitions_math.tex
\makeatletter
\@ifundefined{chapter}
   {\newtheorem{theorem}{Theorem}[section]}
   {}
\makeatother

\newtheorem*{exercise*}{Exercise}

\AtBeginDocument{
    
    \let\obar\bar

    \let\bar\undefined

    \newcommand{\bar}[1]{\ensuremath{\overline{#1}}}						
}

\def\ba{\begin{align}}
\def\ea{\end{align}}

\newcommand{\mbb}[1]{\ensuremath{\mathbb{#1}}}							
\newcommand{\mcal}[1]{\ensuremath{\mathcal{#1}}}							
\newcommand{\mrm}[1]{\ensuremath{\mathrm{#1}}}							

\newcommand{\rsp}{\ensuremath{\mathbb{R}}}					
\newcommand{\zsp}{\ensuremath{\mathbb{Z}}}					

\newcommand{\lie}{\ensuremath{\pounds}}									
\DeclareMathOperator{\hodge}{\star}										
\newcommand{\ext}{\ensuremath{\mathsf{d}}}								

\newcommand{\phy}{\ensuremath{\varphi}}									

\newcommand{\iprod}{\ensuremath{\boldsymbol{\imath}}}					

\newcommand{\pair}[2]{\ensuremath{\left < #1 \, , \, #2 \right >}}
\newcommand{\bracket}[1]{\ensuremath{\left < #1 \right >}}				

